\documentclass[12pt]{amsart}
\usepackage[foot]{amsaddr}
\usepackage[all]{xy}
\usepackage{hyperref}
\usepackage{float}
\usepackage{graphicx}
\usepackage[dvipsnames]{xcolor}
\usepackage{caption}
\usepackage[title]{appendix}
\usepackage{amssymb}
\usepackage{tikz}
\usetikzlibrary{decorations.pathmorphing}
\usepackage{hhline}
\usetikzlibrary{arrows}
\definecolor{singcol}{RGB}{185,45,45}
\definecolor{tgtcol}{RGB}{40,95,170}

\usepackage[latin1]{inputenc}
\usepackage{subcaption}
\usepackage[T1]{fontenc}
\usepackage{multirow}
\usepackage{changepage}
\usepackage{geometry}
\usepackage{tabularx}
\usepackage{comment}
\usepackage{makecell}
\usepackage{longtable}
\usepackage{etoolbox}
\usepackage{dynkin-diagrams}
\usepackage{array}
\usepackage{booktabs}
\usepackage{ragged2e}
\def\row#1/#2!{#1_{\IfStrEq{#2}{}{n}{#2}} & \dynkin{#1}{#2}\\}

\renewcommand{\qed}{\hfill \ensuremath{\Box}}
\reversemarginpar

\usepackage{amsmath}
\usepackage{amssymb}
\usepackage{latexsym}
\usepackage{amsfonts}
\usepackage{epsfig}
\usepackage{picinpar}
\usepackage[osf,sc]{mathpazo}
\usepackage{microtype}
\usepackage{textcomp}
\mathtoolsset{showonlyrefs=true}

\newcommand{\triplehline}{%
	\hline%
	\noalign{\vskip 1.2pt}%
	\hline%
	\noalign{\vskip 1.2pt}%
	\hline%
}

\newcommand{\smat}[1]{%
	\left(
	\begin{smallmatrix}
		#1
	\end{smallmatrix}
	\right)
}

\newcommand{\A}{{\mathbb A}}

\newcommand{\C}{{\mathbb C}}

\newcommand{\F}{{\mathbb F}}

\renewcommand{\H}{{\mathbb H}}
\newcommand{\I}{{\mathbb I}}

\newcommand{\N}{{\mathbb N}}
\renewcommand{\O}{{\mathbb O}}
\renewcommand{\P}{{\mathbb P}}
\newcommand{\Q}{{\mathbb Q}}
\newcommand{\R}{{\mathbb R}}

\newcommand{\T}{{\mathbb T}}

\newcommand{\Z}{{\mathbb Z}}
\newcommand{\lra}{\longrightarrow}

\begin{document}
\title{The platonic elliptic surfaces}
\author{Nutsa Gegelia and Duco van Straten}

\email{gegelian@uni-mainz.de, straten@mathematik.uni-mainz.de}

\begin{abstract}
  We construct certain special rational elliptic surfaces with four reduced singular fibres that are
  naturally attached to the platonic solids and belong to a group of surfaces complementing the
  well-known semi-stable Beauville surfaces. We describe the basic algebraic geometric properties of
  these surfaces, their associated Picard-Fuchs operators, integer sequences, Laurent polynomial representations and Ap\'ery constants. On the arithmetic side,
  we discover the need for a refinement of the Dwork-expansion used to find the Euler factors of the fibres of these fibrations.
This is explained by comparing the $q$-coordinates coming from the elliptic curve and the $q$-coordinate of the Picard-Fuchs equation and is also directly reflected in the shape of the monodromy matrices in the Frobenius basis.

\end{abstract}

\maketitle

\section*{Introduction}
Starting from the three reflection groups of the platonic solids we study certain elliptic surfaces that
are naturally attached to them. These surfaces have four reduced singular fibres, three of which
are semi-stable. It turns out that these  have many properties in common with the well-known Beauville
cases of four semi-simple fibres, but also exhibit some interesting new features. We determine the Picard-Fuchs operators
for these surfaces, which are Heun equations of a special type, from which one obtains a $2$-term recursion relation for
the expansion coefficients, forming integer sequences analogous to the sequence of (small) Ap\'ery numbers. We find that these
surfaces and operators exhibit some unusual arithmetic properties, which is reflected in the monodromy of the operator
and the Dwork-expansion that determines the Euler factors of the fibres of these elliptic surfaces.

This paper, which is partly expository in nature,  is structured as follows.

In Section \ref{ch:platonicsolidsandreflections} we describe how the reflection groups of the platonic solids
give rise to an elliptic threefold, from which we obtain elliptic surfaces as sections, which we call {\em platonic elliptic
  surfaces}.

In Section \ref{ch:ellipticsurfaces} we review some basic properties of elliptic surfaces in general and describe how these platonic elliptic surfaces fit in the classification and point out the differences with some other famous elliptic surfaces that have a clear relation to the platonic solids.

In Section \ref{ch:normalisedperiodintegrals} we describe the Picard-Fuchs operators attached to (normalised) period integrals of our surfaces, which turn out to be special Heun operators with exponents $(0,0)$ at the three finite singular
points and of finite order at $\infty$.

In Section \ref{ch:platonicsequences} we discuss the integer sequences attached to these operators and present Laurent-polynomial expressions for all these sequences, which imply that Lucas-type congruences
hold for them.

In Section \ref{ch:somearithmeticalaspects} we overview some facts about the arithmetic of the elliptic curves appearing as fibres in our elliptic surfaces, using the method that goes back to Dwork. By comparing, at a nodal fibre, the Serre-Tate canonical $q$-parameter with the $q$-coordinate that emerges from the {\em critical integral scaling} of the
Picard-Fuchs operator one finds the appearance of a logarithmic term in the limit Frobenius matrix, which is unlike what happens for the Beauville operators. We find that the monodromy matrices expressed in the Frobenius basis are affected by the logarithmic term in a corresponding way. This will be explained in Section \ref{ch:monodromycalculations}. We
also determine signatures of the resulting
monodromy groups, which shows that, with three exceptions, our surfaces are in fact {\em modular}, but for a non-congruence subgroup. 
In a final Section \ref{ch:outlook} we look
at possible variations and generalisations of
the constructions considered in this paper.

In Appendix \ref{ap:herfurtnersurfaces} some further data attached to our examples is collected, among this some data belonging to other related elliptic surfaces with four singular fibres.

{\em Acknowledgement:} We are grateful to Bradley Klee for calling our attention to his beautiful work \cite{Klee} in which the
platonic operators and some of the associated integer sequences are connected to the vibrational spectrum of molecules with
platonic symmetry. Our desire to understand the relationship of the platonic solids, the elliptic surfaces  and the special Heun operators resulted in this paper. Thanks also to Wadim Zudilin for moral
support and constructive feedback on an earlier version of this paper.
Further thanks to {\tt Claude Fable} and {\tt ChatGPT}, which was of great help with
verification and some extremely helpful hints in the later stage of this project.

We gratefully acknowledge the support by the Deutsche Forschungsgemeinschaft (DFG, German Research Foundation) -- Project-ID 444845124 -- TRR 326.

\section{Platonic Solids and reflection groups}\label{ch:platonicsolidsandreflections}
\subsection{Three reflection groups}

\begin{figure}[ht]
\includegraphics[height=4cm]{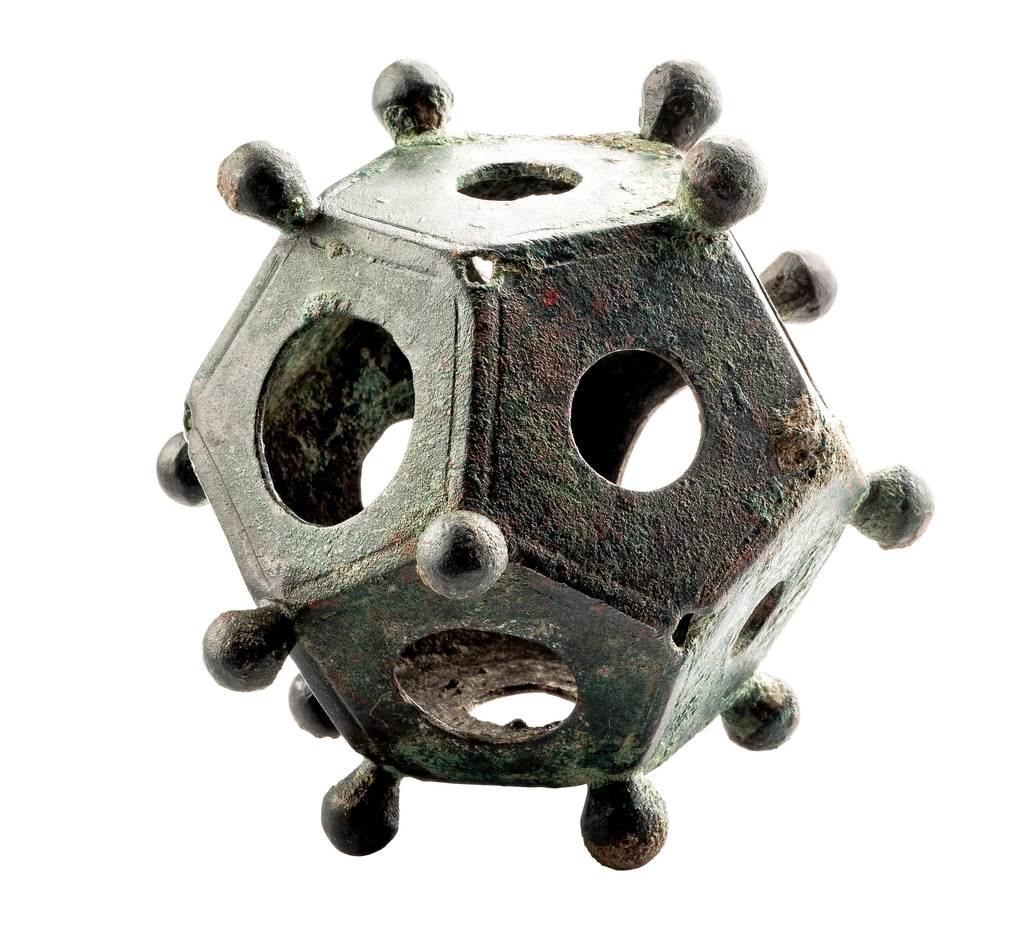}
\caption{Roman dodecahedron.}
\label{fig:romdod}
\end{figure}

The five platonic solids appear ubiquitous in \linebreak mathematics and popular culture. There are good reasons to believe they were known in prehistoric times \cite{Lloyd}, although the details of their discovery remain obscure. In the Scholium to book XIII of the Elements of Euclid, which describes the construction of the platonic solids, the tetrahedron, cube and dodecahedron 
are ascribed to the pythagoreans and the octahedron and icosahedron to Thaetetus. We refer to \cite{Waterhouse} for  a historical acccount. The story of the discovery of the dodecahedron by the Pythagorean Hippasos and his fate is told by Iamblicus \cite{Iamblichus}, Chapter XVIII. The precise use of the so-called {\em Roman dodecahedron} (Figure \ref{fig:romdod}) found all over Europe is as yet unknown.

The basic geometry of the regular polyhedra belongs to the standard canon of elementary mathematics for which we refer to Coxeter \cite{Coxeter}, where among other things, the general theory of real reflection groups (Coxeter groups) and their Coxeter-Dynkin diagrams is developed.

\begin{figure}[ht]
	\centering
	\captionsetup[subfigure]{labelformat=empty}
	
	\begin{subfigure}[t]{0.32\textwidth}
		\centering
		\includegraphics[height=4cm]{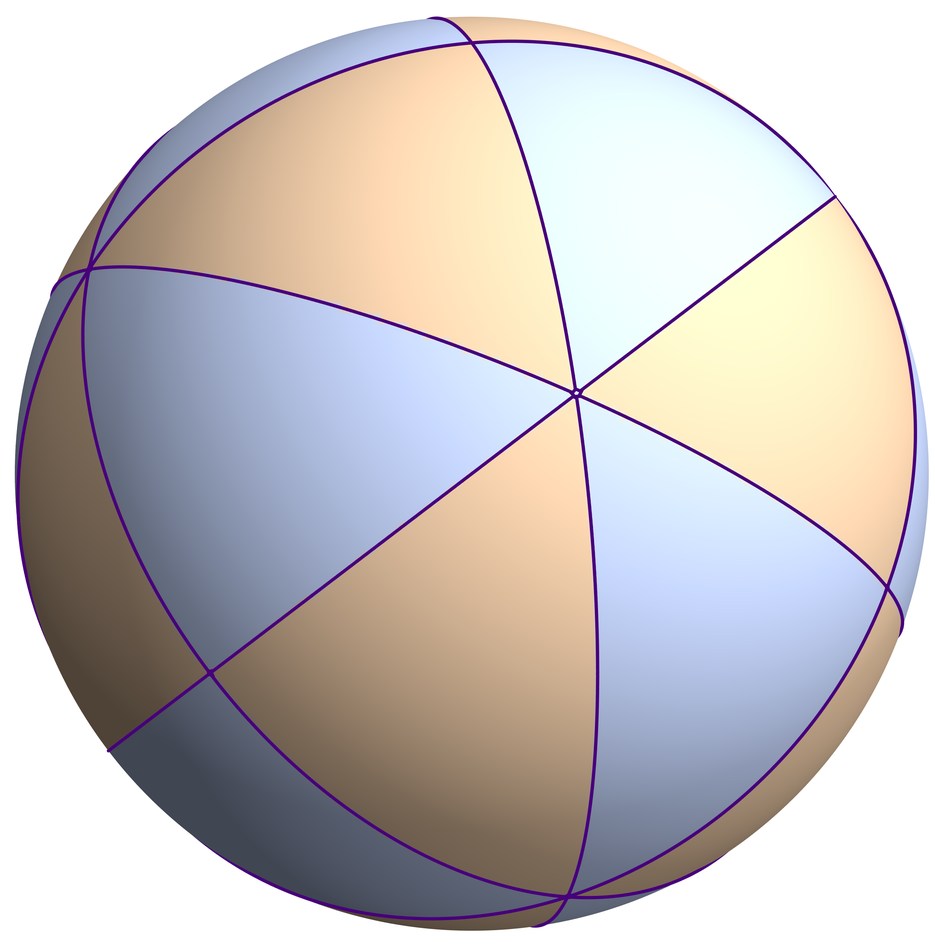}
		\subcaption{\((2,3,3)\)}
		\label{fundamentaldomain1}
	\end{subfigure}
	\hfill
	\begin{subfigure}[t]{0.32\textwidth}
		\centering
		\includegraphics[height=4cm]{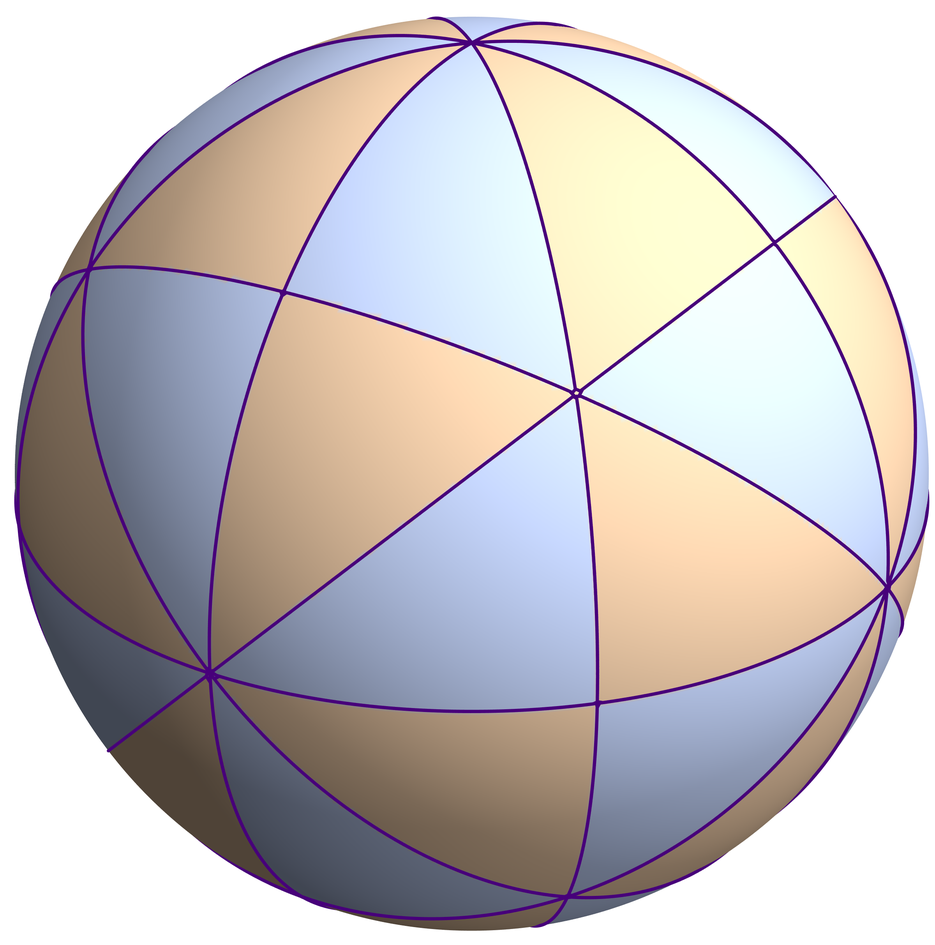}
		\subcaption{\((2,3,4)\)}
		\label{fundamentaldomain2}
	\end{subfigure}
	\hfill
	\begin{subfigure}[t]{0.32\textwidth}
		\centering
		\includegraphics[height=4cm]{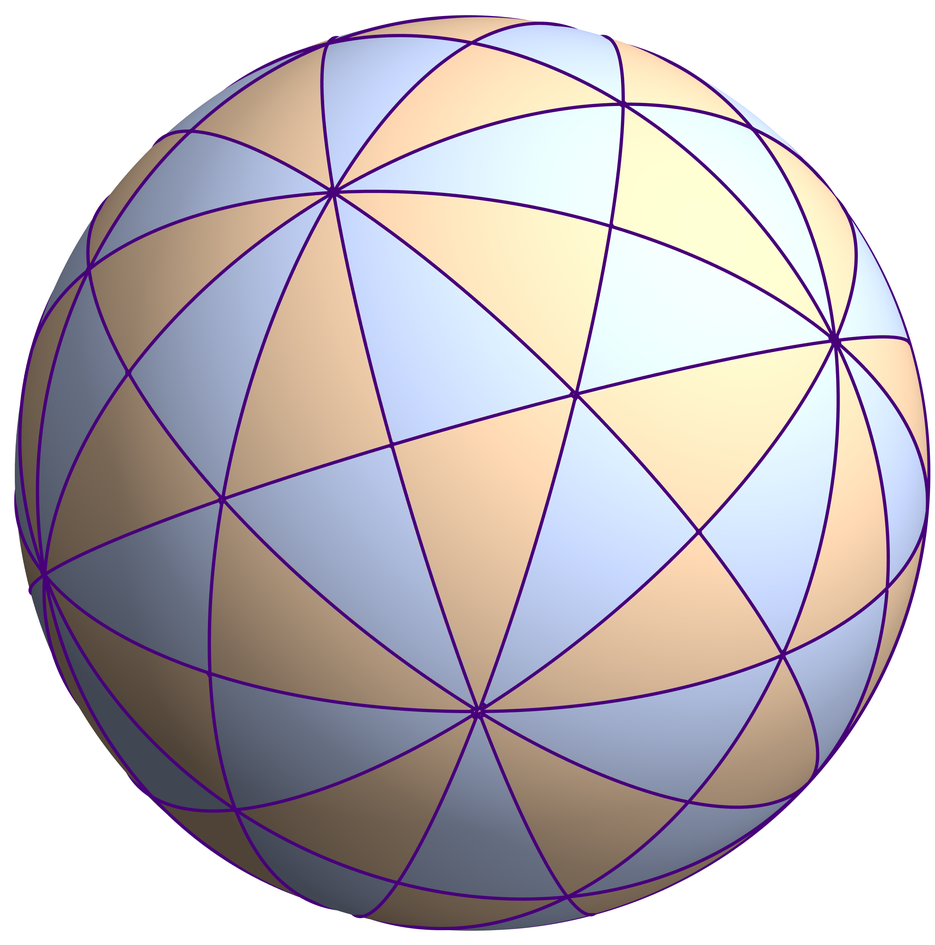}
		\subcaption{\((2,3,5)\)}
		\label{fundamentaldomain3}
	\end{subfigure}
	
	\caption{Spherical tessellations associated with the platonic triples.}
	\label{fig:spherical-tessellations}
\end{figure}

The three-dimensional reflection groups are labelled by {\em platonic triples} $(p,q,r)$, with $\frac{1}{p}+\frac{1}{q}+\frac{1}{r} >1$,
which encode the angles $\pi/p, \pi/q, \pi/r$ of the spherical triangle forming a fundamental domain (Figure \ref{fig:spherical-tessellations}) for the reflection group $G$. The cases of $(p,q,r)=(2,3,3), (2,3,4), (2,3,5)$ lead respectively to the tetrahedron, octahedron/hexahedron, icosahedron/dodecahedron, whereas the cases $(2,2,n)$ correspond to the dihedral groups.

The number of elements of the reflection group $G$ attached to $(p,q,r)$ is given by the uniform formula
\[|G| = \frac{4}{1/p+1/q+1/r-1} .\]

This is also the number of elements of the associated {\em binary groups} ${\bf G}$, inverse image of the rotation group $SG:=SO(3) \cap G$ under the 2-1 covering $SU(2) \lra SO(3)$ that were studied by F. Klein in \cite{Klein}. As recognised by DuVal \cite{DuVal}, the corresponding quotient singularities $\C^2/{\bf G}$ are exactly those that can be resolved by a system of curves organised in a Dynkin graph of type ADE. Thus, the reflection groups of the tetrahedron, octahedron and icosahedron are related in this way to the Dynkin diagrams $E_6$, $E_7$ and $E_8$, where the arm-lengths to the central node coincide with the platonic triples (Table \ref{tab:platonic-reflection-groups}). This relation between these two basic {\em trinities} in the sense of Arnol'd \cite{Arnold2} was more recently reconsidered in  \cite{Dechant}. The ADE list also coincides with the list of simple singularities in the sense of Arnol'd, \cite{Arnold1}. The ADE-graphs and their associated higher-dimensional reflection groups describe the root systems of the simply laced Lie-algebras bearing the same names. For more information highlighting the multitude of relations between reflection groups, algebraic geometry and singularity theory, we refer to the survey papers by I. Dolgachev \cite{Dolgachev}, A. Givental \cite{Givental} and O. Shcherbak \cite{Shcherbak}.

In this paper we will be concerned here with the three reflection groups $G=\T,\O,\I$ of the tetrahedron, octahedron and icosahedron respectively and leave the dihedral cases corresponding to $(2,2,n)$ aside.

\begin{table}[ht]
	\centering
	\renewcommand{\arraystretch}{1.5}
	\begin{tabular}{c|ccccc}
		\hline
		\textbf{Polyhedron}
		& \textbf{Triangle type}
		& \textbf{Group}
		& \textbf{Order}
		& \textbf{Coxeter type}
		& \textbf{ADE type} \\
		\hline
		Tetrahedron
		& $(2,3,3)$
		& $\T$
		& $24$
		& $A_3\;\dynkin[Coxeter]{A}{3}$
		& $E_6\;\dynkin{E}{6}$ \\
		\hline
		Octahedron
		& $(2,3,4)$
		& $\O$
		& $48$
		& $B_3\;\dynkin[Coxeter]{B}{3}$
		& $E_7\;\dynkin{E}{7}$ \\
		\hline
		Icosahedron
		& $(2,3,5)$
		& $\I$
		& $120$
		& $H_3\;\dynkin[Coxeter]{H}{3}$
		& $E_8\;\dynkin{E}{8}$ \\
		\hline
	\end{tabular}
	\caption{Basic data for the Platonic reflection groups: spherical triangle triples, group orders, Coxeter types, and corresponding ADE types.}
	\label{tab:platonic-reflection-groups}
\end{table}
As abstract groups, the rotation subgroups $S{\T}, S{\O}, S{\I}$ are isomorphic
to the alternating group $\mathcal{A}_4$, the symmetric group $\mathcal{S}_4$
and the alternating group $\mathcal{A}_5$.
\subsection{The Goursat Invariants}
In his memoir \cite{Goursat}, Goursat  gives an account of the {\em invariant theory} of the rotation and reflection groups of the regular polyhedra and studies in some detail the beautiful surfaces in three space invariant under these groups. W. Barth's celebrated $65$ nodal sextic \cite{Barth} is a special member of a family of sextics described in Goursat's paper.
This work is complementary to the earlier work of F. Klein \cite{Klein}, which is mainly concerned with the invariants of
the associated binary groups and the resulting surface singularities now commonly named after him.

{\bf Proposition (Goursat):} {\em The ring of invariants of the reflection groups $G=\T,\O,\I$ in $\R^3$ is freely generated by three polynomials  $P,Q,R$:}
\[ \Q[x,y,z]^G=\Q[P,Q,R].\]

In fact Goursat constructs a very specific set of generators $P$, $Q$ and $R$. In all cases, the invariant of lowest degree is
$P=x^2+y^2+z^2$ whereas the invariants $Q$ and $R$ are made up by products of linear forms defining symmetry planes permuted by $G$. He observes that the intersection of the level sets $P=\alpha$, $Q=\beta$, $R=\gamma$ are the orbits of $G$ and consist of
$|G|=deg(P)\cdot deg(Q)\cdot deg(R)$ points, from which he concludes that any
symmetric polynomial is uniquely expressed polynomially in terms of $P,Q,R$. Nowadays that
result is seen as a special case of the theorem of Shephard and Todd \cite{ShephardTodd},
which states that the same holds for any complex reflection group.

There is also an additional {\em semi-invariant} $S$, a product of an {\em odd} number of symmetry planes permuted by the group $G$,
so that it is multiplied by $-1$ under a reflection in $G$, and hence is only an invariant of the associated rotation group $SG$.

For example, for the icosahedron, the sextic $Q$ is the product of the six planes parallel to the opposite pairs of faces of the dual dodecahedron and the dectic $R$ is the product of the ten planes parallel to the pair of faces of the icosahedron. The semi-invariant $S$ has degree $15$ and is obtained as product of the $15$ reflection planes determined by the 15 pairs of opposite edges of the icosahedron. Table \ref{tab:platonic-polynomial-degrees} summarises these facts.

\begin{table}[ht]
	\centering
	\renewcommand{\arraystretch}{1.5}
	\begin{tabular}{c|c|cccc}
		\hline
		\textbf{Group}
		& \textbf{Order}
		& $\boldsymbol{\deg P}$
		& $\boldsymbol{\deg Q}$
		& $\boldsymbol{\deg R}$
		& $\boldsymbol{\deg S}$ \\
		\hline
		$\T$
		& $24$
		& $2$
		& $3$
		& $4$
		& $6$ \\
		\hline
		$\O$
		& $48$
		& $2$
		& $4$
		& $6$
		& $9$ \\
		\hline
		$\I$
		& $120$
		& $2$
		& $6$
		& $10$
		& $15$ \\
		\hline
	\end{tabular}
	\caption{Degrees of the polynomials $P$, $Q$, $R$, and $S$ for the Platonic reflection groups.}
	\label{tab:platonic-polynomial-degrees}
\end{table}

We refer to Appendix \ref{ap:goursatinvariants} for the explicit form of the Goursat invariants $Q,R,S$ that were used in this paper.

\subsection{Level curves on the sphere}
\begin{figure}[ht]
	\centering
	\captionsetup[subfigure]{labelformat=empty}
	
	\begin{subfigure}[t]{0.32\textwidth}
		\centering
		\includegraphics[height=4cm]{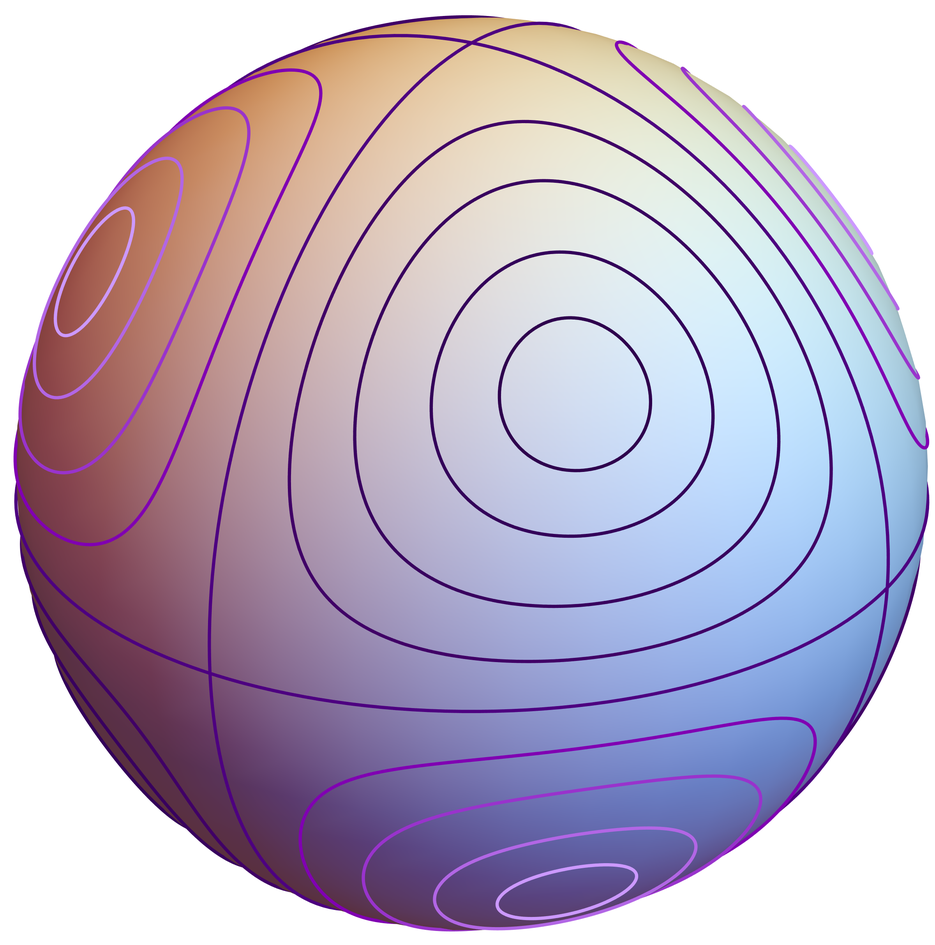}
		\subcaption{Curves of genus \(4\)}
	\end{subfigure}
	\hfill
	\begin{subfigure}[t]{0.32\textwidth}
		\centering
		\includegraphics[height=4cm]{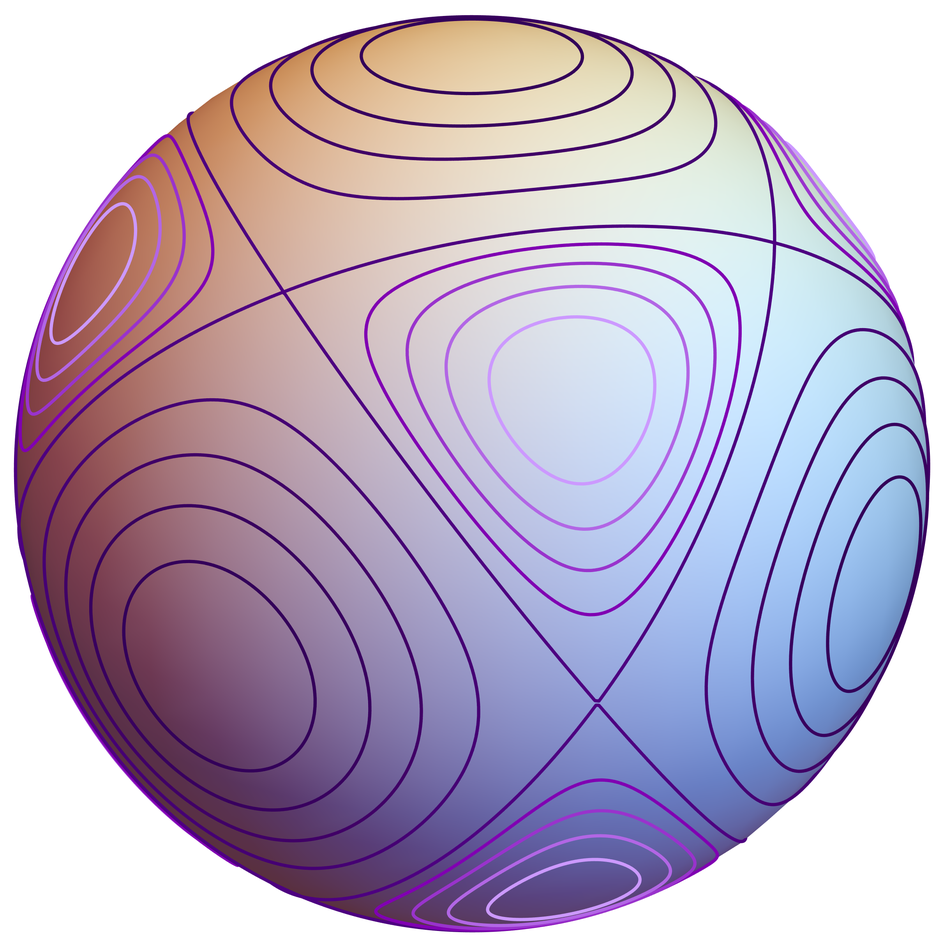}
		\subcaption{Curves of genus \(9\)}
	\end{subfigure}
	\hfill
	\begin{subfigure}[t]{0.32\textwidth}
		\centering
		\includegraphics[height=4cm]{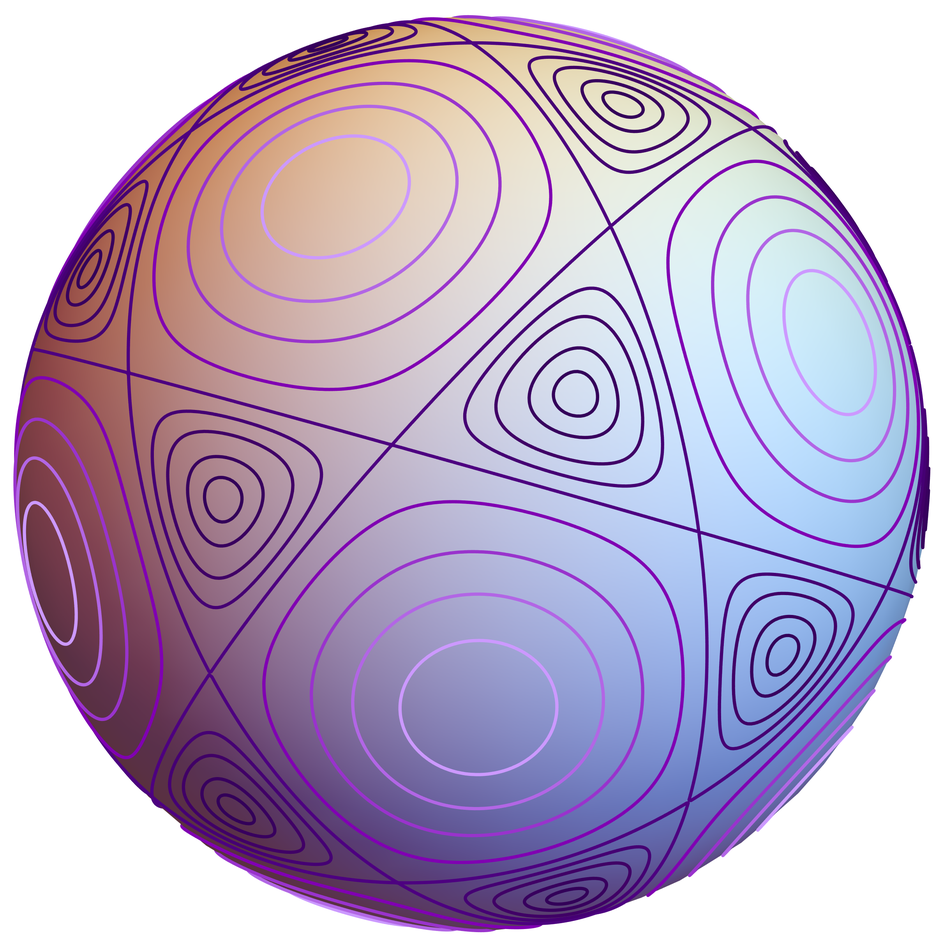}
		\subcaption{Curves of genus \(25\)}
	\end{subfigure}
	
	\caption{Level curves on the sphere $P=r$, $r>0$, cut out by the invariant $Q$.}
	\label{fig:platonic-curves}
\end{figure}

The surface determined by the equation $P=r$, $r>0$ an arbitrary fixed value, is the unit sphere in $\R^3$. The level sets of the
second invariant $Q=t$, $t \in \R$ cut out a family of
algebraic curves
\[   C_{r,t}(\R)^{\circ}:= \{(x,y,z)\in \R^3\;|\;P=r,Q=t\},\]

which belong to the simplest objects that exhibit the symmetry of
the corresponding regular polyhedron (Figure \ref{fig:platonic-curves}).

In the sequel it will be convenient to complexify all spaces involved and
consider the curve $C_{r,t} \subset \P^3$, the projective closure of the complex points $C_{r,t}(\C)^{\circ}$
of the above curve. In this setting we allow values of $r\neq0$.

{\bf Proposition:} {\em  For a general value of $t$, the curve $C_{r,t}$ is an algebraic
curve of genus $4$, $9$ and $25$ for tetrahedron, octahedron and icosahedron
respectively.

There are three values of $t$ for which the curves $C_{r,t}$ are singular, corresponding to the (real) minima,
saddles and maxima of the polynomial $Q$ restricted to the sphere $P=r>0$.}

{\bf proof:} From the adjunction formula it follows that the intersection of
a quadric with a surface of degree $n$ is a curve of genus $(n-1)^2$. For
$\T,\O,\I$ the  degree of $Q$ is $3,4,6$ respectively.

Another way of seeing this is by computing the Euler-characteristic of $C_0$ from the intersection of the sphere with $Q=0$, which consist respectively of $3,4,6$ circles which, projectively and over $\C$, are $2$-spheres. As two circles intersect in two antipodal points, we have $6=2\cdot 3, 12=3\cdot 4, 30=5\cdot 6$ intersection points respectively. Hence the reducible curves consist respectively of $3,4,6$ copies of $\P^1$, with these  intersection points as singular points. For the Euler characteristic of the smoothing $C_{r,t}$ of $C_{r,0}$ we find
$2\cdot 3-2\cdot 6=-6, 2\cdot 4-2\cdot 12=-16, 2\cdot 6-2\cdot 30=-48$, respectively, which gives $g=4,9,25$ by equating
this to $2-2g$.  \qed     

The spherical harmonics associated to the these reflection groups are studied in applied sciences (chemistry and physics in particular) and pure mathematics alike. The papers \cite{Laporte}, \cite{AtiyahSutcliffe}, \cite{Hitchin}, \cite{Klee} contain further information.

\subsection{Invariants of the rotation groups}

Clearly, the polynomials $P,Q,R$ and $S$ are invariants of the rotation groups $SG$. But the square $S^2$ of the semi-invariant $S$ is invariant under the reflection group, so can be expressed polynomially in terms of $P,Q,R$, i.e. we have
a relation of the form
\[S^2=F(P,Q,R)\]
for some unique polynomial $F$. A direct consequence is the following:

{\bf Proposition:} {\em  The ring of invariants of the rotation subgroup $SG:=G \cap SO(3)$ is a hypersurface ring of the form
\[ \Q[x,y,z]^{SG} = \Q[P,Q,R,S]/(S^2-F(P,Q,R))\]

where the relation is given as follows:

{\bf Tetrahedron:}
\begin{equation}\label{eq:relationtetrahedron}
S^2=F(P,Q,R):=-4R^3+P^2R^2+18PQ^2R-(4P^3+27Q^2)Q^2.
\end{equation}

{\bf Octahedron:}
\begin{equation}\label{eq:relationoctahedron}
S^2=F(P,Q,R):=-27R^3+P(18Q-4P^2)R^2+Q^2(P^2-4Q)R.
\end{equation}

{\bf Icosahedron:}
\begin{equation}\label{eq:relationicosahedron}
\begin{aligned}
S^2=F(P,Q,R):=-R^3+P^2(4P^3-65Q)R^2  +  PQ(720Q^2+200 P^6-795 P^3Q)R\\
+500 P^9Q^2-2275 P^6Q^3+3440 P^3Q^4 -1728Q^5.
\end{aligned}
\end{equation}
}

\begin{figure}[ht]
\centering
\captionsetup[subfigure]{labelformat=empty}

\begin{subfigure}[t]{0.32\textwidth}
	\centering
	\includegraphics[width=\textwidth]{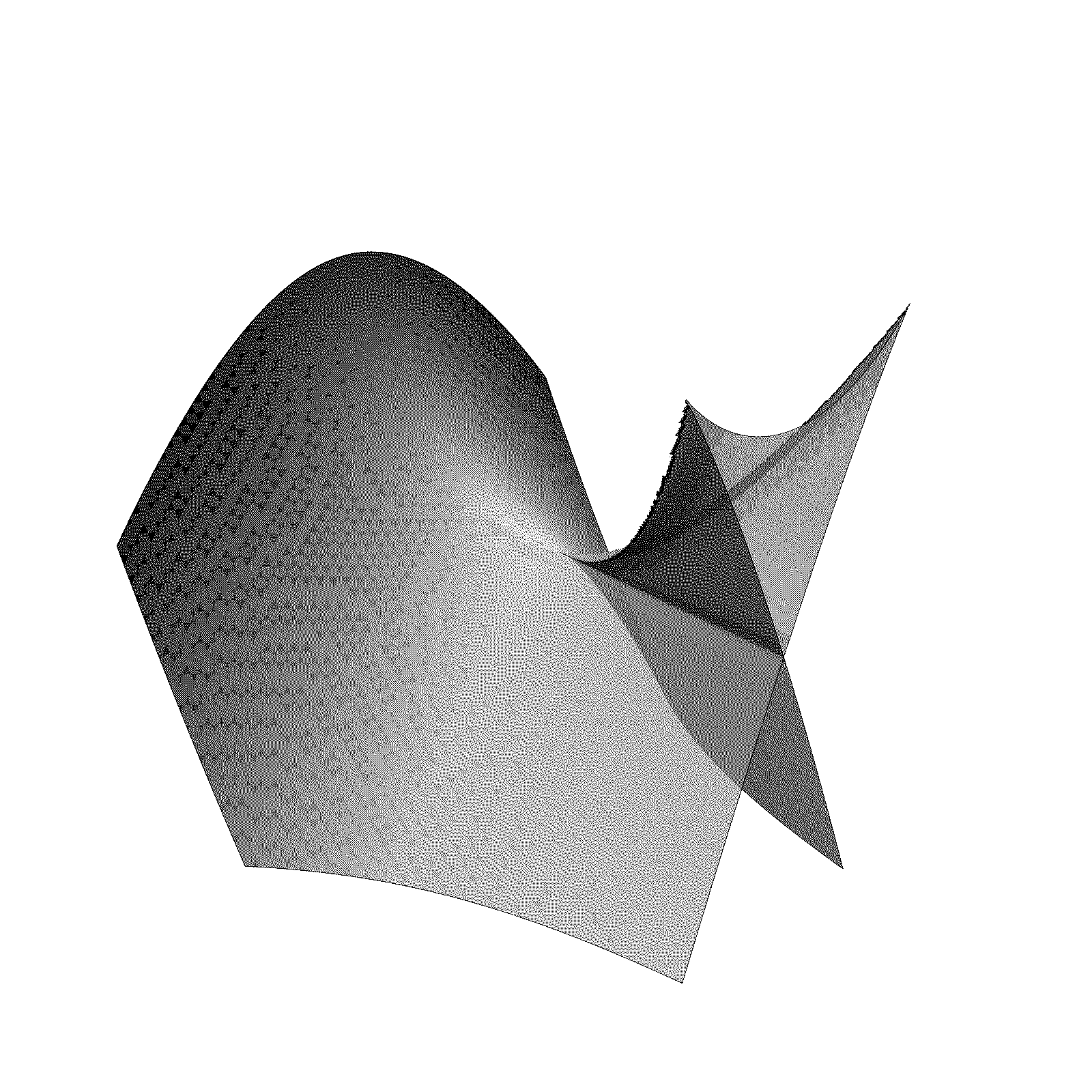}
	\subcaption{Tetrahedral case}
	\label{fig:disctetra}
\end{subfigure}
\hfill
\begin{subfigure}[t]{0.32\textwidth}
	\centering
	\includegraphics[width=\textwidth]{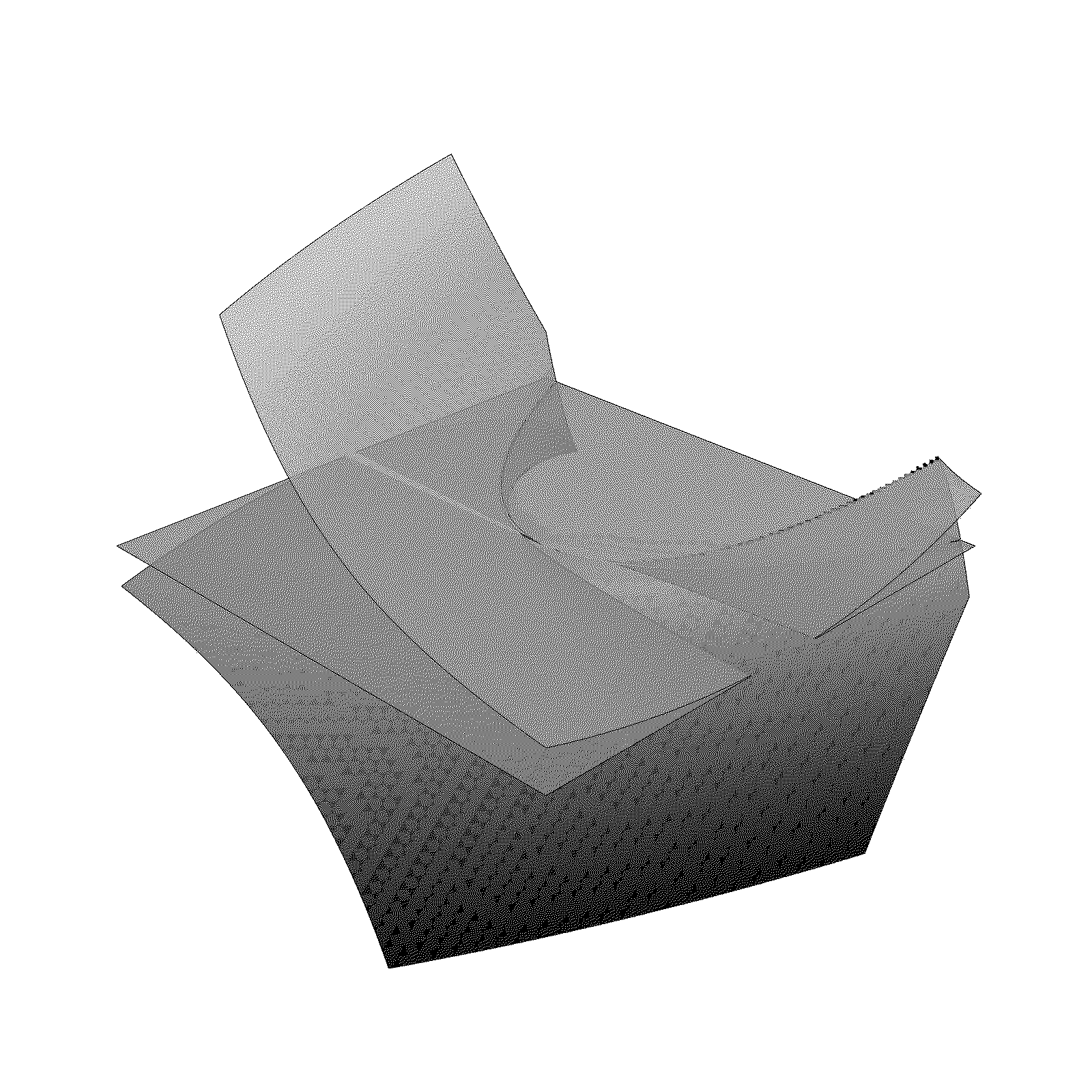}
	\subcaption{Octahedral case}
	\label{fig:discocta}
\end{subfigure}
\hfill
\begin{subfigure}[t]{0.32\textwidth}
	\centering
	\includegraphics[width=\textwidth]{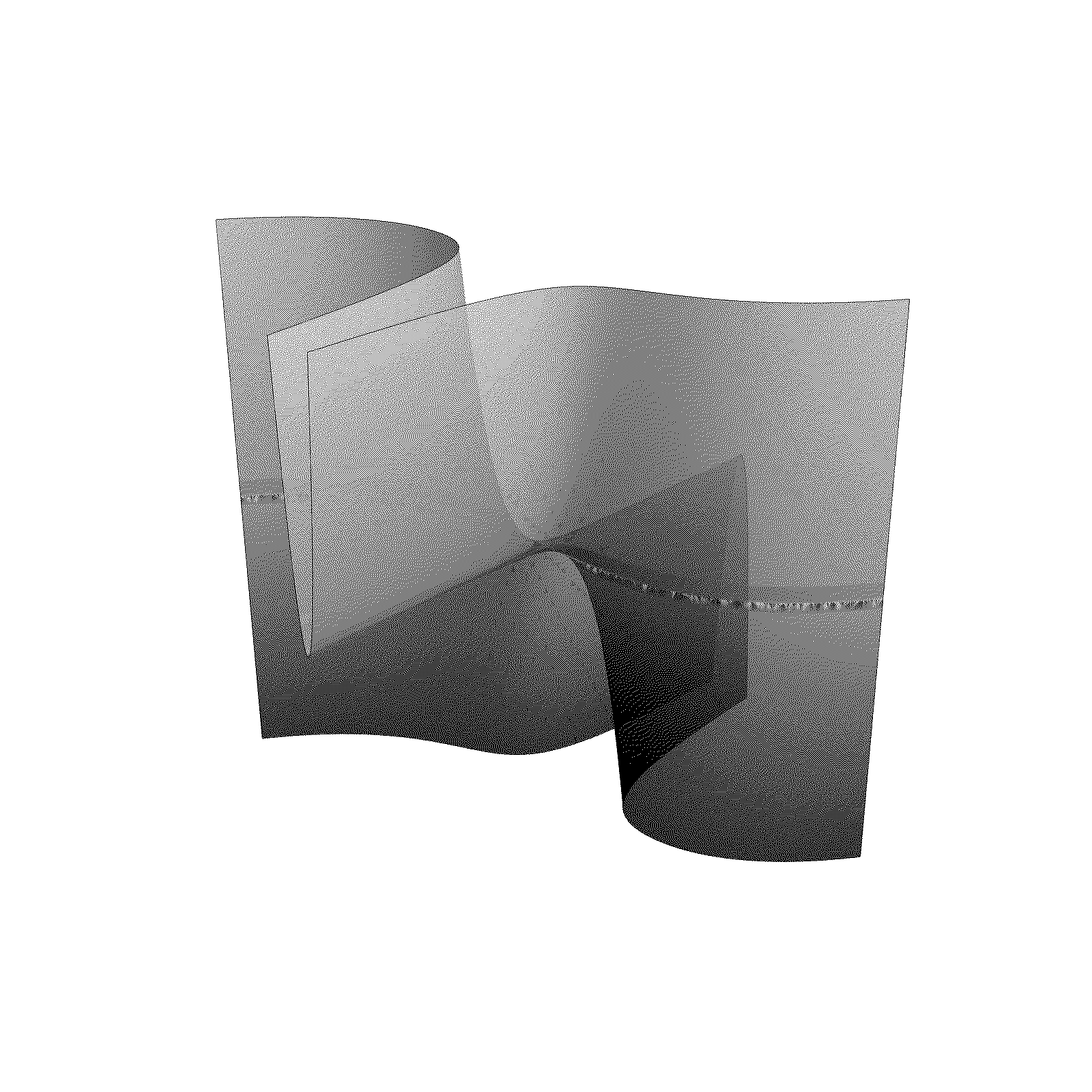}
	\subcaption{Icosahedral case}
	\label{fig:discicosa}
\end{subfigure}

\caption{Discriminant surfaces of the platonic reflection groups.}
\label{fig:discriminant-surfaces}
\end{figure}

Figure \ref{fig:discriminant-surfaces} shows the corresponding {\em discriminant surfaces} $F(P,Q,R)=0$, image of the reflection hyperplanes under the quotient by the reflection groups $G=\T,\O,\I$.

These relations are not recorded in the paper by Goursat, but were determined independently by different authors. We mention \cite{Shcherbak}, but in fact, prior to
Klein and Goursat, Brioschi \cite{Brioschi} computed the discriminant of the binary sextic, which
coincides with the discriminant in the icosahedral case. Certainly, these relations are now routinely
obtained by elimination with the help of a computer algebra system, but these are of little help to understand their true nature. The discriminant of the icosahedral reflection group became known as the {\em caustique mystique} \cite{Bennequin}, as
it was discovered to appear in the ''obstacle problem''. See \cite{Lyashko} for more details.

\subsection{Elliptic threefolds}

It is a remarkable fact that these hypersurfaces can be seen as a {\em Weierstrass
equation in the variables $R$ and $S$}, thus defining a family of elliptic curves, with $P$
and $Q$ as parameters. In other words, the hypersurface given by the equation
$0=S^2-F(P,Q,R)$ is naturally an {\em elliptically fibred threefold}, with the $P-Q$-plane as
basis.

{\bf Corollary:} {\em The quotient of the curve $C_{r,t}$ by the rotation group $SG$ is the elliptic curve $E_{r,t}$ in the $S-R$-plane, defined by the equation
\begin{equation}\label{eq:platonicEG}
E_{r,t}: S^2=F(r,t,R),
\end{equation}
specialised to $P=r$ and $Q=t$, with $F(r,t,R)=$ \eqref{eq:relationtetrahedron}, \eqref{eq:relationoctahedron}, \eqref{eq:relationicosahedron}.
}

\begin{figure}[ht]
	\centering
	\captionsetup[subfigure]{labelformat=empty}
	
	\begin{subfigure}[t]{0.30\textwidth}
		\centering
		\includegraphics[width=\textwidth]{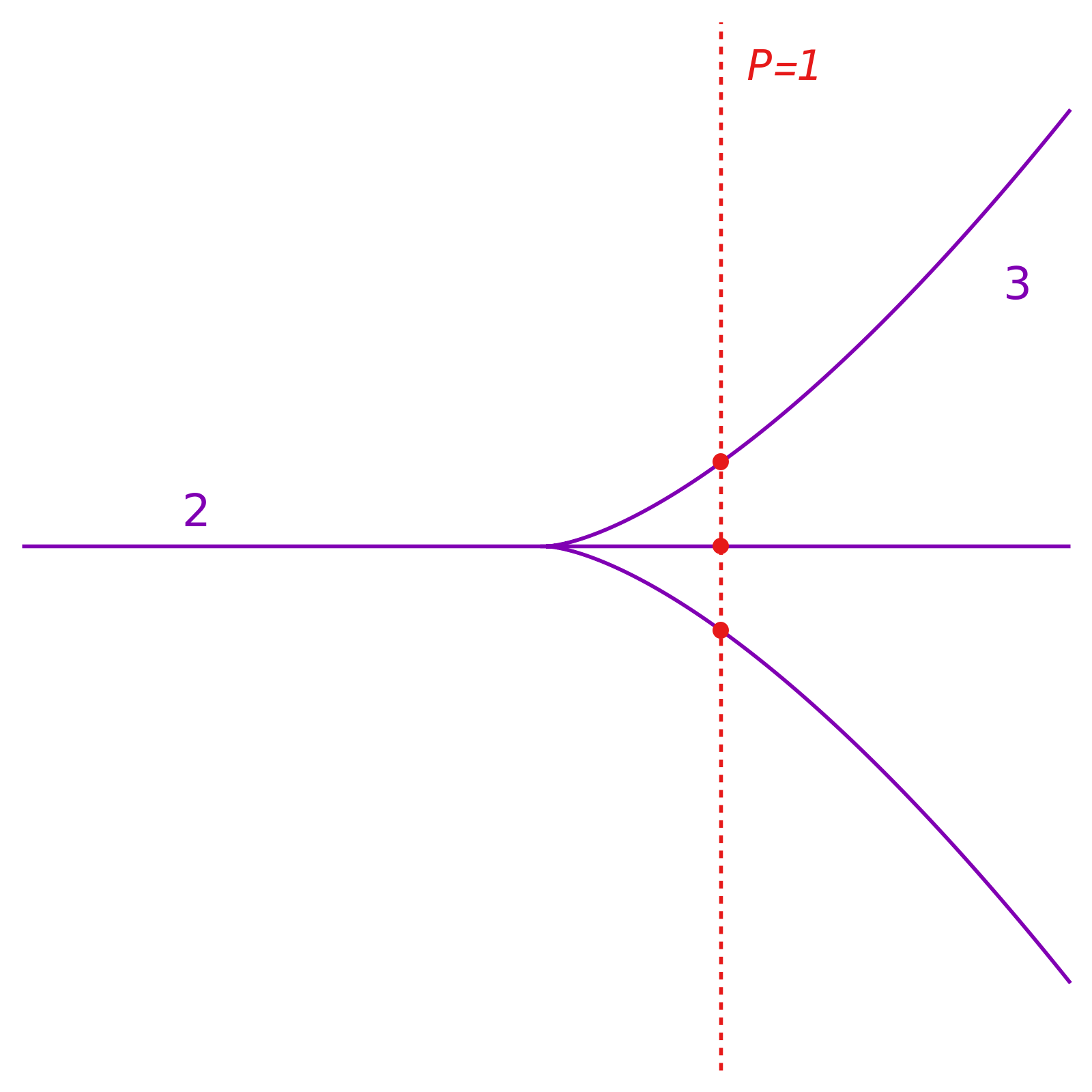}
		\subcaption{Tetrahedral case}
		\label{fig:tetrahedron-curve}
	\end{subfigure}
	\hfill
	\begin{subfigure}[t]{0.30\textwidth}
		\centering
		\includegraphics[width=\textwidth]{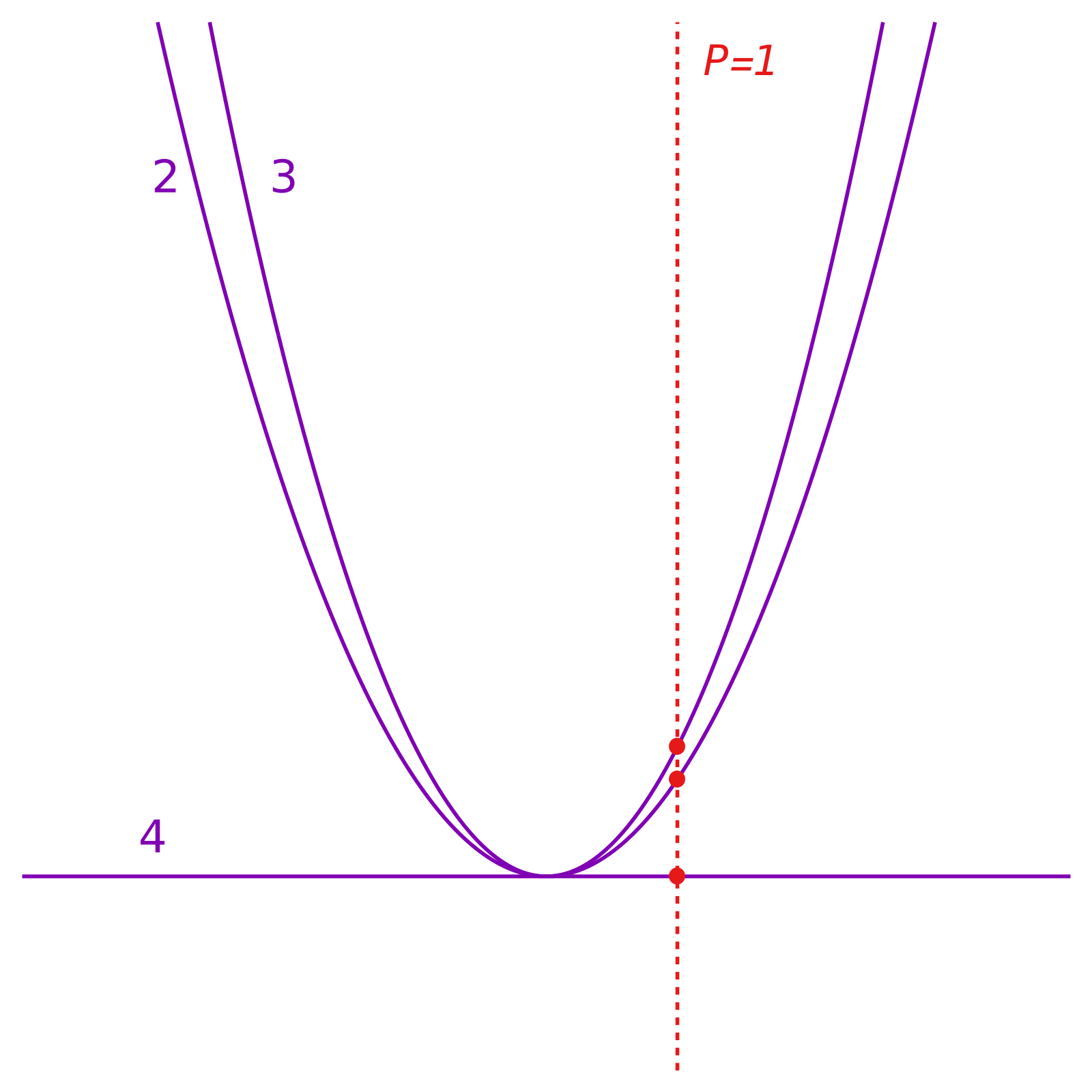}
		\subcaption{Octahedral case}
		\label{fig:octahedron-curve}
	\end{subfigure}
	\hfill
	\begin{subfigure}[t]{0.30\textwidth}
		\centering
		\includegraphics[width=\textwidth]{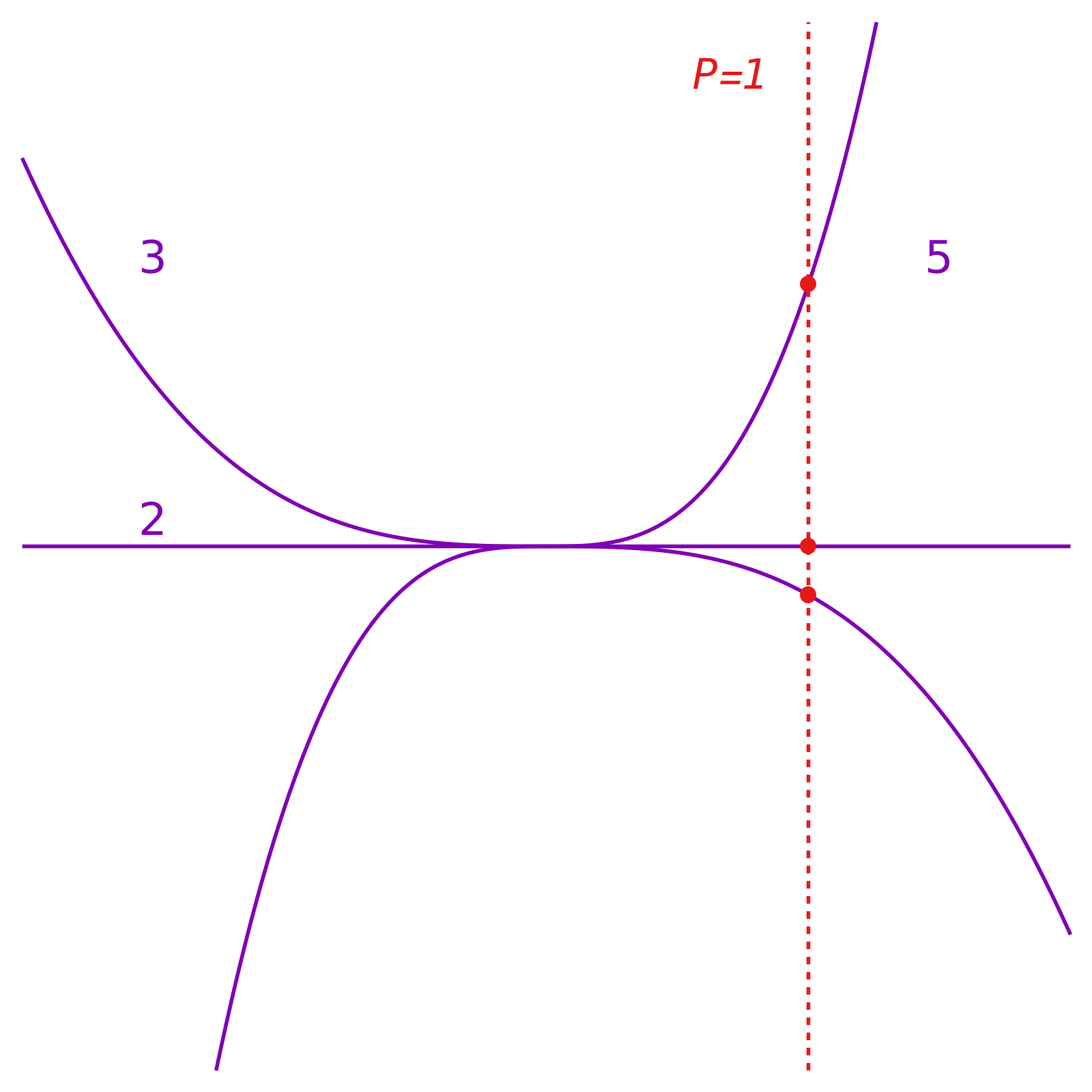}
		\subcaption{Icosahedral case}
		\label{fig:icosahedron-curve}
	\end{subfigure}
	
	\caption{Discriminant curves $\Delta=0$ of $F(P,Q,R)$ in the $P-Q$ plane.}
	\label{fig:platonic-discriminants}
\end{figure}

The singular members in this family are determined by the vanishing of the discriminant $\Delta$ of $F(P,Q,R)$, seen as cubic in $R$ (Figure \ref{fig:platonic-discriminants}). These discriminants are easily computed:

{\bf Tetrahedron} (2,3,3): 
\begin{equation}
\Delta=2^4Q^2(P^3-27Q^2)^3.
\end{equation}

{\bf Octahedron} (2,3,4):  
\begin{equation}
\Delta=2^4Q^4(P^2-4Q)^2(P^2-3Q)^3.
\end{equation}

{\bf Icosahedron} (2,3,5):  
\begin{equation}
\Delta=2^{12}Q^2(5P^3+27Q)^3(P^3-Q)^5.
\end{equation}

\subsection{The platonic elliptic surfaces}\label{ch:platonicellipticsurfaces}

We see that for $P=r$, $r\neq0$, $Q=t$ the discriminant $\Delta$ is of degree $8$ for $\T$, $9$ for $\O$ and $10$ for $\I$ in the variable $t$, with three roots of multiplicities $(2,3,3)$, $(2,3,4)$ respectively $(2,3,5)$, thus recovering the platonic triples (\ref{tab:platonic-reflection-groups}) as multiplicities of the discriminant at the intersection points with the line $P=r$, $r \neq 0$.

Figure \ref{fig:platonic-discriminants} shows the above discriminants with dotted line $P=1$, giving 3 intersection points in the respective cases.

We will see that it is natural to extend these to families of (generalised) elliptic curves  over $\P^1$, that is, consider them as {\em elliptic surfaces}.

The focus on the sections $P=r$, where $r$ is a fixed non-zero value, is quite special. There
are some other interesting sections by lines in the $P-Q$-plane that we will
comment on in Section \ref{ch:outlook}.

\section{Elliptic surfaces}\label{ch:ellipticsurfaces}

By an {\em elliptic surface} over a curve $C$ we mean a smooth projective surface $\mathcal{E}$ equipped with a proper morphism $\pi: \mathcal{E} \lra C$, the general fibre of which is a curve of genus $1$. We always assume that $\pi$ is {\em relatively minimal}, i.e. having no $(-1)$-curves contracted by the map $\pi$. Furthermore, we will assume that a
section $\sigma:C \lra \mathcal{E}$ of $\pi$ is chosen, so that the general fibre has an origin and becomes an elliptic curve. The set of all sections then acquires the structure of an abelian group and is called the {\em Mordell-Weil group} $MW(\mathcal{E})$ of the elliptic surface. The intersection of curves on $\mathcal{E}$ make it into the {\em Mordell-Weil lattice} \cite{OguisoShioda}. Here we will only be dealing with elliptic surfaces over $C=\P^1$.

The elucidation of the basic structure of elliptic surfaces and their degenerate fibres (sometimes called generalised elliptic curves) goes back to Kodaira \cite{Kodaira}. For a modern account we refer to the lecture notes of Miranda \cite{Miranda} and the text book of Sch\"utt and Shioda  \cite{SchuettShioda}.

Elliptic surfaces can be described by {\em Weierstrass models}
\[y^2=x^3+a(t)x+b(t)\textup{ or } y^2=4x^3-g_2(t)x-g_3(t),\]
but often this is not the most natural form.

The singular fibres are determined by the vanishing of the discriminant
$\Delta(t)=g_2(t)^3 -27 g_3(t)^2$; the $J$-map $J: \P^1 \lra \P^1;$ $t \mapsto g_2(t)^3/\Delta(t)$ assigns to $t$ the absolute invariant $J(t)$ that determines (if $\Delta(t) \neq 0$) the elliptic curves $E_t:=\pi^{-1}(t)$ up to $\C$-isomorphism. We consider the restriction
$J:\P^1\setminus S \lra \P^1 \setminus \{0,1,\infty\}$, where $S\supset J^{-1}(\{0,1,\infty\})$ is a finite set, such that each $\pi^{-1}(t)$, $t\in\P^1\setminus S$, is smooth.

We denote by $\widetilde{C}$ the universal cover of $\P^1\setminus S$. For $\widetilde{t}\in \widetilde{C}$ lying above $t\in \P^1\setminus S$,
if we write $E_t=\C/(\Z+\tau(\widetilde{t})\Z)$, $\tau:\widetilde{C}\lra\H=\{z\in\C\mid\mathrm{Im}(z)>0\}$ the period morphism, then $\tau(\widetilde{t})$ is determined up to
a $\mathrm{PSL}_2(\Z)$-transformation. More precisely, we have a map
\[ \rho: \pi_1(\P^1\setminus S)=:G \stackrel{J_*}{\lra} \pi_1(\P^1 \setminus \{0,1,\infty\}) \stackrel{\alpha}{\lra} \mathrm{PSL}_2(\Z),\]
(see \cite{Miranda}, p.61). We call the image $\rho(G)$ in the modular group $\mathrm{PSL}_2(\Z)$ the
{\em (projective) monodromy group}. 

We have the following fact about the monodromy group:
The $J$-map is non-constant if and only if $\rho(G)$ is infinite, in which case it has finite index in $\mathrm{PSL}_2(\Z)$. 

In this case the map $\tau$ induces a factorisation
\begin{equation}\label{eq:factorj}
\P^1\setminus S\simeq\widetilde{C}/G\stackrel{\bar{\tau}}{\lra} \H/\rho(G)\stackrel{J_G}{\lra}\H/\mathrm{PSL}_2(\Z)\simeq\A^1,
\end{equation}
and the composition $J_G\circ\bar{\tau}$ is the $J$-map. Hence \[\deg(\widetilde{\tau})\cdot\deg(J_G)=\deg(J),\]
so we have $[\mathrm{PSL}_2(\Z):\rho(G)]\leq\deg(J)$ (see \cite{UlmerUrzua}). (Actually the index divides the degree.)
If the surface is \emph{modular} in the sense of \cite{TopYui}, then $\deg(\widetilde{\tau})=1$ and we have the equality $[\mathrm{PSL}_2(\Z):\rho(G)]=\deg(J)$.

\subsection{Special elliptic surfaces}
\begin{figure}[ht]
    \centering
    \includegraphics[width=0.30\textwidth]{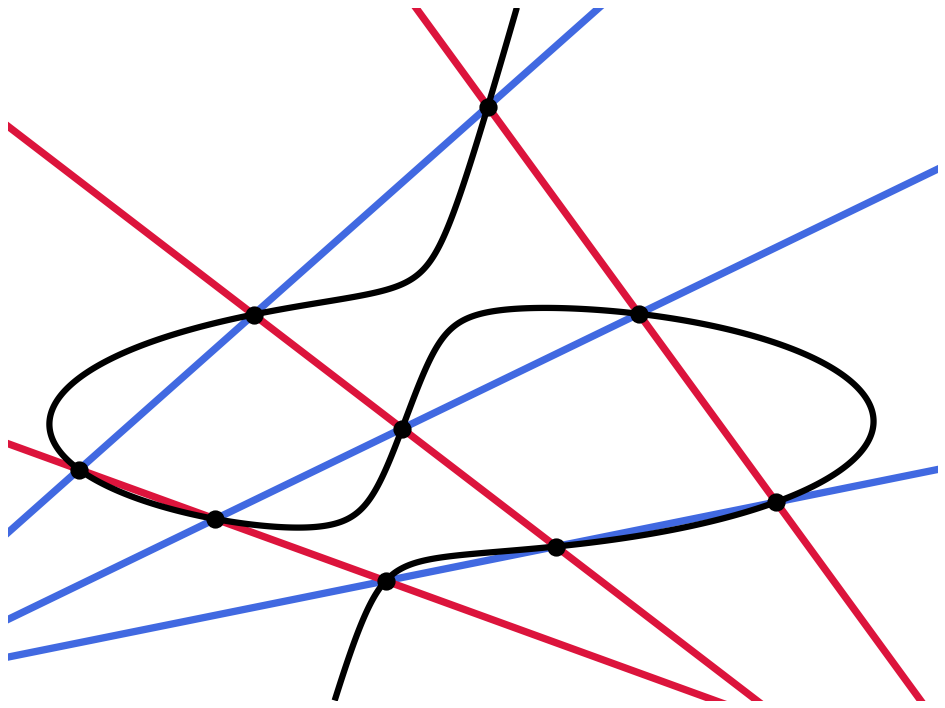}
    \caption{Pencil of cubics spanned by two triples of lines.}
    \label{fig:pencilofcubics}
\end{figure}

Many studies have been made on special classes of elliptic surfaces.
In this paper we will be dealing only with {\em rational elliptic surfaces}
by which one means that $ \mathcal{E}$ is a rational surface, i.e. birational to $\P^2$. The geometrically simplest examples arise from \emph{pencils of cubics}
\[ tf(x,y,z)-sg(x,y,z)=0
\]
by blowing up the $9$ base points defined by the vanishing of $f$ and $g$ (Figure \ref{fig:pencilofcubics}).

For two general cubics, such a surface will have $12$ singular fibres of Kodaira type $I_1$, i.e. ordinary double points. The preimage of the $9$ base points are sections of the fibration and one of them can be chosen as origin of the group law on the fibres; the remaining $8$ sections are generators for the  {Mordell-Weil group} $MW(\mathcal{E})$. The intersection form on divisors gives the Mordell-Weil group the
structure of a \emph{lattice}, which turns out to be of type $E_8$.
By special choice of the cubics $f,g$ one can produce surfaces with fewer singular fibres of various Kodaira types (and other types of  Mordell-Weil groups). Among these the fibres of type $I_n$, which consists of a cycle of $n$-rational curves, are most common and are also called semi-stable. If $e(E_t)$ denotes the Euler number of the fibre $E_t:=\pi^{-1}(t)$, $t \in \P^1$, then one has
$12 =\sum_{t \in \P^1} e(E_t)$, giving the simplest sort of combinatorial restriction on the possible distribution of singular
fibres appearing on a rational elliptic surface.

U. Schmickler-Hirzebruch \cite{Schmickler} compiled the list of rational elliptic surfaces with {\em three singular} fibres and determined the special hypergeometric series that arise as Picard-Fuchs operators attached to these surfaces. Some remarkable surfaces are those related to the so-called {\em euclidean tuples} $(2,2,2,2)$, $(3,3,3)$, $(2,4,4)$ and $(2,3,6)$, whose reciprocals add up to $1$ (Table \ref{tab:euclidean-elliptic-surfaces}).

\begin{table}[ht]
	\centering
	\setlength{\tabcolsep}{4pt}
	\renewcommand{\arraystretch}{1.4}
	
	\begin{tabular}{ccccc}
		\hline
		\textbf{Singular fibres}
		& \textbf{Tuple}
		& \textbf{MW group}
		& \textbf{Modular group}
		& \textbf{Dynkin type} \\
		\hline
		$I_2\,I_2\,I_2^*$
		& --
		& $\mathbb{Z}/2\mathbb{Z}\oplus\mathbb{Z}/2\mathbb{Z}$
		& $\Gamma(2)$
		& -- \\
		\hline
		$I_1\,I_1\,I_4^*$
		& --
		& $\mathbb{Z}/4\mathbb{Z}$
		& $\Gamma_1(4)$
		& -- \\
		\hline
		$I_1\,I_4\,I_1^*$
		& $(2,2,2,2)$
		& $\mathbb{Z}/4\mathbb{Z}$
		& $\Gamma_0(4)$
		& $\widetilde{D}_5$ \\
		\hline\hline
		$I_1\,I_3\,IV^*$
		& $(3,3,3)$
		& $\mathbb{Z}/3\mathbb{Z}$
		& $\Gamma_0(3)$
		& $\widetilde{E}_6$ \\
		\hline
		$I_1\,I_2\,III^*$
		& $(2,4,4)$
		& $\mathbb{Z}/2\mathbb{Z}$
		& $\Gamma_0(2)$
		& $\widetilde{E}_7$ \\
		\hline
		$I_1\,I_1\,II^*$
		& $(2,3,6)$
		& $0$
		& $\Gamma_0(1)$
		& $\widetilde{E}_8$ \\
		\hline
	\end{tabular}
	\caption{Elliptic surfaces attached to the euclidean tuples.}
\label{tab:euclidean-elliptic-surfaces}
\end{table}

The first surface $I_2\;I_2\; I_2^*$ is the {\em Legendre family of elliptic curves} given by the affine equation $y^2=x(x-1)(x-t)$.
The surfaces $I_1\; I_4\; I_1^*$ and $I_1\;I_1\;I_4^*$ are isogenous companions of the Legendre family.
The last three surfaces pull-back to a surface with the nice equation $x^p+y^q+z^r-sxyz=0$ under a cyclic covering $t=s^3, s^4, s^6$.
Note that the starred Kodaira fibres appearing here are exactly the extended Dynkin diagrams of type $\widetilde{D}_5,\widetilde{E}_6,\widetilde{E}_7$ and $\widetilde{E}_8$.

This was followed by the work of Beauville \cite{Beauville}, who gave a description of the six cases with four semi-stable singular fibres, now often called {\em Beauville surfaces}. For convenience of the reader we include a list of the cases in Table \ref{tab:beauville-elliptic-surfaces}.

\begin{table}[ht]
	\centering
	\setlength{\tabcolsep}{6pt}
	\renewcommand{\arraystretch}{1.5}
	
	\begin{tabular}{ccc}
		\hline
		\textbf{Singular fibres}
		& \textbf{Mordell--Weil group}
		& \textbf{Modular group} \\
		\hline
		$I_3\,I_3\,I_3\,I_3$
		& $\mathbb{Z}/3\mathbb{Z}\oplus\mathbb{Z}/3\mathbb{Z}$
		& $\Gamma(3)$ \\
		\hline
		$I_4\,I_4\,I_2\,I_2$
		& $\mathbb{Z}/4\mathbb{Z}\oplus\mathbb{Z}/2\mathbb{Z}$
		& $\Gamma_0(4)\cap\Gamma(2)$ \\
		\hline
		$I_5\,I_5\,I_1\,I_1$
		& $\mathbb{Z}/5\mathbb{Z}$
		& $\Gamma_1(5)$ \\
		\hline
		$I_6\,I_3\,I_2\,I_1$
		& $\mathbb{Z}/6\mathbb{Z}$
		& $\Gamma_0(6)$ \\
		\hline
		$I_8\,I_2\,I_1\,I_1$
		& $\mathbb{Z}/4\mathbb{Z}$
		& $\Gamma_0(8)$ \\
		\hline
		$I_9\,I_1\,I_1\,I_1$
		& $\mathbb{Z}/3\mathbb{Z}$
		& $\Gamma_0(9)$ \\
		\hline
	\end{tabular}
	
	\caption{The Beauville elliptic surfaces with four semi-simple fibres.}
	\label{tab:beauville-elliptic-surfaces}
\end{table}

For these the degree of the $J$-map is $12$, which coincides with the index of the monodromy group
in $\Gamma$. These surfaces were found independently by C. Schoen \cite{Schoen}. All these surfaces have finite Mordell-Weil group and appeared subsequently in the list of {\em extremal surfaces} in the work of Persson-Miranda \cite{MirandaPersson}. These families are modular and can be identified with the universal elliptic curve with the indicated level structure. The number of papers related to these examples is large and quite impossible to overview here.

The classification was pushed further by Herfurtner \cite{Herfurtner}, who compiled the complete list of all surfaces with four singular fibres, given in terms of the Weierstrass covariants $g_2(t)$ and $g_3(t)$. At the same time, U. Persson \cite{Persson} and R. Miranda \cite{Miranda2} gave the complete list of all $279$ possible combinations of Kodaira fibres that may occur on a rational elliptic surface. In a precise sense, all cases with up to $12$ singular fibres of type $I_1$, can be obtained by deformation from surfaces with $\le 4$ fibres, see \cite{Miranda2}.

\subsection{The platonic elliptic surfaces}

What about the three elliptic surfaces obtained in section 1.6 from the Goursat invariants of the $3$-dimensional reflection groups?

{\bf Proposition:}  {\em The fibre types of the three platonic elliptic surfaces \eqref{eq:platonicEG} are summarised in Table \ref{tab:platonic-elliptic-surfaces}}.
\begin{table}[ht]
	\centering
	\renewcommand{\arraystretch}{1.5}
	\setlength{\tabcolsep}{10pt}
	
	\begin{tabular}{c|cccc}
		\hline
		\textbf{Polyhedron}
		& \multicolumn{4}{c}{\textbf{Singular fibres}} \\
		\hline
		Tetrahedron
		& $I_2(0)$
		& $I_3((r/3)^{3/2})$
		& $I_3(-(r/3)^{3/2})$
		& $IV(\infty)$ \\
		\hline
		Octahedron
		& $I_2(r^2/4)$
		& $I_3(r^2/3)$
		& $I_4(0)$
		& $III(\infty)$ \\
		\hline
		Icosahedron
		& $I_2(0)$
		& $I_3(-5r^3/27)$
		& $I_5(r^3)$
		& $II(\infty)$ \\
		\hline
	\end{tabular}
	
	\caption{The platonic elliptic surfaces.}
	\label{tab:platonic-elliptic-surfaces}
\end{table}

{\bf proof}: As we can read off $g_2(t)$ and $g_3(t)$ from the Goursat equations \eqref{eq:platonicEG}, this is a straightforward calculation by hand, PARI or MAGMA. We already have seen that the finite singular fibres arose from a transverse intersection of the
discriminant. The multiplicities of the components were exactly the platonic triples $(p,q,r)$ and thus give semi-stable fibres
$I_p,\;I_q,\;I_r$. By analysing the fibre over $t=\infty$ we find the result as stated. \qed

So we see we nicely recover the platonic triples in the fibre types; the
fibres $II$ (rational curve with cusp), $III$ (two tangent rational curves), $IV$ (three rational curves intersecting in a single point) appearing at infinity are precisely the three other reduced Kodaira-fibres.
Their starred Kodaira-duals $IV^*, III^*, II^*$ correspond to Kodaira fibres that consist of curve configurations of the extended $E_6$, $E_7$ and $E_8$ graphs respectively, which are precisely the Dynkin labels corresponding to $\T,\O,\I$!

As rational elliptic surfaces with four singular fibres, these surfaces appear in the list compiled by Herfurtner \cite{Herfurtner}. Indeed, our three surfaces belong to the group of the $11$ rational elliptic surfaces with {\em reduced fibres only, one of which is unstable}. We call them the \emph{cousins} of our platonic elliptic surfaces and we obtain in this way a rather strange {\em platonic family} of elliptic surfaces of which it is not clear how the members are related precisely (Table \ref{tab:platonic-family}).

\begin{table}[ht]
	\centering
	\setlength{\tabcolsep}{5pt}
	\renewcommand{\arraystretch}{1.5}
	
	\begin{tabular}{c|cccc|c|c|c}
		\hline
		\textbf{Case}
		& \multicolumn{4}{c|}{\textbf{Fibre types}}
		& \textbf{MW lattice}
		& \textbf{Nm.\ in \cite{OguisoShioda}}
		& $\boldsymbol{\deg(J)}$ \\
		\hline
		Tetrahedron
		& $2$ & $3$ & $3$ & $IV$
		& $\langle 1/6\rangle\oplus\mathbb{Z}/3\mathbb{Z}$
		& $61$ & $8$ \\
		
		& $1$ & $2$ & $5$ & $IV$
		& $\langle 1/30\rangle$
		& $56$ & $8$ \\
		
		& $1$ & $1$ & $6$ & $IV$
		& $\langle 1/2\rangle\oplus\mathbb{Z}/3\mathbb{Z}$
		& $51$ & $8$ \\
		\hline
		
		Octahedron
		& $2$ & $3$ & $4$ & $III$
		& $\langle 1/12\rangle\oplus\mathbb{Z}/2\mathbb{Z}$
		& $59$ & $9$ \\
		
		& $1$ & $3$ & $5$ & $III$
		& $\langle 1/30\rangle$
		& $56$ & $9$ \\
		
		& $1$ & $2$ & $6$ & $III$
		& $\langle 1/6\rangle\oplus\mathbb{Z}/2\mathbb{Z}$
		& $53$ & $9$ \\
		
		& $1$ & $1$ & $7$ & $III$
		& $\langle 1/14\rangle$
		& $47$ & $9$ \\
		\hline
		
		Icosahedron
		& $2$ & $3$ & $5$ & $II$
		& $\langle 1/30\rangle$
		& $56$ & $10$ \\
		
		& $1$ & $4$ & $5$ & $II$
		& $\langle 1/20\rangle$
		& $55$ & $10$ \\
		
		& $1$ & $2$ & $7$ & $II$
		& $\langle 1/14\rangle$
		& $47$ & $10$ \\
		
		& $1$ & $1$ & $8$ & $II$
		& $\langle 1/8\rangle$
		& $45$ & $10$ \\
		\hline
	\end{tabular}
	
	\caption{The platonic family ($\langle m\rangle$ means $\Z\sigma$, with $\sigma^2=m$).}
	\label{tab:platonic-family}
\end{table}

We will refer to these elliptic surfaces by their abbreviated Kodaira-label and speak of ''the $127II$ elliptic surface''
for the penultimate case in this list. For example for the three platonic cases we refer to the $234III$ elliptic surface
as the $\O$-elliptic surface etc. Furthermore we use the notation $23\underline{4}II$ to emphasise that the fibre $I_4$ lies above the singular point $0$. The importance of this will become clear later in Section \ref{ch:platonicsequences}.

As in all these $11$ cases the Mordell-Weil group is of rank one, we see that there is a non-torsion section.
It was remarked by Elkies \cite{Elkies} that the elliptic surface $235II=\I$ is the one with the lowest height $=1/30$,
but its direct connection to the icosahedron was not mentioned there. Note that the surface $2 3 3 IV$ has a $3$-torsion section, connecting it with $1 1 6 IV$, and $2 3 4 III$ has a $2$-torsion section, connecting it with $1 2 6 III$. More on this in Section \ref{ch:isogenies}.

\subsection{The Kleinian elliptic surfaces}
There are three further remarkable elliptic surfaces that are naturally attached to the platonic solids via the corresponding  {\emph{binary}} polyhedral groups ${\bf T}$, ${\bf O}$ and ${\bf I}$. They can be understood in terms of the Beauville surfaces as follows.
The Beauville surface $I_3 I_3 I_3 I_3$, also known as the {\em Hesse pencil}, has its four singular fibres lying on the vertices of the regular tetrahedron in the $2$-sphere,  identified with $\P^1$. Three
faces of the tetrahedron come together at a vertex, where we find fibres of type $I_3$. The monodromy group is the full congruence group $\Gamma(3)$ and the surface is identified with the universal elliptic curve with full level $3$ structure over the modular curve $X(3)=\P^1$.

When we pull-back the Beauville surface $I_4 I_4 I_2 I_2$ by a degree $2$ map $\P^1 \lra \P^1$ that ramifies at the $I_2$ fibres, we obtain an elliptic surface with $2+2\times 2=6$ fibres of type $I_4$, which are located on the vertices of a regular octahedron,
where $4$ faces come together. The monodromy group is the full congruence group $\Gamma(4)$ and the
surface can be identified with the universal elliptic curve with full level $4$ structure over $X(4)=\P^1$. However, it is no longer a {\em rational} elliptic surface, as its geometrical genus $p_g=6\cdot 4/12-1=1$; in fact it is a K3 surface.\\
When we pull-back the Beauville surface $I_5 I_5 I_1 I_1$  by a cyclic degree $5$ map $\P^1 \lra \P^1$ ramified over the $I_1$ fibres
we obtain an elliptic surface with $2+5 \times 2=12$ fibres of type $I_5$, which are located on the vertices of a regular icosahedron,
where $5$ faces come together. The monodromy group is the full congruence group $\Gamma(5)$ and the surface can be identified with the universal elliptic curve with full level $5$ structure over the modular curve $X(5)=\P^1$. But of course, this is not a {\em rational} elliptic surface. In fact $p_g=12\cdot 5/12-1=4$.

These three cases relate to the {\em binary} groups ${\bf G}$ rather than to the rotation groups $SG$ we were considering earlier.
If we denote the extended modular group by $\tilde{\Gamma}:=\mathrm{SL}_2(\Z)$, one has ${\bf T} \simeq \mathrm{SL}_2(\Z/3)=\tilde{\Gamma}/\tilde{\Gamma}(3)$, ${\bf O} \simeq \mathrm{SL}_2(\Z/4)=\tilde{\Gamma}/\tilde{\Gamma}(4)$ and ${\bf I}=\mathrm{SL}_2(\Z/5)=\tilde{\Gamma}/\tilde{\Gamma}(5)$. Note that for $n > 5$ the modular curve $X(n)$ is no longer rational, and one obtains elliptic surfaces over curves of higher genus. We refer to \cite{TopYui} for explicit equations for other modular elliptic surfaces.

\section{Normalised Period integrals}\label{ch:normalisedperiodintegrals}

In fact, it was the systematic study of normal forms of elliptic curves by R. Fricke und F. Klein \cite{Klein} that led to the study of families of elliptic curves with level structures over modular curves. Fricke and Klein also described in some detail the Picard-Fuchs differential equations, satisfied by the period integrals attached to the surface. We briefly recall the facts, for more information we refer to \cite{Stiller}.
\subsection{Pull-back of the universal elliptic curve}
Writing the equation of an elliptic surface in Weierstrass form as
\[ y^2=4x^3-g_2(t)x-g_3(t),\]
we extract as before the discriminant $\Delta(t)$ and the $J$-function $J(t)$ as
\[\Delta(t):=g_2(t)^3-27g_3(t)^2,\;\;\;\; J(t):=g_2(t)^3/\Delta(t).\]

We can form period integrals by integrating the canonical differential
\[\omega:= \frac{dx}{y},\;\;\;y=\sqrt{4x^3-g_2(t)x-g_3(t)}\]
over a family of horizontal closed cycles $\gamma_t \in H_1(E_t,\Z)$:
\[ \phi(t):=\int_{\gamma_t} \omega .\]
However, as explained in \cite{Klein}, sometimes it is useful to use a slightly different normalisation which arises by using the {\em normalised differential}
\[ \Omega:=\sqrt{\frac{g_3(t)}{g_2(t)}}\omega,\]
as the corresponding {\em normalised period integrals}
\[\phi(t):=\int_{\gamma(t)}\Omega=\sqrt{\frac{g_3(t)}{g_2(t)}}\int_{\gamma(t)}\frac{dx}{\sqrt{4x^3-g_2(t)x-g_3(t)}}\]
now can be expressed in terms of $J=J(t)$ as
\begin{equation}\label{eq:ellint}
\int_{\gamma(t)} \frac{dx}{\sqrt{4x^3-\frac{27 J}{J-1}(x+1)}},
\end{equation}
and hence can be checked to satisfy the differential equation
\begin{equation}\label{eq:jdiffeq}
	\frac{d^2 \phi}{dJ^2}+\frac{1}{J}\frac{d \phi}{dJ}+\frac{1}{144}\frac{31 J-4}{J^2(J-1)^2}\phi=0.
\end{equation}
The elliptic curve 
\[
y^2=4x^3-\frac{27J}{J-1}(x+1)
\]
appearing in the integral \eqref{eq:ellint} defines an elliptic surface with singular fibres at $0$, $1$ and
$\infty$, with fibres of type $II$, $III^*$ and $I_1$ respectively.
It is called the \emph{universal elliptic curve over the J-line}, as the fibre at $J$ has $J$ as $J$-invariant and thus
each value of $J$ appears precisely once in the family.

The differential equation for the periods of this universal family can be written in terms of the {\em logarithmic derivation}
\[\theta:= t\frac{d}{dt}\]
in the following 
so-called {\em $\theta$-form}:
\[ \left(\theta-\frac{1}{12}\right)\left(\theta-\frac{5}{12}\right)-J\theta^2 = 0, \quad \theta=J\frac{d}{dJ}, \]
which we recognise as hypergeometric operator, with extended Riemann symbol
\[
  \left\{
    \begin{array}{ccc}
      II&III^*&I_1\\
      \hline
      0&1&\infty\\
      \hline
      {1}/{12}&0&0\\
      {5}/{12}&{1}/{2}&0\\
  \end{array}
\right\}.
\]

In order to use this to get a differential operator in the original variable $t$, we have to replace the derivative with respect to $J$ by the derivative with respect to $t$, i.e. we have to set
\begin{equation}\label{eq:dJ}
 \frac{d}{dJ}=\frac{1}{\partial J/\partial t}\frac{d}{dt}.
\end{equation}

\subsection{Some examples}

For the Legendre elliptic curve $y^2=x(x-1)(x-t)$, which defines the rational elliptic surface with singular fibres of type $I_2$, $I_2$, $I_2^*$, one has 
\begin{equation}
	g_2(t)=\frac{4}{3}(t^2 - t + 1), \quad g_3(t)=\frac{4}{27}(2t - 1)(t - 2)(t + 1)
\end{equation}
and
\begin{equation}
	\Delta(t)=16t^2(t-1)^2, \quad J(t)=\frac{4}{27}\frac{(t^2 - t + 1)^3}{t^2(t-1)^2}.
\end{equation}
In order to compute the operator in the variable $t$ one has to replace $J(t)$ in \eqref{eq:jdiffeq} by \eqref{eq:dJ} and scale the result. We get the well-known hypergeometric operator
\begin{equation}\label{eq:legendreddtform}
	\frac{d}{dt}^2+\frac{2t-1}{t(t-1)}\frac{d}{dt}+\frac{1}{4t(t-1)}
\end{equation}
or in $\theta$-form
\begin{equation}\label{eq:legendrethetaform}
\theta^2- t\left(\theta+\frac{1}{2}\right)^2.
\end{equation}
The extended Riemann symbol of this operator is 
\begin{equation}
	\left\{
          \begin{array}{cccc}
            I_2&I_2&I_2^*\\
            \hline
		0 & 1 &\infty\\
		\hline
		0&0&1/2\\
		0&0&1/2\\
	\end{array}
	\right\}.
\end{equation}

Another example is the elliptic surface $y^2-(t+1)xy-ty=x^3+tx^2$, which is the rational elliptic surface with singular fibres of type $I_5$, $I_1$, $I_1$, $I_5$, which has
\begin{equation}
	g_2(t)=\frac{1}{12}(t^4 + 12t^3 + 14t^2 - 12t + 1), \quad g_3(t)=\frac{1}{216}(t^2 + 1)(t^4 + 18t^3 + 74t^2 - 18t + 1)
\end{equation}
and
\begin{equation}
	\Delta(t)=t^5(-t^2-t+1), \quad J(t)=\frac{(t^4 + 12t^3 + 14t^2 - 12t + 1)^3}{1728\cdot t^5(-t^2-t+1)}.
\end{equation}
As differential equation we get
\begin{equation}
	\frac{d}{dt}^2+\frac{3t^2+22t-1}{t(-t^2-t+1)}\frac{d}{dt}-\frac{t+3}{t(-t^2-t+1)}
\end{equation}
or in $\theta$-form
\begin{equation}\label{eq:aperyop}
	\theta^2- t\left(11\theta^2+11\theta+1\right)-t^2(\theta+1)^2.
\end{equation} 

The extended Riemann symbol of this operator is 
\begin{equation}
	\left\{
	\begin{array}{cccc}
    I_5&I_1&I_1&I_5\\
    \hline
		0 & -\frac{11-5\sqrt5}{2} & -\frac{11+5\sqrt5}{2}&\infty\\
		\hline
		0&0&0&1\\
		0&0&0&1\\
	\end{array}
	\right\}.
\end{equation}

\subsection{Operators for the platonic surfaces}

If we apply this to the families $\mathcal{E}_G$ \eqref{eq:platonicEG}, by putting $P=r$ and $Q=t$ and bringing the cubic in Weierstrass form $4x^3-g_2(t)x-g_3(t)$, we read off
the covariants $g_2$ and $g_3$ which then immediately give the differential equations for the normal integrals. The result is:

{\bf Proposition:} {\em  The Picard-Fuchs differential equations 
 satisfied by the normalised period integral for the elliptic surfaces $\T, \O, \I$ are:
\[
\renewcommand{\arraystretch}{1.5}
\begin{array}{cc}
\hline
\textbf{Case}&\textbf{Operator}\\
\hline
  \text{Tetrahedron}& r^3\theta^2 - 3^3t^2(\theta+2/3)(\theta+4/3)\\[1mm]
  \hline
  \text{Octahedron}&4r^4\theta^2 - r^2 t(28\theta^2+28\theta+9)+ 3\cdot 4^2 t^2(\theta+3/4)(\theta+5/4)\\[1mm]
  \hline
\text{Icosahedron}& 20r^6\theta^2+ r^3t(88\theta^2+88\theta+25) -3\cdot 6^2 t^2(\theta+5/6)(\theta+7/6)\\[1mm]
\hline
\end{array}
\]

The extended Riemann symbols of these operators are:

{\bf Tetrahedral operator $\T$:}
\[
\left\{
  \begin{array}{cccc}
    I_2 &I_3 &I_3 &IV\\
    \hline
    0 &(r/3)^{3/2} &-(r/3)^{3/2} &\infty\\
    \hline
    0&0&0&2/3\\
    0&0&0&4/3\\
  \end{array}
\right\}
\]
{\bf Octahedral operator $\O$:}
\[
\left\{
  \begin{array}{cccc}I_4 &I_2 &I_3 &III\\
  \hline
    0 &r^2/4 &r^2/3 &\infty\\
    \hline
    0&0&0&3/4\\
    0&0&0&5/4\\
  \end{array}
\right\}
\]
{\bf Icosahedral operator $\I$:}
\[
\left\{
  \begin{array}{cccc}
    I_2 &I_3 &I_5 &II\\
    \hline
    0 &-5r^3/27 &r/3 &\infty\\
    \hline
    0&0&0&5/6\\
    0&0&0&7/6\\
  \end{array}
\right\}
\]
}

\begin{figure}[H]
\centering
\resizebox{0.8\textwidth}{!}{%
\begin{tikzpicture}[
  x=1cm, y=1cm,
  singfib/.style   ={circle, fill=singcol, inner sep=1.7pt},
  smoothfib/.style ={circle, draw=singcol, fill=white, line width=0.6pt, inner sep=1.4pt},
  tgt/.style       ={circle, fill=tgtcol, inner sep=1.9pt},
  singline/.style  ={singcol, line width=0.7pt},
  smoothline/.style={gray, line width=0.5pt, densely dashed},
  ram/.style       ={font=\footnotesize, color=singcol!85!black, fill=white, inner sep=0.8pt},
  ftype/.style     ={font=\small},
  tlab/.style      ={font=\small},
  jval/.style      ={font=\footnotesize, gray},
]

\draw[gray!50, line width=0.9pt] (0.3,5) -- (12.5,5);
\draw[gray!50, line width=0.9pt] (0.3,0) -- (12.5,0);
\node[font=\small, gray!70, anchor=east] at (0.2,5) {$\mathbb{P}^1_t$};
\node[font=\small, gray!70, anchor=east] at (0.2,0) {$\mathbb{P}^1_J$};

\node[tgt] (T0) at (2,0)    {};
\node[tgt] (T1) at (6.5,0)  {};
\node[tgt] (Ti) at (11,0)   {};
\node[tlab, anchor=north] at (2,-0.16)   {II};
\node[tlab, anchor=north] at (6.5,-0.16) {III$^*$};
\node[tlab, anchor=north] at (11,-0.16)  {I$_1$};
\node[jval, anchor=north] at (2,-0.60)   {$J=0$};
\node[jval, anchor=north] at (6.5,-0.60) {$J=1$};
\node[jval, anchor=north] at (11,-0.60)  {$J=\infty$};

\node[singfib]   (a0) at (0.7,5) {};
\node[smoothfib] (b0) at (1.4,5) {};
\node[smoothfib] (c0) at (2.1,5) {};
\node[smoothfib] (d0) at (2.8,5) {};
\draw[singline]   (a0) -- (T0);
\draw[smoothline] (b0) -- (T0);
\draw[smoothline] (c0) -- (T0);
\draw[smoothline] (d0) -- (T0);
\node[ftype, anchor=south] at (0.7,5.14) {II};
\node[ram] at (0.7,4.72) {$1$};
\node[ram] at (1.4,4.72) {$3$};
\node[ram] at (2.1,4.72) {$3$};
\node[ram] at (2.8,4.72) {$3$};

\node[smoothfib] (p1) at (5.3,5) {};
\node[smoothfib] (p2) at (5.9,5) {};
\node[smoothfib] (p3) at (6.5,5) {};
\node[smoothfib] (p4) at (7.1,5) {};
\node[smoothfib] (p5) at (7.7,5) {};
\draw[smoothline] (p1) -- (T1);
\draw[smoothline] (p2) -- (T1);
\draw[smoothline] (p3) -- (T1);
\draw[smoothline] (p4) -- (T1);
\draw[smoothline] (p5) -- (T1);
\node[ram] at (5.3,4.72) {$2$};
\node[ram] at (5.9,4.72) {$2$};
\node[ram] at (6.5,4.72) {$2$};
\node[ram] at (7.1,4.72) {$2$};
\node[ram] at (7.7,4.72) {$2$};

\node[singfib] (i2) at (10.2,5) {};
\node[singfib] (i3) at (11,5)   {};
\node[singfib] (i5) at (11.8,5) {};
\draw[singline] (i2) -- (Ti);
\draw[singline] (i3) -- (Ti);
\draw[singline] (i5) -- (Ti);
\node[ftype, anchor=south] at (10.2,5.14) {I$_2$};
\node[ftype, anchor=south] at (11,5.14)   {I$_3$};
\node[ftype, anchor=south] at (11.8,5.14) {I$_5$};
\node[ram] at (10.2,4.72) {$2$};
\node[ram] at (11,4.72)   {$3$};
\node[ram] at (11.8,4.72) {$5$};

\node[font=\small] at (4,3.15) {$J:\ t\mapsto J(t)$};
\draw[->, >=stealth, gray!70, line width=0.6pt] (4,2.8) -- (4,2.05);
\node[font=\footnotesize, gray] at (4,1.6) {$\deg J = 10$};

\node[singfib]   at (2.15,-1.5) {};
\node[font=\footnotesize, anchor=west] at (2.30,-1.5) {singular fibre};
\node[smoothfib] at (5.05,-1.5) {};
\node[font=\footnotesize, anchor=west, align=left]
at (5.20,-1.5) {smooth fibre\\(branch point)};
\node[font=\footnotesize, anchor=west] at (9.35,-1.5) {digit $=$ ram.\ index};

\end{tikzpicture}
}
\caption{The $J$-map of the icosahedral surface $I_2\,I_3\,I_5\,\mathrm{II}$ as a
degree$-10$ branched cover of the $J$-line. Filled dots are singular fibres of the
surface, open dots smooth fibres (branch points of $J$); the digit at each preimage
is its ramification index.}
\label{fig:icosa-Jmap}
\end{figure}
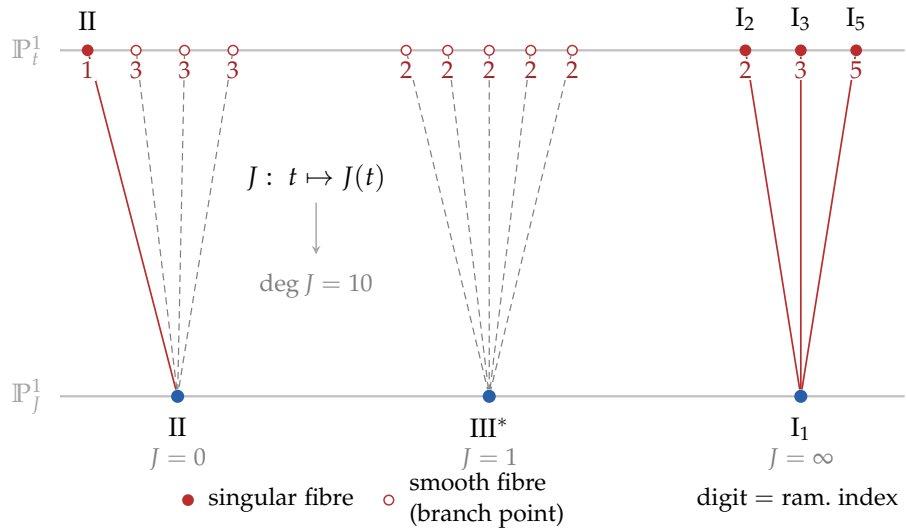

The $J$-map for $\I$ is illustrated in Figure \ref{fig:icosa-Jmap}.

The value of $r$ may be put equal to $1$ or any other value $\neq 0$: the
operators only differ by a scaling of the variable $t$.
We see that the operator for the tetrahedron is obtained as a quadratic pull-back of a corresponding hypergeometric operator. This operator belongs to the elliptic surface with three singular fibres $I_1, I_3, IV^*$ by a covering map
$\P^1 \rightarrow \P^1$, $t\mapsto t^2$, of degree $2$ ramifying over the $I_1$ and $IV^*$ fibres.

Later we will need the notion of the \emph{discriminant} $\Delta(t)=(t-s_1)\dots(t-s_k)$ of an operator (not to be confused with the discriminant of the corresponding elliptic curve family), where $s_i$ denote the singular points $\neq0,\infty$.

\section{The platonic sequences}\label{ch:platonicsequences}

\subsection{The Ap\'ery case}

First discovered in \cite{Apery}, the sequence of (small) Ap\'ery numbers  $1, 3, 19, 147, 1251, 11253,\ldots$, whose $n$th term is the binomial sum
\[ A_n=\sum_{k=0}^n \binom{n}{k}^2\binom{n+k}{k}\]
and satisfies the two term recursion relation
\[(n+1)^2 A_{n+1}= (11n^2+11n+3)A_n+(n-1)^2A_{n-1},\]
which implies that the generating series
\[ \phi(t) =\sum_{n=0}^\infty A_n t^n\] satisfies the second order differential equation
\[\mathcal{A}\phi(t)=0,\;\;\;\mathcal{A}:=\theta^2-t(11\theta^2+11\theta+3)-t^2(\theta+1)^2.\]
It was recognised by F. Beukers \cite{Beukers1} that $\mathcal{A}$ is the Picard-Fuchs operator of the Beauville
elliptic surface $I_5\;I_1\;I_1\;I_5$, which we already mentioned in Section \ref{ch:normalisedperiodintegrals}, \eqref{eq:aperyop}, and $\phi(t)$ is the expansion at the unique holomorphic period integral at a $I_5$ cusp of this surface.

\subsection{Period sequences}
We have seen that the period integrals on any elliptic surface satisfy a second order Picard-Fuchs equation. By a translation we can
shift any singularity to the origin and study the solutions of the differential equation at $0$.
Near a singularity where the elliptic surface acquires a Kodaira fibre of type $I_n$, there is a preferred basis of solutions, called \emph{Frobenius basis}, defined on a slit disc around $0$,
consisting of the {\em holomorphic solution} $\phi$ and the {\em logarithmic solution} $\psi$. Here $\phi(t)$ expands
in a power series, whereas $\psi(t)$ has a logarithmic singularity:
\begin{equation}
\begin{aligned}
\phi(t)&=\sum_n a_nt^n,\;\;a_0=1,
\\
\psi(t)&= \phi(t) \log(t)+f(t),\;\;f(t)=\sum_n b_nt^n,\;\;b_0=0.
\end{aligned}
\end{equation}
Such a point is sometimes called a {\em MUM-point}.
A Picard-Fuchs operator $\mathcal{P}$ of the form ($P_0 \neq 0$, $P_r \neq 0$)
\begin{equation}\label{eq:picardfuchsoperator}
\mathcal{P}=P_0(\theta)+tP_1(\theta)+t^2P_2(\theta)+\dots+t^rP_r(\theta)
\end{equation}
is said to be of {\em degree $r$}. The sequence $a_n$ will then satisfy an $r$-term recursion relation: 
\begin{equation}
	- P_0(n+1)a_{n+1}=P_1(n)a_n+P_2(n-1)a_{n-1}+\dots+P_r(n-r)a_{n-r}.
\end{equation}
It is a theorem formulated by Christol \cite{Christol} that in such a geometrical situation the series $\phi(t)$
has the {\em Eisenstein property}: there exists a natural number $N$ such that $\phi(Nt) \in \Z[[t]]$, i.e. if
we scale the parameter by a factor of $N$, we obtain a series with integral coefficients. 

We call a differential operator and the corresponding holomorphic solution $\phi(t)$ \emph{critically scaled} if
\begin{equation}
\phi(t) \in \Z[[t]], \; \; \phi(t/k) \notin \Z[[t]] \text{ for all } k>1, \text{ and if } a_1\geq0.
\end{equation}
It is customary for an operator
(although not always useful) to replace the parameter $t$ by a parameter $t'=\lambda t$ such that $\phi(t')$ is critically scaled. In this way we can attach a unique integer sequence $a_0=1,a_1,a_2,\ldots$ to each $I_n$ fibre of any elliptic surface.
If the Picard-Fuchs operator is of degree $2$, the sequence satisfies a $2$-term recursion relation.

\subsection{The platonic cases}
We followed this program for each platonic operator. Those have four singular fibres, among these three fibres of
type $I_n$. When we shift such a point $\alpha$ to the origin by $t \mapsto t-\alpha$, there is a unique holomorphic solution to the (shifted) and critically scaled Picard-Fuchs equation of the form
\[\phi(t)=\sum_n a_n t^n,\;\;a_0=1,\;\;\;\phi(t) \in \Z[[t]].\]
{\bf Proposition:} {\em After critical scaling we obtain the $8$ differential operators and period sequences which are summarised in Table \ref{tab:platonicoperators}.}

\clearpage
\newgeometry{
    left=1.2cm,
    right=1.2cm,
    bottom=1.0cm
}

\begingroup

\renewcommand{\arraystretch}{2.2}
\setlength{\tabcolsep}{2pt}
\fontsize{10.2pt}{11pt}\selectfont

\centering
	\medskip
    
	\textbf{Tetrahedral sequences}
	
	\medskip
	
	\noindent\makebox[\textwidth][c]{%
		\begin{tabular}{c|l|c}
			\hline
			
			\multirow{2}{*}{ \begin{tabular}{c} $\T(3)$\\(or $23\underline{3}IV$) \end{tabular}}
			&
			$1,12,198,3720,75690,1626912,36376704,837364176,\dots$
			&
			OEIS: no entry, see \cite{StraubZudilin}
			\\[1mm]
			\cline{2-3}
			
			&
			$\displaystyle
			2\theta^2
			-3t(27\theta^2+27\theta+8)
			+3^6t^2
			\left(\theta+\frac23\right)
			\left(\theta+\frac43\right)$
			&
			$\displaystyle
			I_2\left(\frac{1}{27}\right)\;\;
			I_3\left(\frac{2}{27}\right)\;\;
			I_3(0)\;\;
			IV(\infty)$
			\\[1mm]
			\hline\hline
			
			$2\underline{3}3IV$
			&
			Same as $\T(3)$
			&
			Same as $\T(3)$
			\\[1mm]
			\hline\hline
			
			\multirow{2}{*}{\begin{tabular}{c} $\T(2)$\\(or $\underline{2}33IV$) \end{tabular}}
			&
			$1,0,18,0,810,0,45360,0,2806650,0,183891708,0,\dots$
			&
			OEIS: transformed  \href{https://oeis.org/A006480}%
     {\nolinkurl{A006480}}
			\\[1mm]
			\cline{2-3}
			
			&
			$\displaystyle
			\theta^2
			-3^4t^2
			\left(\theta+\frac23\right)
			\left(\theta+\frac43\right)$
			&
			$\displaystyle
			I_2(0)\;\;
			I_3\left(\frac19\right)\;\;
			I_3\left(-\frac19\right)\;\;
			IV(\infty)$
			\\[1mm]
			\hline
		\end{tabular}%
	}
	
	\medskip
    \medskip
	
	\textbf{Octahedral sequences}
	
	\medskip
	
	\noindent\makebox[\textwidth][c]{%
		\begin{tabular}{c|l|c}
			\hline
			
			\multirow{2}{*}{\begin{tabular}{c} $\O(4)$\\(or $23\underline{4}III$) \end{tabular}}
			&
			$1,12,180,2928,49860,875952,15754704,288722880,\dots$
			&
			OEIS: \href{https://oeis.org/A318245}{\nolinkurl{A318245}}
			\\[1mm]
			\cline{2-3}
			
			&
			$\displaystyle
			3\theta^2
			-4t(28\theta^2+28\theta+9)
			+2^{10}t^2
			\left(\theta+\frac34\right)
			\left(\theta+\frac54\right)$
			&
			$\displaystyle
			I_2\left(\frac{3}{64}\right)\;\;
			I_3\left(\frac{1}{16}\right)\;\;
			I_4(0)\;\;
			III(\infty)$
			\\[1mm]
			\hline\hline
			
			\multirow{2}{*}{\begin{tabular}{c} $\O(3)$\\(or $2\underline{3}4III$) \end{tabular}}
			&
			$1,6,54,588,7110,91476,1224636,16849944,236523078,\dots$
			&
			OEIS: \href{https://oeis.org/A186375}{\nolinkurl{A186375}}
			\\[1mm]
			\cline{2-3}
			
			&
			$\displaystyle
			\theta^2
			-2t(10\theta^2+10\theta+3)
			+2^6t^2
			\left(\theta+\frac34\right)
			\left(\theta+\frac54\right)$
			&
			$\displaystyle
			I_2\left(\frac{1}{16}\right)\;\;
			I_3(0)\;\;
			I_4\left(\frac14\right)\;\;
			III(\infty)$
			\\[1mm]
			\hline\hline
			
			\multirow{2}{*}{\begin{tabular}{c} $\O(2)$\\(or $\underline{2}34III$) \end{tabular}}
			&
			$1,12,612,25392,1298340,67041072,3633970704,\dots$
			&
			OEIS: no entry
			\\[1mm]
			\cline{2-3}
			
			&
			$\displaystyle
			3\theta^2
			-4t(32\theta^2+32\theta+9)
			-2^{12}t^2
			\left(\theta+\frac34\right)
			\left(\theta+\frac54\right)$
			&
			$\displaystyle
			I_2(0)\;\;
			I_3\left(\frac{1}{64}\right)\;\;
			I_4\left(-\frac{3}{64}\right)\;\;
			III(\infty)$
			\\[1mm]
			\hline
		\end{tabular}%
	}
	
	\medskip
    \medskip
	
	\textbf{Icosahedral sequences}
	
	\medskip
	
	\noindent\makebox[\textwidth][c]{%
		\begin{tabular}{c|l|c}
			\hline
			
			\multirow{2}{*}{\begin{tabular}{c} $\I(5)$\\(or $23\underline{5}II$) \end{tabular}}
			&
			$1,10,120,1540,20500,279480,3876600,54496200,\dots$
			&
			OEIS: \href{https://oeis.org/A318495}{\nolinkurl{A318495}}
			\\[1mm]
			\cline{2-3}
			
			&
			$\displaystyle
			2\theta^2
			-t(59\theta^2+59\theta+20)
			+2^4 3^3t^2
			\left(\theta+\frac56\right)
			\left(\theta+\frac76\right)$
			&
			$\displaystyle
			I_2\left(\frac{1}{16}\right)\;\;
			I_3\left(\frac{2}{27}\right)\;\;
			I_5(0)\;\;
			II(\infty)$
			\\[1mm]
			\hline\hline
			
			\multirow{2}{*}{\begin{tabular}{c} $\I(3)$\\(or $2\underline{3}5II$) \end{tabular}}
			&
			$1,30,1440,85260,5606100,391231080,28360117800,\dots$
			&
			OEIS: \href{https://oeis.org/A318496}{\nolinkurl{A318496}}
			\\[1mm]
			\cline{2-3}
			
			&
			$\displaystyle
			10\theta^2
			-3t(333\theta^2+333\theta+100)
			+24\cdot3^6t^2
			\left(\theta+\frac56\right)
			\left(\theta+\frac76\right)$
			&
			$\displaystyle
			I_2\left(\frac{2}{432}\right)\;\;
			I_3(0)\;\;
			I_5\left(\frac{2}{27}\right)\;\;
			II(\infty)$
			\\[1mm]
			\hline\hline
			
			\multirow{2}{*}{\begin{tabular}{c} $\I(2)$\\(or $\underline{2}35II$) \end{tabular}}
			&
			$1,20,1140,68240,4572100,321427920,23420014800,\dots$
			&
			OEIS: no entry
			\\[1mm]
			\cline{2-3}
			
			&
			$\displaystyle
			5\theta^2
			-4t(88\theta^2+88\theta+25)
			-2^8 3^3t^2
			\left(\theta+\frac56\right)
			\left(\theta+\frac76\right)$
			&
			$\displaystyle
			I_2(0)\;\;
			I_3\left(\frac{5}{432}\right)\;\;
			I_5\left(-\frac{1}{16}\right)\;\;
			II(\infty)$
			\\[1mm]
			\hline
		\end{tabular}%
	}

\captionsetup{hypcap=false}
\captionof{table}{Platonic operators and integral sequences.}
\label{tab:platonicoperators}

\endgroup

\restoregeometry

\medskip 
These operators and all but two of the integer sequences appear in the work of  B. Klee \cite{Klee}.
For the integer sequences associated to the 8 cousins (Table \ref{tab:platonic-family})  we refer to  Appendix \ref{ap:herfurtnersurfaces}, Table \ref{tab:herfurtner-operators}.

\subsection{Laurent polynomial representations.}
We say a sequence of integers \\ $a_0,a_1,a_2,\dots$ has a {\em Laurent polynomial
representation} if there exists a Laurent polynomial $F \in \Z[x,y,1/x,1/y]$
such that
\[ a_n=[F^n]_0,\]
where $[-]_0$ denotes the operation of taking the constant term.
Many interesting sequences of integers admit such a representation, which is rather useful
for establishing congruence properties of the numbers appearing in the sequence.

Inspired by the findings of Gorodetsky in \cite{Gorodetsky}, in the first attempts of finding Laurent polynomials for the platonic sequences (Table \ref{tab:platonicoperators}) we made the simple ansatz of assuming a Laurent polynomial $F$ to have Newton polygon of genus $1$ corresponding to one of the $16$ {\em reflexive polygons} given in \cite{Fernando}. Taking general integer coefficients for $F$, iterating over them up to a small bound, computing the powers $F^n$ for the first few
values of $n$ and comparing $[F^n]_0$ to $a_n$  produces candidate Laurent polynomials which can be verfied to produce the sequence $a_n$ for higher values of $n$.

With this ansatz we obtained Laurent polynomial representations for several of the sequences that emerge from the platonic surfaces:

{\bf Proposition:}

{\em The sequence $\T(3)$: $1, 12, 198, 3720, 75690, 1626912, 36376704, 837364176,\dots$ has Laurent
polynomial representation with
\begin{equation}\label{eq:laurentT3}
    F=\frac{(y + 2 + x) (y  - 2 x - 1) (x-2 y-1)}{x y}.
\end{equation}
The sequence $\T(2)$: $1,0,18, 0, 810, 0, 45360, 0, 2806650, 0, 183891708, 0, \dots$ has Laurent
polynomial representation with
\begin{equation}\label{eq:laurentT2}
F=\frac{3 x^2 y^2 + x^2 + 3 x y^2 + 2 x + 1}{x y}.
\end{equation}
The sequence $\O(4)$: $1, 12, 180, 2928, 49860, 875952, 15754704, 288722880,\dots$
has Laurent polynomial representation with
\begin{equation}
F=\frac{x^2 y^2 + 2 x^2 y - 3 x^2 + 2 x y^2 + 12 x y + 2 x - 3 y^2 + 2 y + 1}{x y}.
\end{equation}
The sequence  $\O(3)$: $1, 6, 54, 588, 7110, 91476, 1224636, 16849944, 236523078, \dots$
has two non-equivalent (in the sense of \cite{Fernando}) Laurent polynomial representations with
\begin{align}
    F^{(1)}&=\frac{(x + y - 1)(2x^2 -4xy + x + 2y^2 + y - 1)}{xy},\label{eq:laurentO3_1} \\
    F^{(2)}&=\frac{(y + 2 x + 1) (x y + x + 2 y)}{x y}.\label{eq:laurentO3_2}
\end{align}
The sequence ${\I(5)}$: $1, 10, 120, 1540, 20500, 279480, 3876600, 54496200,\dots$ has Laurent polynomial representation with
\begin{equation}
F=\frac{2 x^2 y^2 + 3 x^2 y - 2 x^2 - x y^2 + 10 x y - 3 x + y + 2}{x y}.
\end{equation}
}

A Laurent polynomial $F$ defines a family of curves $E_t^\circ$ in the torus $(\C^*)^2$:
\[ E_t^\circ :=\{(x,y) \in (\C^*)^2\;\;|\;\; 1-t F=0\}.\]
This curves $E_t^\circ$ can be compactified in $\P^2$ by rationalising and homogenising the equation $1-tF=0$. With this geometric viewpoint one can see that a Laurent polynomial generates an integer sequence, by computing the period integral:
\begin{align}
	\phi(t) 
	&=\frac{1}{(2\pi i)^2}\oint \frac{1}{1-tF}\frac{\mathrm{d}x\mathrm{d}y}{xy} \\
	&=\sum_{n=0}^\infty [F(x,y)^n]_0.
\end{align}

\subsection{Elliptic surfaces from Laurent polynomials}
It is useful to replace $t$ by $t/s$, so that we end up with a \emph{pencil} over $\P^1 \ni (t:s)$.

\begin{figure}[htbp]
	\centering
	\makebox[\textwidth][c]{%
		\includegraphics[height=12cm]{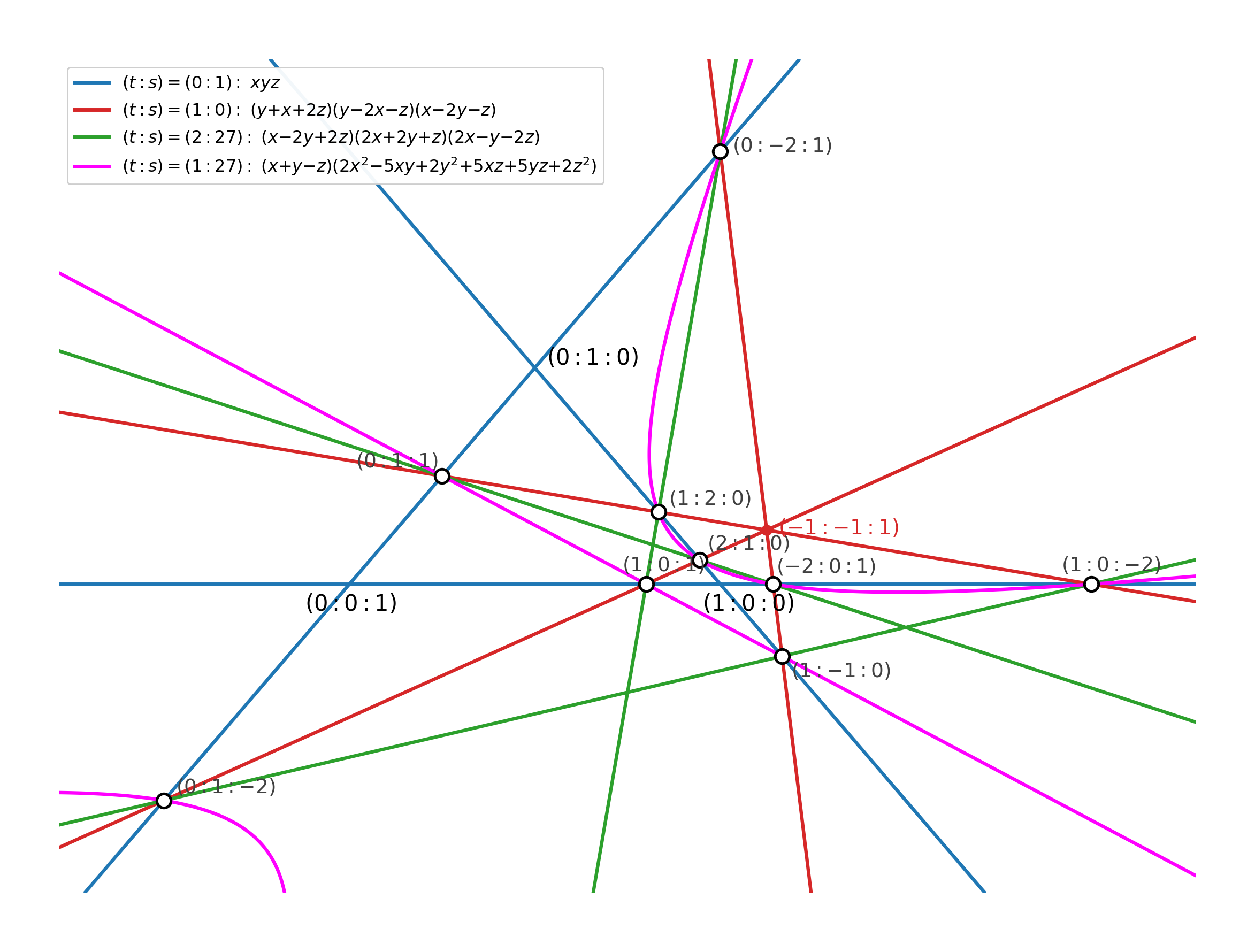}%
	}
	\caption{Reducible fibres of the pencil\\
		\(t\cdot(x+y+2z)(-2x+y-z)(x-2y-z)-s\cdot xyz\).}
	\label{fig:reducible-fibres}
\end{figure}

As an example we consider the Laurent polynomial $F_{\T(3)}$ \eqref{eq:laurentT3}, denoting the one that generates the sequence $\T(3)$. This
leads to a pencil of cubics
\begin{equation}\label{eq:pencilforT3}
C(t,s): \;  t (y + x+ 2z) (y - 2 x - z) (x-2 y - z)-s xyz=0.
\end{equation}
The following observations are illustrated in (Figure \ref{fig:reducible-fibres}): 

For the point $(t:s)=(0:1)$ we obtain the coordinate lines $xyz=0$ (\emph{blue}), whereas for $(t:s)=(1:0)$ we obtain three further
lines  $(y + x+ 2z)(y - 2 x - z)(x-2 y - z)=0$ (\emph{red}), all passing through the point $(-1:-1:1)$
and intersecting the coordinate triangle in $9$ distinct points.
All cubics $C(t,s)$ of the pencil pass through 
these nine points. 

There are two further reducible elements in the pencil: For the point $(t:s)=(2:27)$ we obtain three further lines
$(-2 y + 2z + x) (2 y + z + 2 x) (-y - 2z + 2 x)$ (\emph{green})
and for $(t:s)=(1:27)$ the curve degenerates to a line and a conic: $(y - z + x) (2 x^2  - 5 x y + 2 y^2  + 5 xz + 5 yz + 2z^2)$ (\emph{magenta}).

Hence from the pencil $C(t,s)$ \eqref{eq:pencilforT3} we obtain an elliptic surface with Kodaira fibres $I_3$ at $0$, $I_3$ at $2/27$, $I_2$ at $1/27$ and $IV$ at $\infty$. Note that we have made a complete circle: starting from the platonic elliptic surface
$2 3 \underline{3} IV$, we determined its Picard-Fuchs equation from which we obtained an integer sequence $\T(3)$ for which
we found a Laurent polynomial, which gave us back the elliptic surface we started with:
\[
\tikz[baseline=(current bounding box.center), scale=1.2, transform shape]{
  \node (A) at (90:1.3cm) {\scriptsize $23\underline{3}IV$};
  \node (B) at (180:1.3cm) {\scriptsize $PF$};
  \node (C) at (-90:1.3cm) {\scriptsize $\T(3)$};
  \node (D) at (0:1.3cm) {\scriptsize $F_{\T(3)}$};

  \draw[->, decorate, decoration={snake, amplitude=0.4mm, segment length=3mm}] 
    (A) to[bend right=18] (B);
  \draw[->, decorate, decoration={snake, amplitude=0.4mm, segment length=3mm}] 
    (B) to[bend right=18] (C);
  \draw[->, decorate, decoration={snake, amplitude=0.4mm, segment length=3mm}] 
    (C) to[bend right=18] (D);
  \draw[->, decorate, decoration={snake, amplitude=0.4mm, segment length=3mm}] 
    (D) to[bend right=18] (A);
}
\]

\subsection{Isogenies}\label{ch:isogenies}
Somewhat to our surprise, {\em this does not happen always} with our Laurent polynomials. For the sequence $\T(2)$ it turns out that the Laurent polynomial $F_{\T(2)}$ \eqref{eq:laurentT2} given above defines an elliptic surface of fibre type $I_6$ at $0$,  $I_1$ at $\pm \frac{1}{9}$ and $IV$ at $\infty$. Instead of the circle we have here
\[ \underline{2} 3 3 IV \rightsquigarrow PFE \rightsquigarrow \T(2) \rightsquigarrow F_{\T(2)} \rightsquigarrow 11\underline{6} IV,\]
so that we arrive at a {\em different} elliptic surface! The issue here is that the surfaces
$2 3 3 IV$ and $1 1 6IV$ are connected by \emph{degree $3$-isogeny}, which turns
$I_2$ into $I_6$ and $I_3$ into $I_1$. Two such isogenous surfaces have the same Picard-Fuchs operator and thus produce the same integer sequences. 

For the sequence $\O(3)$ we have two non-equivalent Laurent polynomials. It turns out that they describe the isogenous surfaces $2\underline{3}4 III$ ($F^{(1)}$, \eqref{eq:laurentO3_1}) and $12\underline{6} III$ ($F^{(2)}, \eqref{eq:laurentO3_2})$: 
\[
\tikz[baseline=(current bounding box.center), scale=1.2, transform shape]{
  \node (A) at (90:1.3cm) {\scriptsize $2\underline{3}4\,III$};
  \node (B) at (180:1.3cm) {\scriptsize $PFE$};
  \node (C) at (-90:1.3cm) {\scriptsize $\O(3)$};
  \node (D) at (0:1.3cm) {\scriptsize $F^{(1)}_{\O(3)}$};
  \node (E) at ($(C)+(1.7cm,0)$) {\scriptsize $F^{(2)}_{\O(3)}$};
  \node (F) at ($(E)+(1.7cm,0)$) {\scriptsize $12\underline{6}\,III$};

  \draw[->, decorate, decoration={snake, amplitude=0.4mm, segment length=3mm}]
    (A) to[bend right=18] (B);
  \draw[->, decorate, decoration={snake, amplitude=0.4mm, segment length=3mm}]
    (B) to[bend right=18] (C);
  \draw[->, decorate, decoration={snake, amplitude=0.4mm, segment length=3mm}]
    (C) to[bend right=18] (D);
  \draw[->, decorate, decoration={snake, amplitude=0.4mm, segment length=3mm}]
    (D) to[bend right=18] (A);

  \draw[->, decorate, decoration={snake, amplitude=0.4mm, segment length=3mm}]
    (C) -- (E);
    \draw[->, decorate, decoration={snake, amplitude=0.4mm, segment length=3mm}]
    (E) -- (F);
}
\]

The elliptic surfaces $234III$ and $126III$ are connected by an isogeny of degree $2$.

Since we have the two pairs of isogenous surfaces 
\[233IV\stackrel{\text{isog.}}{\longleftrightarrow}116IV \; \text{ and } \; 234III\stackrel{\text{isog.}}{\longleftrightarrow}126IV\] 
and two non-equivalent Laurent polynomials for $\O(3)$ (i.e. $2\underline{3}4III)$, the question arises whether we can find two Laurent polynomials for each sequence $\T(3),\T(2),\O(4),\O(3)$ and $\O(2)$. This will be answered in Section \ref{ch:remainingcases}.

\subsection{From quartics to cubics}
The Laurent polynomials $F_{\T(2)}$, $F_{\O(4)}$ and $F_{\I(5)}$ produce pencils of {\em quartics}.
For example, for $F_{\T(2)}$ we obtain the pencil of quartic curves
\[Q(t,s):  s x y z^2-t (x+z)(3xy^2z+xz^2+z^3)=0.\]
The general homogeneous quartic defines a curve of genus $3$, but it is easily
verified that for general $(t:s) \in \P^1$ the curve $Q(t,s)$ has two ordinary
double points located at coordinate points $(0:1:0)$ and $(0:0:1)$, so that the
genus of $C(t,s)$ equals $1$. These curves can be transformed into cubics by a
Cremona transformation. For this one considers the linear system of conics through
the two singular points and one of the remaining base points of the pencil.
In concrete terms, after the substitution
\[ x=YZ,\;\;\;y = X(X+Z),\;\;\;z = XY\]
in the equation, a factor $X^2Y^2(Z+X)$ can be removed and one obtains a pencil
of cubics
\[C(t,s): s XYZ-t (X+Z)(3XZ+Y^2+3Z^2)=0.\]

It is readily verified that this produces the alternative Laurent polynomial
\[ F_{\T(2)}':=\frac{(X+Z)(3XZ+Y^2+3Z^2)}{XY}\]
for the sequence $\T(2)$.

The same can be done for the Laurent polynomial $F_{\O(4)}$ and $F_{\I(5)}$.
For $F_{\O(4)}$ we make the substitution
\[ x=YZ,\;\;\;y=X(3Z+X),\;\;\;z=XY,\]
which transforms the quartic pencil of $F_{\O(4)}$ into the cubic pencil
\[D(t,s):=sXYZ -t (3 X^2  + X Y + 10 X Z - Y Z + 3 Z^2) (Y+Z-X)\]
and hence we get an alternative representation of the sequence $\O(4)$ by the Laurent
polynomial
\[ F_{\O(4)}':=\frac{(3 X^2  + X Y + 10 X - Y + 3 ) (Y+1-X)}{XY}.\]

In a similar vein, we use the substitution
\[ x=YZ,\;\;y=X(2Z-X),\;\;\;z=XY\]
to transform the quartic pencil of $F_{\I(5)}$ into the cubic pencil
\[s XYZ -t (X^2Y+X^2Z-2XY^2+10XYZ-4XZ^2-Y^2Z+3YZ^2+4Z^3)\]
and we obtain the corresponding alternative
Laurent representation for the sequence $\I(5)$.
\[F_{\I(5)}':=\frac{X^2Y+X^2-2XY^2+10XY-4X-Y^2+3Y+4}{X Y}.\]

\subsection{The remaining cases}\label{ch:remainingcases}

What about the remaining sequences $\O(2)$, $\I(3)$ and $\I(2)$? And what about the alternative Laurent polynomials for $\T(3)$, $\T(2)$ and $\O(4)$?

As we have by now aquired a deeper insight into the nature of
Laurent representations, we can directly settle case $\I(3)$.

{\bf Proposition:}
{\em The integer sequence $\I(3)$ does admit a Laurent polynomial representation:}
\[
F=\frac{2U^3-3U^2V-12UV^2+20V^3-3U^2+30UV-21V^2-12U-21V+20}{U V}.
\]

{\bf proof:} We start with the cubic Laurent polynomial $F_{\I(5)}'$ for the sequence $\I(5)$. The associated pencil of cubics has six base-points, three of them with multiplicity $2$: $(0:1:0),(1:0:0),(2:0:1)$ and three of them with multiplicity $1$: $(0:-1:1),(0:4:1),(2:1:0)$. The sequence $\I(3)$ is associated to the expansion of the normalised period at $t=2/27$.
The pencil determined by $F_{\I(5)}'$ at $(t:s)=(2:27)$ factorises as
\[(2 X  Y + 2 X Z + Y Z - 4 Z^2 ) (X - 2 Y - 2 Z)\]
which defines a line $L: X-2Y-2Z=0$ and a conic $Q: 2 X  Y + 2 X Z + Y Z - 4 Z^2=0$, which intersect at $(2:0:1)$
and $(-10:-9:4)$.
We now apply a quadratic Cremona transformation on the three base points
$(0:1:0)$, $(2:0:1)$, $(0:4:1)$. This transforms cubics to cubics ($3=2\times 3-3$),
the conic $Q$ to a line ($1=2\times2-3)$ and the line $L$ into another line ($1=2\times 1 -1$).
The point $(2:0:1)$ is blown up and produces a third line, so after the transformation we have
three lines over $2/27$.
To compute it, we first replace $X$ be $X+2Z$ and then make the substitution
$X=VW,\;\;\; Y=WU,\;\;\; Z=UV$. After dividing by $UVW$ we obtain a cubic pencil which for
$t=2, s=27$ factors into three linear factors. Now we make the substitution $t\mapsto27t+2s$, $s\mapsto27s$,
shifting the $2/27$ fibre to $0$. By a linear transformation in $U,V,W$ we transform the pencil
into the form
\begin{align*} t(2U^3-3U^2V-12UV^2+20V^3-3U^2W+30UVW-21V^2W \\
-12UW^2-21VW^2+20W^3)-sUVW=0
\end{align*}
So we obtain the desired Laurent polynomial. \qed

For the remaining two cases $\O(2)$ and $\I(2)$ there is no Laurent polynomial representation of genus $1$. With the help of ChatGPT and Claude we were able to find genus-$3$ and genus-$10$ Laurent polynomials for those cases and for $\T(3)$, $\T(2)$ and $\O(4)$.

{\bf Proposition:} 
{\em
The sequence $\T(3)$: $1, 12, 198, 3720, 75690, 1626912, 36376704, 837364176,\dots$ has a second Laurent polynomial representation corresponding to the surface $23\underline{3}IV$ with
\begin{equation}
\begin{aligned}
F={}&
\frac{
-1-6y+3y^2+52y^3+12y^4-96y^5-64y^6
+3xy+21xy^2}{
\phantom{x^2y}
}
\\[-1.8ex]
&\quad\frac{+18xy^3-96xy^4-96xy^5+12x^2y+9x^2y^2-72x^2y^3-96x^2y^4
}{
\phantom{x^2y}
}
\\[-1.8ex]
&\quad
\frac{
-2x^3+3x^3y-24x^3y^2-56x^3y^3
-6x^4y-24x^4y^2-6x^5y-x^6
}{
x^2y
}.
\end{aligned}
\end{equation}
}

{\em
The sequence $\T(2)$: $1,0,18, 0, 810, 0, 45360, 0, 2806650, 0, 183891708, 0, \dots$ has a second Laurent polynomial representation corresponding to the surface $\underline{2}33IV$ with
\begin{equation}
F=\frac{
9x^4-3x^2y^2+3x^2y+6x^2+y^4+y^3+y+1
}{
xy^2}.
\end{equation}
}

{\em
The sequence $\O(4)$: 
$1,12,180,2928,49860,875952,15754704,288722880,\dots$ has a second Laurent polynomial representation corresponding to the surface $1\underline{2}6III$ with
\begin{equation}
F=
\frac{
\bigl(x^2-2xy-x-3y^2+y\bigr)
\bigl(x^3+x^2y-x^2-5xy^2-2xy+3y^3-y^2\bigr)
}{
x^2y^2
}.
\end{equation}
}

{\em
The sequence $\O(2)$: $1,12,612,25392,1298340,67041072,3633970704,\dots$ has two Laurent polynomial representations corresponding to the isogeneous surfaces $\underline{2}34III$ and $\underline{1}26III$ with

\begin{equation}
\begin{aligned}
F_{\underline{2}34III}&=-
\frac{
\bigl(3x^2-2xy-x+3y^2-y\bigr)
\bigl(9x^3-5x^2y-3x^2-5xy^2+2xy+9y^3-3y^2\bigr)
}{
x^2y^2
},\\
F_{\underline{1}26III}&=
\frac{-9x^6+5x^4y^2+6x^4y-27x^4+5x^2y^4+12x^2y^3-2x^2y^2
      }{\phantom{x^2y^3}}\\
\\[-1.8ex]
&\quad
\frac{+12x^2y-27x^2-9y^6+6y^5-7y^4+20y^3
      -7y^2+6y-9
}{
x^2y^3
}.
\end{aligned}
\end{equation}
}

{\em
The sequence $\I(2)$: $1,20,1140,68240,4572100,321427920,23420014800,\dots$ has Laurent polynomial representation with

\begin{equation}
\begin{aligned}
F
={}&
\frac{
125x^4-60x^3y-100x^3
+6x^2y^2+4x^2y+70x^2-12xy^3
}{
\phantom{xy^2}
}
\\[-1.8ex]
&\quad
\frac{
+20xy^2-4xy-20x
+5y^4-4y^3-2y^2-4y+5
}{
xy^2
}.
\end{aligned}
\end{equation}
}

\begin{figure}[htbp]
	\centering
	\captionsetup[subfigure]{labelformat=empty}
	\noindent\makebox[\textwidth][c]{%
    \begin{subfigure}[t]{0.30\textwidth}
        \centering
        \includegraphics[width=\textwidth]{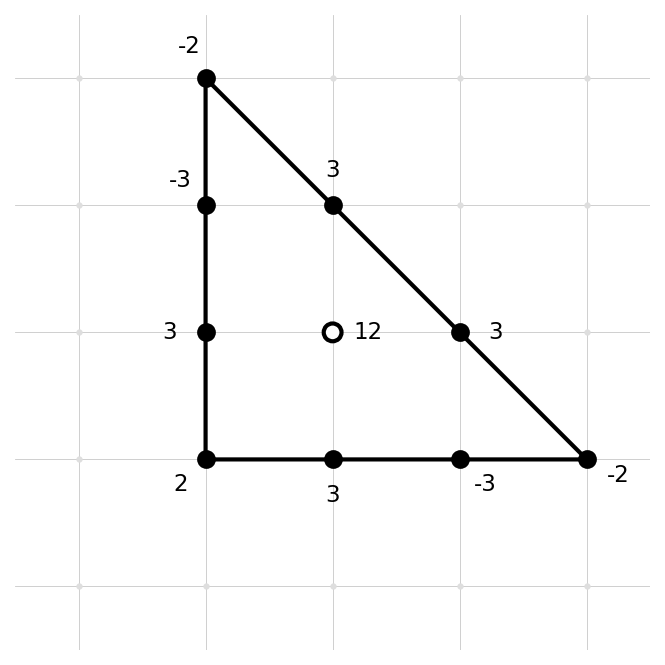}
        \subcaption{Newton polygon of $F_{\T(3)}^{23\underline{3}IV}$}
    \end{subfigure}
    \hspace{0.04\textwidth}
    \begin{subfigure}[t]{0.30\textwidth}
        \centering
        \includegraphics[width=\textwidth]{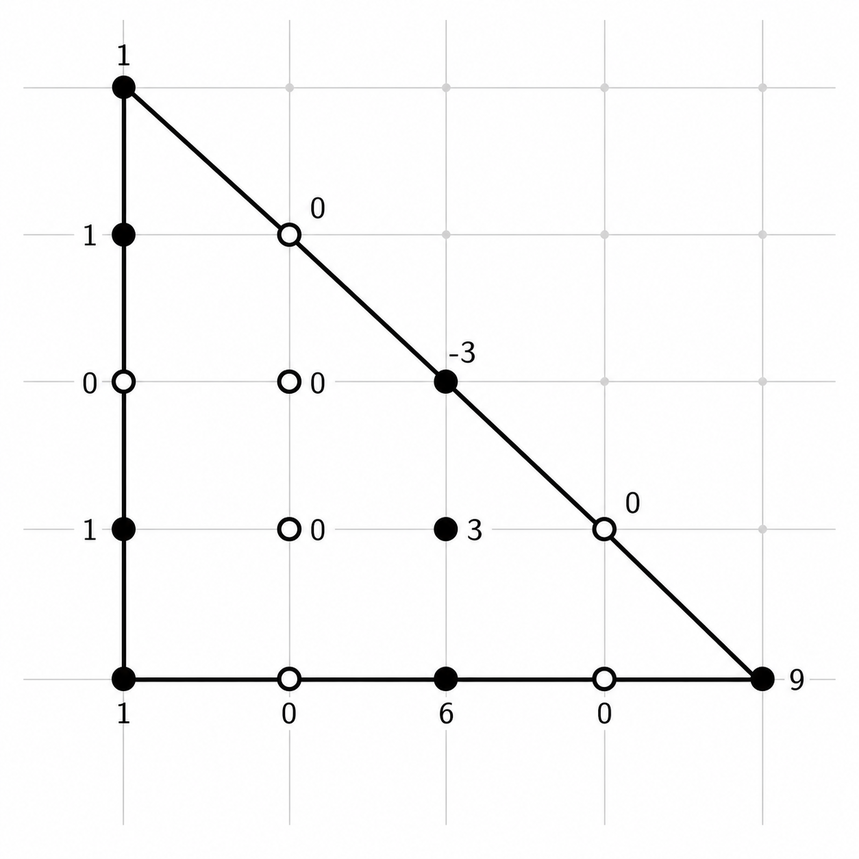}
        \subcaption{Newton polygon of $F_{\T(2)}^{\underline{2}33IV}$}
    \end{subfigure}%
}
\medskip
	\noindent\makebox[\textwidth][c]{%
    \begin{subfigure}[t]{0.30\textwidth}
        \centering
        \includegraphics[width=\textwidth]{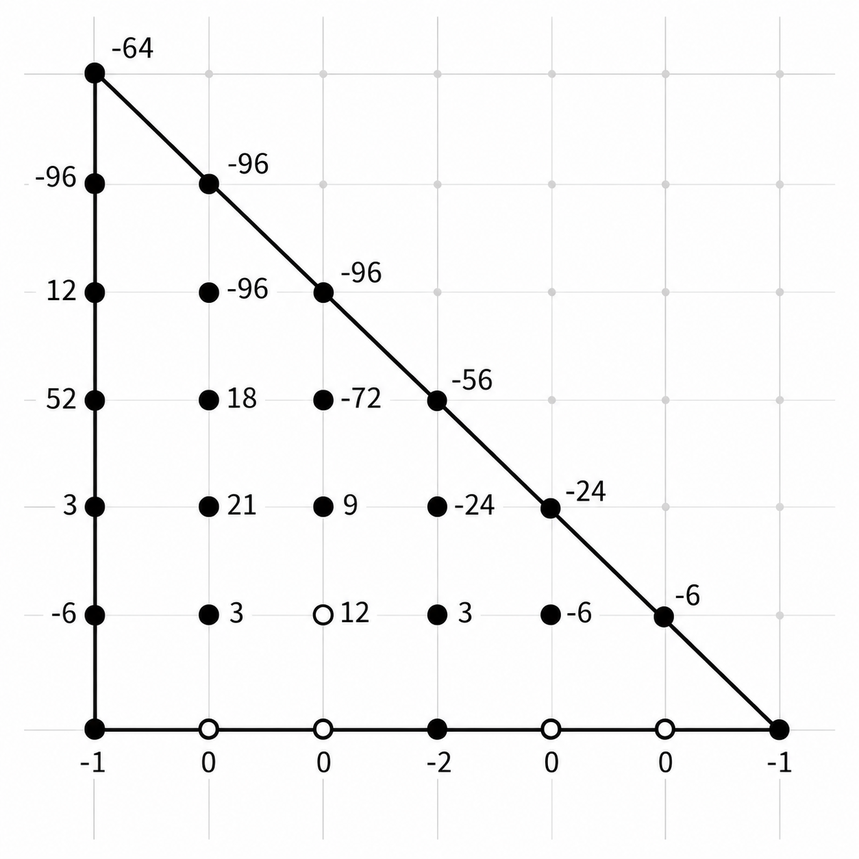}
        \subcaption{Newton polygon of $F_{\T(3)}^{1\underline{1}6IV}$}
    \end{subfigure}
    \hspace{0.04\textwidth}
    \begin{subfigure}[t]{0.30\textwidth}
        \centering
        \includegraphics[width=\textwidth]{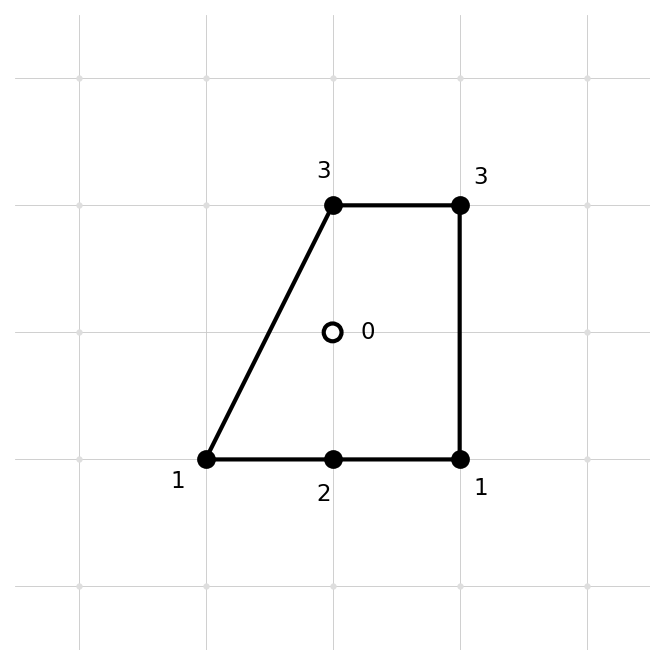}
        \subcaption{Newton polygon of $F_{\T(2)}^{11\underline{6}IV}$}
    \end{subfigure}%
}
	\caption{Newton polygons of the Laurent polynomials for $\T$.}
	\label{fig:newton-polygons_T}
\end{figure}

\begin{figure}[htbp]
	\centering
	\captionsetup[subfigure]{labelformat=empty}
	\noindent\makebox[\textwidth][c]{%
	\begin{subfigure}[t]{0.30\textwidth}
		\centering
		\includegraphics[width=\textwidth]{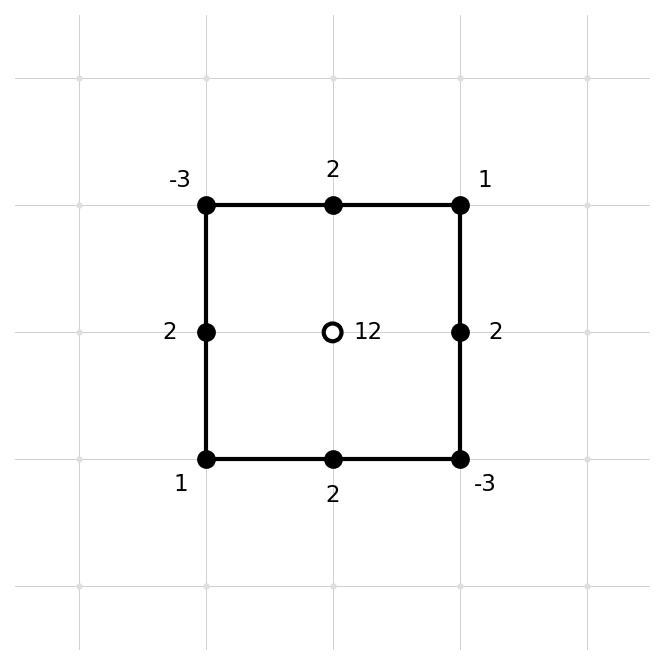}
		\caption{Newton polygon of $F_{\O(4)}^{23\underline{4}III}$}
	\end{subfigure}
	\hfill
	\begin{subfigure}[t]{0.30\textwidth}
		\centering
		\includegraphics[width=\textwidth]{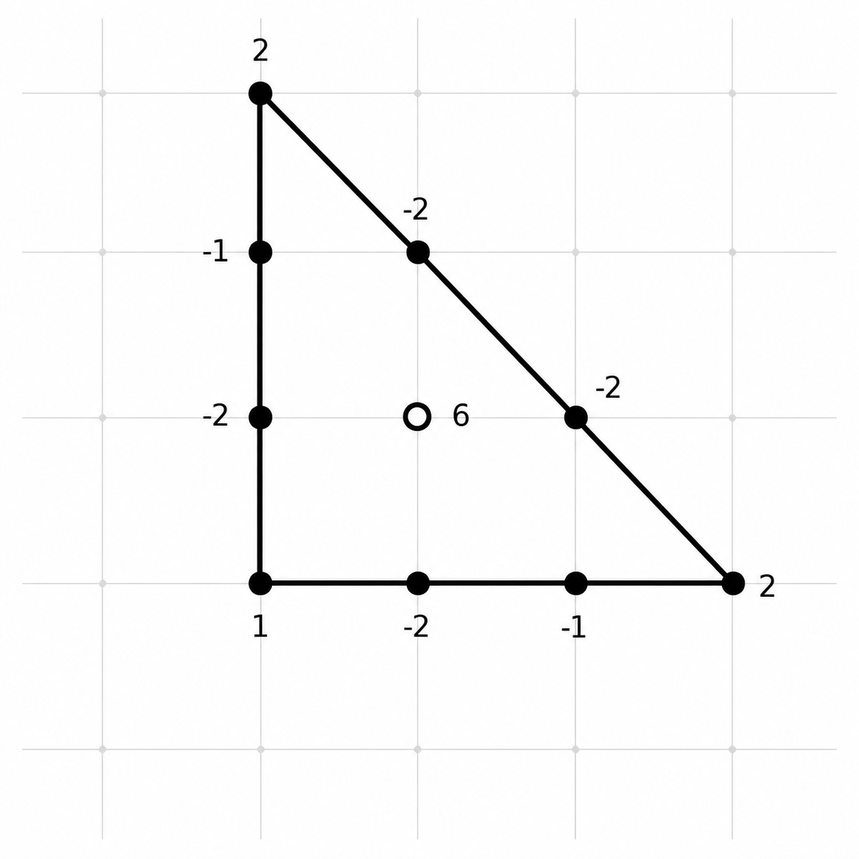}
		\caption{Newton polygon of $F_{\O(3)}^{2\underline{3}4III}$}
	\end{subfigure}
    \hfill
    \begin{subfigure}[t]{0.30\textwidth}
		\centering
		\includegraphics[width=\textwidth]{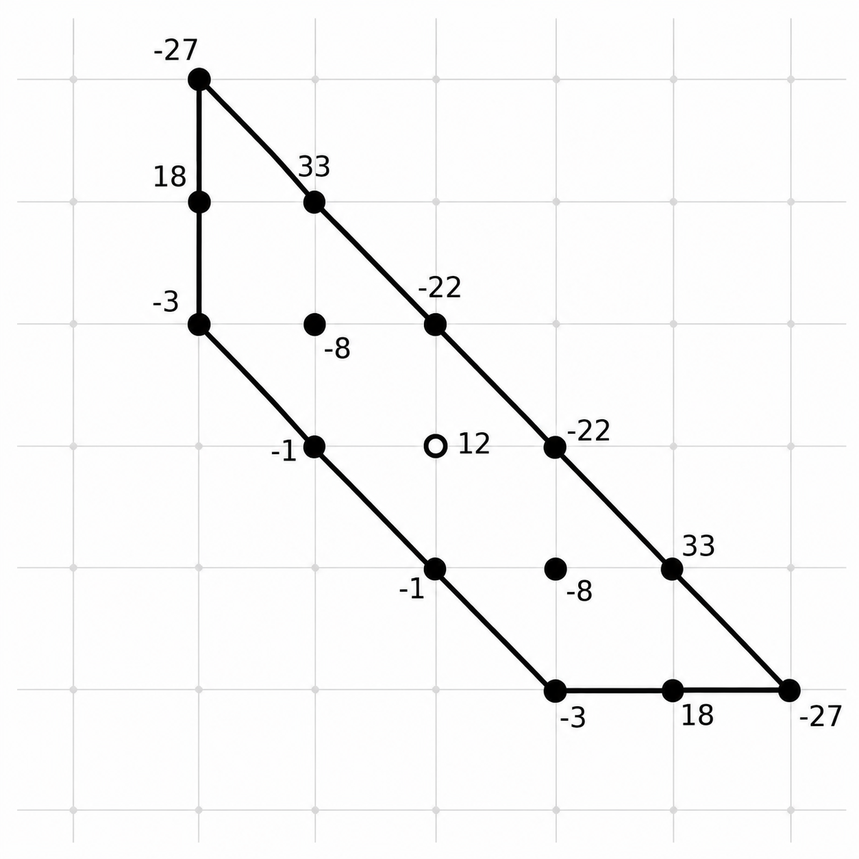}
		\caption{Newton polygon of $F_{\O(2)}^{\underline{2}34III}$}
	\end{subfigure}
}

\medskip
\noindent\makebox[\textwidth][c]{%
    \begin{subfigure}[t]{0.30\textwidth}
		\centering
		\includegraphics[width=\textwidth]{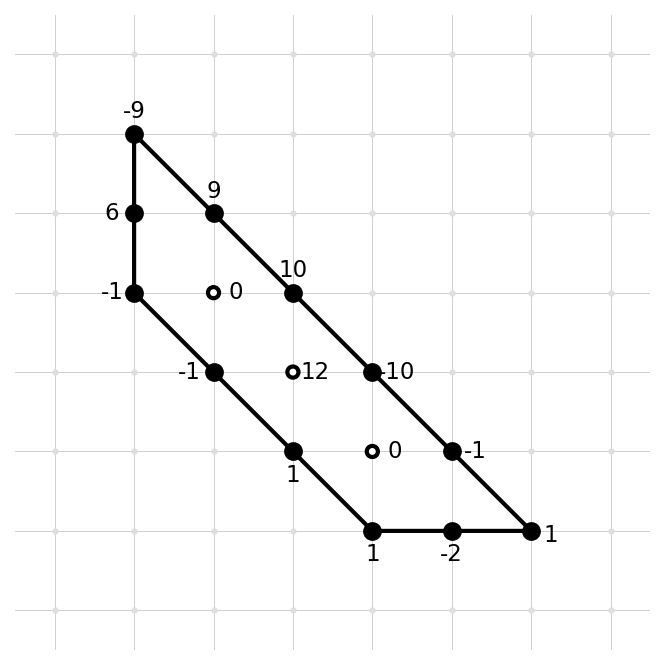}
		\caption{Newton polygon of $F_{\O(4)}^{1\underline{2}6III}$}
	\end{subfigure}
    \hfill
	\begin{subfigure}[t]{0.30\textwidth}
		\centering
		\includegraphics[width=\textwidth]{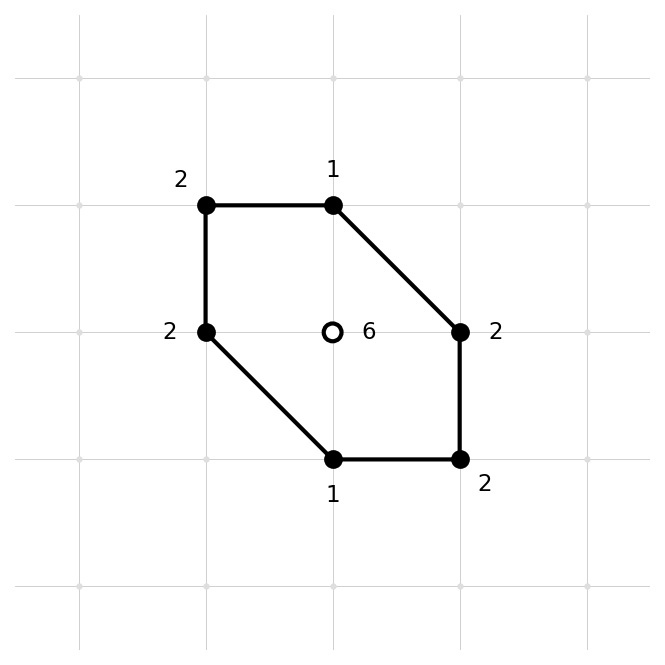}
		\caption{Newton polygon of $F_{\O(3)}^{12\underline{6}III}$}
	\end{subfigure}\hfill
	\begin{subfigure}[t]{0.30\textwidth}
		\centering
		\includegraphics[width=\textwidth]{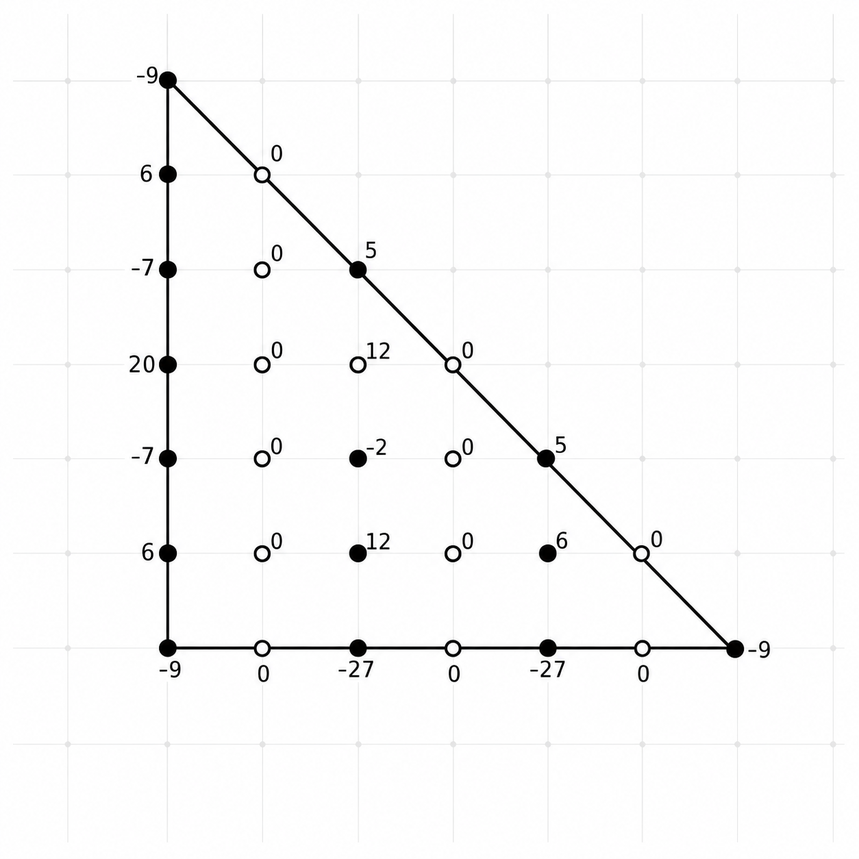}
		\caption{Newton polygon of $F_{\O(2)}^{\underline{1}26III}$}
	\end{subfigure}
    
    }
	\caption{Newton polygons of the Laurent polynomials for $\O$.}
	\label{fig:newton-polygons_O}
\end{figure}

\begin{figure}[htbp]
	\centering
	\captionsetup[subfigure]{labelformat=empty}
	\noindent\makebox[\textwidth][c]{%
	\begin{subfigure}[t]{0.30\textwidth}
		\centering
		\includegraphics[width=\textwidth]{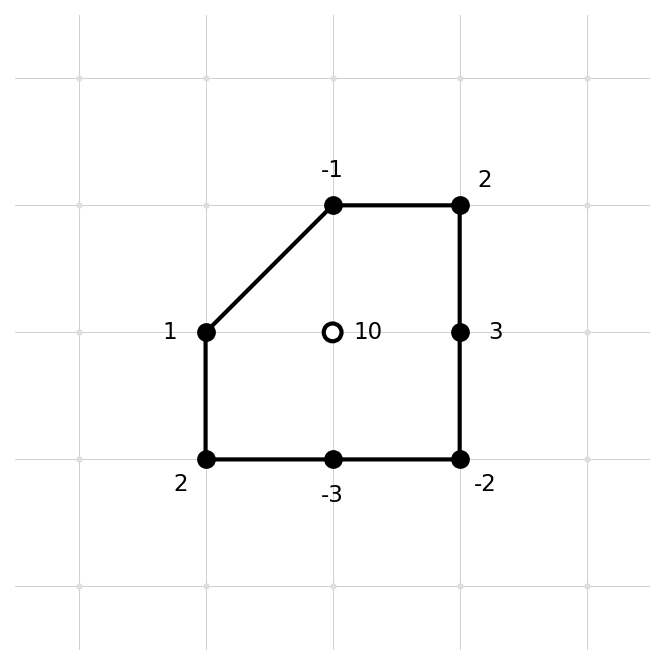}
		\caption{Newton polygon of $F_{\I(5)}$.}
	\end{subfigure}
    \hfill
	\begin{subfigure}[t]{0.30\textwidth}
		\centering
		\includegraphics[width=\textwidth]{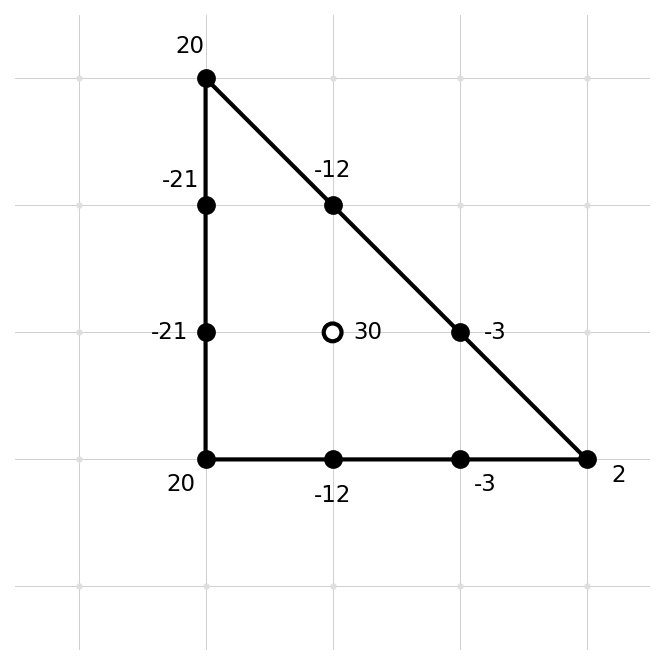}
		\caption{Newton polygon of $F_{\I(3)}$.}
    \end{subfigure}
        \hfill
	\begin{subfigure}[t]{0.30\textwidth}
		\centering
		\includegraphics[width=\textwidth]{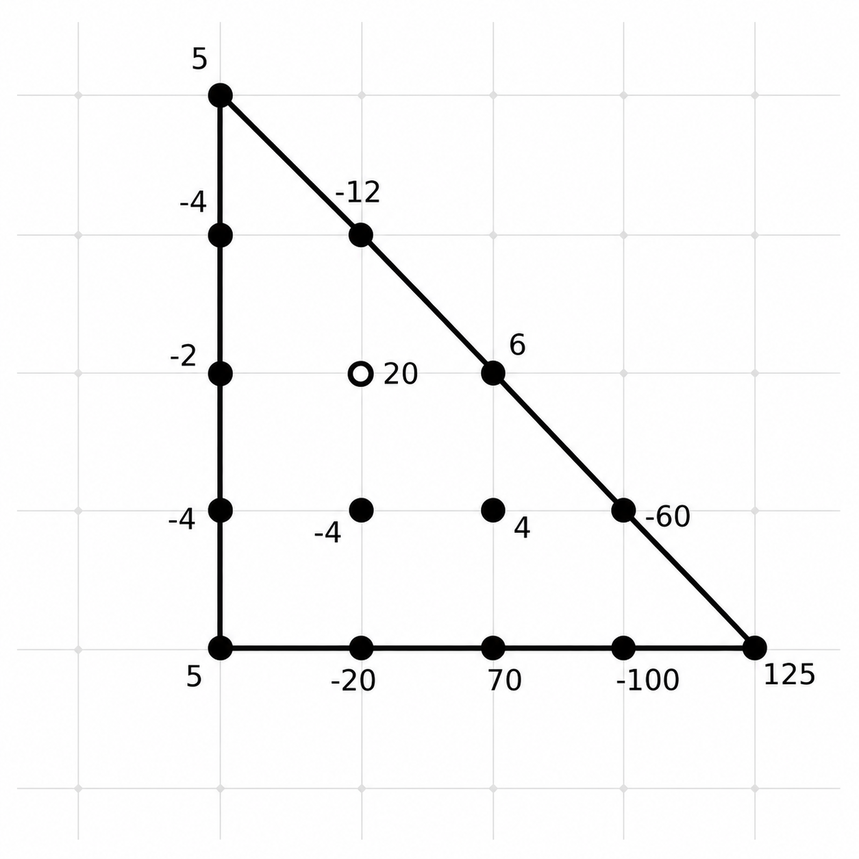}
		\caption{Newton polygon of $F_{\I(2)}$.}
	\end{subfigure}
     }
	\caption{Newton polygons of the Laurent polynomials for $\I$.}
	\label{fig:newton-polygons_I}
    
\end{figure}

The Figures \ref{fig:newton-polygons_T}, \ref{fig:newton-polygons_O} and \ref{fig:newton-polygons_I} show the corresponding Newton polygons where the arithmetic genus of each case can directly be read off by counting the interior points.

\newpage
For some of the sequences a nice binomial sum representation could be found, using the algorithm developed in \cite{Bostan}:

$\T(3):$ \[a_n=
\sum_{k=\lceil n/2\rceil}^{n}
(-27)^{\,n-k}2^{\,2k-n}
\binom{3k}{k}
\binom{2k}{k}
\binom{k}{n-k},
\]

$\T(2)$: \[a_{2m}=3^m\binom{2m}{m}\binom{3m}{m},
\quad
a_{2m+1}=0,\]

$\O(3)$: \[a_n=
\sum_{r=0}^{n}
\sum_{s=0}^{n-r}
4^r
\binom{n}{r}^2
\binom{n-r}{s}^2.\]
\\

A list of Laurent polynomials for the integer sequences arising from the cousins of the Platonic surfaces can be found in Appendix \ref{ap:herfurtnersurfaces} in Table \ref{tab:herfurtner-operators}. The existence of Laurent polynomial representations suggest that Lucas-type congruences should hold for those integer sequences.

\section{Some arithmetical aspects}\label{ch:somearithmeticalaspects}
In this section, we will develop the topic of integer sequences stemming from the expansion of holomorphic periods of our elliptic surfaces further. We will see how one can use
the $p$-adic properties of series expansions of the holomorphic and logarithmic solution to
determine a \emph{Frobenius structure}.

An elliptic surface over $\P^1$ determines an elliptic curve over the function field $\C(t)$.
Conversely, any  elliptic curve over $\C(t)$ has an elliptic surface as minimal model, so
that these two notions are basically the same. In all our examples it appeared that the equations we are
dealing with have coefficients in $\Q$, so that we have elliptic surfaces over $\Q$ and
corresponding elliptic curves over $\Q(t)$. By 'clearing out denominators' we in fact
can arrange having an elliptic curve defined over the $2$-dimensional ring $\Z[t]$. In
algebraic geometric terms we have a family
\[ \mathcal{E} \lra Spec(\Z[t])\]
whose general fibre is an elliptic curve, so $\mathcal{E}$ is an {\em arithmetic threefold}.
There are two sorts of specialisations: We can fix a value of $t$ and obtain, apart from
finitely many singular values, an elliptic curve over $Spec(\Z)$. Or we can fix a prime $p$,
reduce modulo $p$ and get an elliptic surface over $\F_p[t]$. When we fix both $t$ and $p$, we
get an elliptic curve $E_{t,p}$ over the finite field $\F_p$. The arithmetical properties of elliptic curves is a huge, beautiful and deep subject. We try to retrace some historical lines of discovery and will use only the most basic notions and refer to \cite{Roquette} for more details and the textbooks, e.g. \cite{Silverman}.

\subsection{Hasse, Witt and Deuring}
It was known after the works of E. Artin and F. K. Schmidt that the Zeta-function
of a curve $C$ of genus $g$ over a finite field $\F_q$, $q=p^f$, factors in a specific way:
\[Z(C;T):=\text{exp}\left(\sum_{k=1}^\infty \frac{\#C(\F_q)}{k}T^k\right)=\frac{P(T)}{(1-T)(1-qT)}.\]
It was conjectured that an analogue of the Riemann hypothesis holds true: the reciprocal roots of the polynomial $P$ all have size $\sqrt{q}$.
For an elliptic curve ($g=1$) one has
\[ P(T)=1-a_qT+qT^2\]
where $\#C(\F_q)=1-a_q+q$ and the conjecture is equivalent to $|a_q| \le 2\sqrt{q}$.
In a series of papers \cite{Hasse}, \cite{Hasse1-3}, Hasse developed the basic geometrical properties of an elliptic curve $E$ over a (perfect)
field $k$ of characteristic $p$, purely in the algebraic language of the function  field $K=k(E)$. He showed that the group
\[E[n]:=\{x \in E(\overline{k})\;|\;n\cdot x=0\}\]
of $n$-torsion points of $E$ (over the algebraic closure $\overline{k}$) 
has $n^2$ elements if $n$ is prime to $p$, $n$ elements if $n=p^{\nu}$ and a
certain quantity $A$ (now called the Hasse invariant) is non-zero, and
$1$ element if $n=p^{\nu}$ and $A=0$ (now called the super-singular case).
Furthermore, he showed that the endomorphism ring $End(E)$ contains $\Z$ as
a subring and hence is a ring of characteristic $0$. For each (algebraic) endomorphism 
$\mu$ outside the subring $\Z$, he showed that it is {\em imaginary quadratic},
i.e. satisfies an equation of the form
\[\mu^2+\nu \mu+m=0,\]
where $\nu, m \in \Z$ and $\nu^2 \le 4 m$. He applied this to the Frobenius
endomorphism $\pi \in End(E)$, induced by $f \mapsto f^q$ on the function field $K$, 
and concluded that the Riemann hypothesis holds for elliptic curves over a finite field, i.e.
$a_p^2 \le 4 p$.
We will call the number $a_p$ {\em the Frobenuis trace}, as it is interpreted nowadays as the trace of the map induced by Frobenius on any Weil cohomology theory.

Over an algebraically closed field $k$ the structure of the torsion points $E[n]$ of $E$ is remarkable: as long as $n$ is prime to $p$ the result is the same as in characteristic zero, forming the well-known pattern of division points in the fundamental period parallelogram of $E$. But if we specialise to characteristic $p$, the $p^2$-torsion points come together in groups of $p$ (if $A \neq 0$) to form the group of $p$-torsion points $E[p]$, which is of order $p$, and if $A=0$  all $p^2$ torsion points specialise to the neutral element, so $E[p]$ consists of one element.

If we divide out a point of order $p$ we obtain a separable, cyclic unramified field extension of degree $p$ of an elliptic function field and it is in this context that Hasse \cite{Hasse} introduced his invariant.

In \cite{HasseWitt} Hasse and Witt reformulated the definition of $A$ and generalised it to curves $C$ of arbitrary genus $g$ over $\F_p$. They start by choosing $g$ $k$-rational points $P_i$ on $C$, such that the divisor $D=P_1+P_2+\ldots+P_g$ is non-special. For each point they pick a local uniformiser $t_i$
at $P_i$. They consider the space $P(*D) \supset P(D)$ of principal parts supported at $D$ and the subgroup of those elements of pole order $\le 1$.
Furthermore, one considers $K(*D) \supset K(D)$ of elements of the function field with poles in $D$ and the subspace of those elements with pole of order $\le 1$.
The factor space $V:=P(*D)/K(*D) \simeq P(D)/K(D)$ is of dimension $g$, with
the elements $r_i$, defined by $(r_i)_Q=[1/t_i]$ if $Q=P_i$ and $0$ else, as basis.
Applying the Frobenius map $f \mapsto f^p$ to $r_i$ and reexpressing the result
in terms of the basis $r_i$ leads to a matrix $A=(A_{ij})$:
\[r_i^p=\sum_{j=1}^g A_{ij} r_j \in V. \]
The matrix $A$ is of size $g \times g$ with entries in $k$, nowadays called \emph{Hasse-Witt matrix}. It is an object of semi-linear algebra, as under a base change matrix $S$, the matrix $A$ transforms to $S^{-1}AS^{(p)}$, where $S^{(p)}$ denotes the matrix obtained from $S$ by taking the $p$-th power of each entry. For $g=1$ the matrix $A$ is a single number and scaling the basis by a factor $r$ 
changes $A$ to $r^{-1}A r^p=r^{p-1}A$. If $k$ is algebraically closed we can scale $A$ to $0$ or $1$; 
for the prime field $k=\F_p$, $A$ is well-defined as element of $k$.
The procedure of Hasse and Witt can be interpreted in newspeak by remarking that
the vector space $V=H^1(C,\mathcal{O}_C)$ and $A$ is identified as matrix of the
natural $p$-linear Frobenius map $H^1(C,\mathcal{O}_C) \lra H^1(C,\mathcal{O}_C)$ induced by $\pi$, with respect to a basis defined by the divisor $D$.\footnote{The exact sequence of sheaves on $C$ ($\mathcal {O}:=\mathcal{O}_C$): $\mathcal{O} \hookrightarrow \mathcal{O}(*D)\twoheadrightarrow \mathcal{O}(*D)/\mathcal{O}$ gives, by taking global sections, 
$H^1(\mathcal{O}) =Coker(H^0(\mathcal{O}(*D))\lra H^0(\mathcal{O}(*D)/\mathcal{O})=P(*D)/K(*D))=V.$
Similarly one sees from the exact sequence $\mathcal{O} \hookrightarrow \mathcal{O}(D) \twoheadrightarrow \mathcal{O}(D)/\mathcal{O}$ and Riemann-Roch that if $D$ is
non-special one has $P(D)/K(D)=H^1(\mathcal{O})=V$.}

Deuring classified the possible endomorphism rings of elliptic curves in  \cite{Deuring} and showed that $A=0$
is equivalent to $End(E)$ being an order in a quaternion algebra. 
He computed $A$ for an elliptic curve explicitly with respect to the point at infinity and the uniformiser $x/y$, using the definition in \cite{HasseWitt}. If the elliptic curve is given in the form
\begin{equation}
	y^2=f(x)=x^3+a_2x^2+a_4x+a_6,
\end{equation}
the invariant $A$ can be computed as the coefficient of $x^{p-1}$ in $f^\frac{p-1}{2}$. If $E$ is given in Legendre normal form $y^2=x(x-1)(x-t)$, the Hasse invariant $A$ depends polynomially on the
parameter $t$ and Deuring finds (for $p\ge 3)$:
\begin{equation}\label{eq:hasselegendre}
	A=A(t)=(-1)^{\frac{p-1}{2}}\sum_{n=0}^\frac{p-1}{2}\binom{\frac{p-1}{2}}{n}^2 t^n. 
\end{equation}

\subsection{Igusa and Manin}
We have seen that the  Picard-Fuchs operator of the Legendre family $y^2=x(x-1)(x-t)$ is
\begin{equation}
	\mathcal{L}:=t(1-t)\frac{d^2}{dt^2}+(1-2t)\frac{d}{dt}-\frac{1}{4},
\end{equation}
The unique holomorphic solution to $\mathcal{L}\phi=0$ near $0$ is given by the power series 
\begin{equation}\label{eq:legendresol}
	\phi(t)=\sum_{n=0}^\infty\binom{\frac{-1}{2}}{n}^2t^n,
\end{equation}
that is, 
\[\phi(t)=1+\left(\frac{1}{2}\right)^2 t+\left(\frac{1\cdot 3}{2\cdot 4}\right)^2 t^2+\left(\frac{1\cdot3 \cdot 5}{2\cdot 4 \cdot 6}\right)^2 t^3+\ldots=\sum_{n=0}^{\infty} \binom{2n}{n}^2 \left(\frac{t}{16}\right)^n.\]

Igusa notices in \cite{Igusa} that the expression \eqref{eq:hasselegendre} coincides with the $(p-1)$-truncation
\[\phi^{(p-1)}(t):=\sum_{n=0}^{p-1}\binom{\frac{-1}{2}}{n}^2t^n\] 
of \eqref{eq:legendresol} taken modulo $p$:
\[ (-1)^{(p-1)/2}\phi^{(p-1)}(t)=A(t) \bmod p.\]
As one can show
\[ \phi(t)=\phi^{(p-1)}(t)\phi(t^p) \bmod p,\]
we see that $A(t)$ solves the Picard-Fuchs equation $\bmod p$:
\[\mathcal{L} A(t) = 0 \bmod p,\]
and that as a result of this fact, the roots of the Hasse polynomial $A(t)$ are all distinct.

By an elementary argument, Deuring's expression can also be seen to arise in the determination of the number of points of the Legendre curve $E_t$ considered as family over $\F_p$, taken $\bmod p$: 
\[ \#E_t(\F_p)=1-A(t) \bmod p.\]
Hence we see that 
\[ a_p(t)=A(t) \bmod p.\]
In \cite{Manin} Manin showed this directly for arbitrary elliptic curves over $k=\F_p$ and showed
\begin{equation}
	P(T)\equiv 1-AT \bmod p,
\end{equation}
where $P(T)=1-a_pT+pT^2$ denotes as before the numerator polynomial of the Zeta-function of the curve.
It seems strange that  Hasse and Deuring have overlooked this fundamental fact or found it too trivial to mention. Keeping in mind that $A$ is the trace of Frobenius on $H^1(E,\mathcal{O}_E)$, this statement is now an easy consequence  of the general Fulton trace formula \cite{Fulton}:
\[\#X(\F_p) \equiv\sum_{i=0}^d (-1)^i Tr(F | H^i(X,\mathcal{O}_X)) \bmod p,\] which was not available at the time. 

Manin also showed in \cite{Manin} that Igusa's observation was not a coincidence for the Legendre family, but in fact is a general phenomenon. For that, he developed in great generality the theory of differentiation of deRham cohomology classes with respect to a parameter, now called the {\em Gauss-Manin connection}.

We consider  curves $C$ to be  over the field $K:=k(t)$ of rational functions with coefficients in $k=\F_p$. 
The deRham cohomology $H^1_{dR}(C)$ sits in a short
exact sequence
\[ 0 \lra H^0(C,\Omega^1_C) \lra H^1_{dR}(C) \lra H^1(C,\mathcal{O}_C) \lra 0,\]
and the spaces $H^1(C,\mathcal{O}_C)$ and $H^0(C,\Omega^1_{C})$ are $g$-dimensional $K$-vector spaces, dually paired by an anti-symmetric pairing
\[ \langle -,- \rangle: H^1_{dR}(C) \times H^1_{dR}(C) \lra K,\]
which extends the Serre-duality pairing between 
$H^0(C,\Omega^1_C)$  and $H^1(C,\mathcal{O}_C)$.

The space $Der_k(k(t), k(t))$ of
derivations is generated by  $\partial_t=\frac{\partial}{\partial t}$
and acts on the first deRham cohomology $H^1_{dR}(C)$:
\[\nabla=\nabla_{\partial_t}: H^1_{dR}(C) \lra H^1_{dR}(C).\]
By composition
of derivations it extends to an action of the universal enveloping algebra, i.e. the ring $D:=k[t,\partial_t]$ of differential operators. 

If $C=E$ is an elliptic curve and if $r$ is a generator of $H^1(E,\mathcal{O}_E)$ and $\omega\in H^0(E,\Omega^1_{E})$ its dual, then for the Frobenius $F:H^1(E,\mathcal{O}_E)\rightarrow H^1(E,\mathcal{O}_E)$ we have
\begin{equation}
	A=\langle \omega,F(r)\rangle \in k(t).
\end{equation}

For $\omega\in H^0(E,\Omega^1_{E})$ we let  $\overline{\omega}$ denote its image in  $H^1_{dR}(E)$.

{\bf Theorem} (Manin) If for some $\mathcal{L} \in D$ one has $\mathcal{L}(\overline{\omega})=0$, then $\mathcal{L}(A)=0$.

{\bf proof sketch:} $\nabla_{\mathcal{L}}(\overline{\omega})=0$ means that we have
\begin{equation}
	\nabla_{\mathcal{L}}\omega=d\phi,
\end{equation}
for some function $\phi\in K(E)$, so that
\begin{equation}
	\langle\nabla_\mathcal{L}(\omega),F(r)\rangle = \langle d\phi,F(r)\rangle.
\end{equation}
In \cite{Cartier}, Cartier had defined an operator $C$, with the property that $C(d\phi)=0$ and that it is adjoint to $F$ in the sense that $\langle \omega,F(r)\rangle=\langle C(\omega),r\rangle ^p$. It follows
\begin{equation}
	\langle C(\nabla_\mathcal{L}(\omega)),r\rangle^p=0.
\end{equation}
Manin's Cartier compatibility shows that applying $\mathcal{L}$ to the Hasse--Witt entry is the same as pairing the Cartier transform of $\nabla_{\mathcal L}(\omega)$ with $r$:
\begin{equation} 0=\langle C(\nabla_\mathcal{L}(\omega)),r\rangle^p=\mathcal{L}(\langle C(\omega),r\rangle^p)=\mathcal{L}(\langle\omega,F(r)\rangle)=\mathcal{L}(A).
\end{equation}
Therefore the Hasse invariant $A$ satisfies the same differential equation as the deRham class
of the differential form $\omega$.
\qed

We refer to the great book \cite{Clemens} for a lucid account of this important chain of ideas.

We recall that for an elliptic surface $\pi: \mathcal{E} \to \P^1$ there is a naturally attached {\em fundamental line bundle $\mathcal{L}$} (\cite{Miranda}, section II.4). It can be defined as the Hodge bundle $\pi_*(\omega_{\mathcal{E}/\P^1})$, with fibres $H^0(\omega_{E_t})$, which is
dual to $R^1\pi_*(\mathcal{O}_{\mathcal{E}})$, with fibres $H^1(\mathcal{O}_{E_t})$. The degree of this bundle relates to the geometric genus of the elliptic surface (\cite{Miranda}, section IV.1): 
\[\deg(\mathcal{L}) =p_g(\mathcal{E})+1.\]
For rational elliptic surfaces we have $p_g=0$, so $\deg(\mathcal{L}) =1$.
Frobenius induces a $p$-linear map 
from $\mathcal{L}$ to itself, and as a result the Hasse-invariant can be seen as a section of $\mathcal{L}^{p-1} \simeq \mathcal{O}_{\P^1}(p-1)$,
i.e. a polynomial of degree $\le p-1$.

According to Manin's theorem, the polynomial $A(t)$
is a solution of the Picard-Fuchs equation $\bmod p$.
In case the fibre over $0$ is semi-stable, the Picard-Fuchs operator has a MUM-point at $0$ and 
has a normalised period function $\phi(t)=1+u_1t+\ldots \in\F_p[[t]]$. This solution is unique up to
multiplication by a series  $Q(t^p) \in \F_p[[t^p]]$ (which has derivative zero in characteristic $p$).
It follows that one must have 
\[ A(t) =Q(t^p)\phi(t).\]
If we truncate this relation at $t^p$ we find:
\[A(t)=A(t)^{(p-1)}=Q(0) \phi^{(p-1)}(t).\]
Note that $Q(0)=\pm 1$, as the Hasse invariant $A(0)$
of the fibre at $0$ is $\pm 1$ (\cite{Silverman}, p. $449$) and by construction $\phi(0)=1$.

{\bf Corollary} {\em Let $\mathcal{E}\rightarrow\P^1$ be a rational elliptic surface over $\Z$ with a semistable fibre at $0$. Then, at each smooth fibre $t$, if the reduction $E_t$ mod $p$ is smooth, we have the trace formula
\begin{equation}
	A(t)=a_p(t) \equiv \pm \phi(t)^{(p-1)}\bmod p.
\end{equation}
}
So this statement allows us to determine $\pm a_p(t) \bmod p$ from the Picard-Fuchs equation alone.
We remark that the sign $\pm$ is usually dependent on $p$. For the Legendre family it is $(-1)^{(p-1)/2}$,
determined by split/non-split multiplicative reduction at $t=0$. It can not be determined from the Picard-Fuchs operator alone. If we perform a quadratic twist on the fibres $E_t/\F_p$, then $a_p(t)$ changes by a sign of the quadratic character, but the Picard-Fuchs operator does not change.

\subsection{The Dwork expansion}
We have seen that for an elliptic curve the Frobenius trace $a_p \bmod p$ can be computed by the
Hasse invariant and that in a family it satisfies the Picard-Fuchs equation $\bmod p$. It is natural to try to obtain $a_p$ modulo higher powers of $p$. Can one use the complete series $\phi(t)$ to 
obtain congruences for $a_p \bmod p^2, \bmod p^3, \ldots$? As we will see, there are indeed such formulas, which find their natural place in $p$-adic analysis and $p$-adic cohomology theories.
\footnote{Of course,
from $|a_p| \le 2\sqrt{p}$, since $2\sqrt{p}<p/2$ for $p>13$, the residue class mod $p$ uniquely determines
the integer $a_p$, so from a practical standpoint this is not very urgent. But for curves of higher genus or higher-dimensional varieties and from a theoretical standpoint it turns out to be very useful.}

\subsubsection{Unit root formula}  
The numerator polynomial $P(T)=1-a_pT+pT^2 \in \Z[T]$ of an elliptic curve $E$ over $\F_p$ in general does not factor over $\Z$. If we have a factorisation over the $p$-adic numbers $\Z_p$ as
\[ P(T)=(1-uT)(1-vT) \in \Z_p[T],\]
then the $p$-adic order of say $u$ must be $0$, that of $v$ must be one. One calls $u \in \Z_p$ the {\em unit root} of $E$; then $v=p/u$. 
There  exists a unit root precisely when the Hasse invariant $A \neq 0$.
For the Legendre family $E_t:\; y^2=x(x-1)(x-t)$ of elliptic curves, viewed as family over $\F_p[t]$,
Dwork showed the following.
 
{\bf Theorem } (Dwork's Unit root formula, \cite{Dwork1}) 
Let $p \neq 2$ and $t_0\in \F_p$ such that $E_{t_0}$ is smooth ($\Delta(t_0) \neq 0$) and not supersingular ($A(t_0) \neq 0$). Then the unit root is given by
\begin{equation}\label{eq:unitroot}
	u_{t_0}=(-1)^{(p-1)/2}\frac{\phi(z)}{\phi(z^p)}\mid_{z=[t_0]},
\end{equation}
where  $\phi(z)\in\Z_p[[z]]$ is as in \eqref{eq:legendresol} and $[t_0]\in\Z_p$ the multiplicative lift of $t_0$. 
(Implicit in this statement is that the quotient $\phi(z)/\phi(z^p) \in \Z_p[[z]]$ allows analytic continuation to the point $[t_0]$).

Dwork showed for general $\phi(z)=\sum_{n=0}^\infty a(n)z^n$ and the corresponding sequence $1,a(1),a(2),\dots$ that if for all $m,s\in\N$ the \emph{Dwork congruences}
\begin{equation}\label{eq:dworkcongruences}
a(n)a\left(\left\lfloor\frac{n+mp^s}{p}\right\rfloor\right)\equiv a(n+mp^s)a\left(\left\lfloor\frac{n}{p}\right\rfloor\right)\bmod p^s
\end{equation}
are satisfied, then for all $s\in\N$ the unit root can be computed $\bmod p^s$:
\begin{equation}\label{eq:unitrootapprox}
\frac{\phi(z)}{\phi(z^p)}\mid_{z=[t_0]}\equiv \frac{\phi^{(p^s-1)}([t_0])}{\phi^{(p^{s-1}-1)}([t_0]^p)}\bmod p^s.
\end{equation}
This makes it possible to effectively compute the unit root, but for larger values of $p$ and higher powers of $s$ it becomes prohibitive.

{\bf Example} We take the Legendre family $y^2=x(x-1)(x-t)$ defined over $\F_5$ and consider the fibre at $t_0=4\bmod 5$ to illustrate \eqref{eq:unitroot}. The discriminant of the family is $\Delta(t)=16t^2(1-t)^2$, hence $\Delta(4)\neq 0\bmod 5$ and the fibre $E_{4}$ is smooth. Furthermore $A(4)\equiv\sum_{n=0}^{\frac{5-1}{2}}\binom{\frac{-1}{2}}{n}^24^n=1+\left(\frac{1}{2}\right)^24+\left(\frac{1\cdot3}{2\cdot4}\right)^24^2\equiv 3\not\equiv 0\bmod 5$, so $E_4$ is not supersingular and the conditions of the theorem are satisfied. Since the sequence given by \eqref{eq:legendresol} satisfies the Dwork congruences, we can compute the $p$-adic digits by using \eqref{eq:unitrootapprox}:
\begin{align*}
s&=1: \quad u_{4}\equiv 3 \bmod 5, \\
s&=2: \quad u_{4}\equiv 3+2\cdot5 \bmod 5^2, \\
s&=3: \quad u_{4}\equiv 3+2\cdot5+4\cdot5^2 \bmod 5^3, \\
s&=4: \quad u_{4}\equiv 3+2\cdot5+4\cdot5^2+2\cdot5^3 \bmod 5^4, \\
&\dots
\end{align*}
We see that the approximation tends to a value $u_4$ satisfying $u_{4}+p/u_{4}=-2$, which is indeed the correct value for $a_p(4)$, which can easily be verified  by a direct point count.

\subsubsection{Deformation method}
It is a strange thing to try to compute a single root of a polynomial of higher degree, and use only one of the solutions to the differential equation, and having to exclude the super-singular cases. It would be much more natural to look for a matricial version of the unit root formula that uses a basis of solutions. We
will describe such a method, discovered by Dwork that he
described in \cite{Dwork1} and \cite{Dwork2}.
He shows that the Zeta-functions of a general hypersurface varying in a $1$-parameter family (defined say over $\Z$) can be obtained in terms of the Frobenius endomorphism on the solution space of the Picard-Fuchs equation, if the Frobenius morphism is 'known'  at a single point.

The point is that the Frobenius action $F$ and Gauss-Manin differentiation $\nabla$ basically 'commute'. Imagine we have a family $\mathcal{E} \lra \P^1$ and pick a prime $p$ of good reduction;
by completion we obtain a family $\mathcal{E}_{\Z_p}$ over $\P^1_{\Z_p}$.
One can construct the relative deRham cohomology ${\bf H} =H^1_{dR}(\mathcal{E}_{\Z_p})$ of this family.
By restriction to the affine line we obtain a free $R$-module of rank $2$,
where $R:=\Z_p[[t]]$. As before, there is a duality
pairing 
\[\langle-,-\rangle : {\bf H} \times {\bf H} \lra R,\]
extending the pairing we had $\bmod p$.
On the ring $R$ there are two operations
\[\theta: R \lra R, \; a(t) \mapsto t\frac{da(t)}{dt}\;\text{ and }\;\sigma: R \lra R, \; a(t) \mapsto a(t^p),\]
and it holds, as operation on $R$, the commutation
\[\theta \sigma = p \sigma \theta.\]
These operations lift to operations
\[\nabla: {\bf H}\lra {\bf H},\;\;F:{\bf H}\lra {\bf H}.\]
Here $\nabla=\nabla_{\theta}$, the Frobenius lift $F$ is $\sigma$-linear, and they satisfy the commutation
\[ \nabla F =pF \nabla.\]
After choosing a basis for ${\bf H}$, this translates into a differential equation for the matrix $F=F(t)$.
Namely, if $B$ denotes the connection matrix of $\nabla$ in the given basis, the differential equation reads $ \theta(F)=pF(t)B(t^p)-B(t)F(t)$.
Dwork showed that there exists a constant matrix $F(0)$ such that, if $Y(t)$ denotes a fundamental solution matrix for
$\nabla$, the differential equation is solved by
\begin{equation}
F(t) =Y(t)^{-1} F(0) Y(t^p).
\end{equation}
In this way we get, starting from $F(0)$, a power series expansion of $F(t)$ in powers of $t$, in terms of data uniquely attached to the Gauss-Manin connection, i.e. the Picard-Fuchs operator. Dwork showed that if $t_0 \in \P^1(\F_p)$
is a regular fibre, then the series $F(t)$ admits $p$-adic analytic continuation to the multiplicative lift $[t_0]$ where it can be evaluated and that
$F([t_0])$ is the matrix of Frobenius of the fibre
at $t_0$, in particular $a_p(t_0)=tr(F([t_0])) \in \Z_p.$

For technical reasons we use $U(t)=p F(t)^{-1}$ that satisfies 
\begin{equation}\label{eq:umatrix}
U(t) =Y(t^p)^{-1} U(0) Y(t).
\end{equation}
Then the formula is
\[a_p(t_0)=tr(U([t_0])) \in \Z_p.\]
In \cite{Dwork2} also the case that $0$ is a singular
point of the differential equation is covered.

We spell this out for the Picard-Fuchs equation of elliptic surfaces with semi-stable fibre at $0$. First we redefine the discriminant of the operator:
\begin{equation}
\Delta(t)=c(t-s_1)\dots(t-s_r), 
\end{equation}
$s_i$ the singular points $\neq 0,\infty$ of the operator, $c$ minimal, such that $\Delta(t)$ has integral coefficients, $c\cdot s_1\cdot\dots \cdot s_r>0$.
In the following $p$ is an arbitrary prime not dividing the polynomial discriminant of $\Delta(t)$
and not dividing the numerators and denominators of the local exponents $\alpha$ and $\beta$ at infinity. 

We recall the Frobenius basis of solutions near $0$:
\begin{equation}\label{eq:frobbasis}
\begin{aligned}
	\phi(t)&=f_0(t), \\
	\psi(t)&=\log(t)\cdot f_0(t)+f_1(t),
\end{aligned}
\end{equation}
with $f_0(t)\in\Z[[t]]$, $f_1(t)\in t\Q[[t]]$. The fundamental solution matrix $Y(t)$ is:
\begin{equation}
	Y(t)=\begin{pmatrix}
		\phi(t) & \theta (\phi(t)) \\
		\psi(t) & \theta (\psi(t)) 
	\end{pmatrix}
=\begin{pmatrix}
	1 & 0 \\
	\log(t) & 1 
\end{pmatrix}
\begin{pmatrix}
	f_0(t) & \theta (f_0(t)) \\
	f_1(t) & f_0(t)+\theta (f_1(t)) 
\end{pmatrix}.
\end{equation}

In this case the matrix $U(0)$ in \eqref{eq:umatrix} has to be lower triangular:
\begin{equation}\label{eq:umatrixshape}
	U(t)=Y(t^p)^{-1}\begin{pmatrix}
		\epsilon & 0 \\
		b & \epsilon p 
	\end{pmatrix} Y(t),
\end{equation}
with $\epsilon=\pm1$ and $b$ some constant, which needs to be specified.
One can then write
\begin{equation}
	\begin{aligned}
		 \Delta(t)^{p\cdot\delta_{0}}\cdot U(t) & = U_0(t) \bmod p, \\
		 \Delta(t)^{p\cdot\delta_{1}}\cdot U(t) & = U_1(t) \bmod p^2, \\
		\Delta(t)^{p\cdot\delta_{2}}\cdot U(t) &= U_2(t) \bmod p^3, \\
        \Delta(t)^{p\cdot\delta_{3}}\cdot U(t) &= U_3(t) \bmod p^4, \quad \dots
	\end{aligned}
\end{equation}
for some matrices $U_N\in Mat(2\times2,(\Z/p^{N+1})[t])$. Hence one gets a $p$-adic expansion for $U(t)$ that we write as:

\begin{equation}\label{eq:umatrixisrational}
	U(t) = \frac{V_0(t)}{\Delta(t)^{p\cdot\delta_0}} + p\cdot\frac{V_1(t)}{\Delta(t)^{p\cdot\delta_1}} + p^2\cdot \frac{V_2(t)}{\Delta(t)^{p\cdot\delta_2}} + p^3 \cdot \frac{V_3(t)}{\Delta(t)^{p\cdot\delta_3}} + \dots,
\end{equation}

with $V_N\in Mat(2\times2,\{0,1,\ldots,p-1\}[t])$. 

\subsection{Hypergeometric and Beauville Examples}
First we consider the four well-known hypergeometric families with semi-stable fibre at $0$ (Table \ref{tab:euclidean-elliptic-surfaces}), among which are the Legendre family and the Hesse pencil, and the six Beauville families of elliptic curves (Table \ref{tab:beauville-elliptic-surfaces}). For each of those families the Picard-Fuchs equation has a  Frobenius basis of solutions near $0$.

Dwork analyses the Legendre family and differential operator \eqref{eq:legendrethetaform} in \cite{Dwork2} in great detail. He proves that in that concrete case we have  in \eqref{eq:umatrixshape}
\begin{equation}
b=(p-1)\log_p(16)
\end{equation}
and $\epsilon=(-1)^{\frac{p-1}{2}}$. Looking at the holomorphic solution \eqref{eq:legendresol} we see that the operator is \emph{not critically scaled}. Critical scaling is obtained by rescaling with the factor $16$, which is exactly the constant appearing in $b$, which leads to the hypergeometric operator
\begin{equation}\label{eq:legendrecriticallyscaled}
	\theta^2-16t\left(\theta+\frac{1}{2}\right)^2.
\end{equation}
It turns out experimentally that, assuming all the differential operators of the above mentioned cases are critically scaled, the correct value\footnote{The differential equation is formally solved by $U(t)$ with arbitrary values of $b$ but then the power series entries do not admit global $p$-adic analytic continuation. Experimentally the correct value of $b$ can be determined by checking if $U(t)$ can be evaluated at every $[t_0]$.} in \eqref{eq:umatrixshape} is $b=0$:
\begin{equation}\label{eq:U0naiv}
	U(0)=\begin{pmatrix}
		1 & 0 \\
		0 & p 
	\end{pmatrix}.
\end{equation}
Furthermore, it is experimentally possible to determine the degrees of $V_N$ (meaning the maximal degree of the polynomial entries) and the powers $\delta_N$ in \eqref{eq:umatrixisrational}. In all our examples we have
\begin{equation}\label{eq:degreeofVN}
\deg(V_0)=\lfloor(p-1)\alpha\rfloor, \quad \deg(V_N)=\lfloor(p-1)\beta\rfloor+p\cdot\delta_N\cdot\deg(\Delta), \; N\geq1,
\end{equation}
where the ordering of the local exponents at infinity $(\alpha,\beta)$ is given by $\alpha\leq\beta$.

The $\delta_N \in \N$ are the \emph{minimal exponents} such that
\begin{equation}\label{eq:wantthistobepol}
	\Delta(t)^{p\cdot\delta_{N}} \cdot U(t) \equiv U_N(t) \bmod p^{N+1}
\end{equation}
is a matrix with polynomial entries. Their determination led to interesting observations, which are demonstrated in the example of the critically scaled Legendre family.

We consider the operator \eqref{eq:legendrecriticallyscaled}, which has discriminant $\Delta(t)=1-16t$. For $p=3$ we get the following concrete expansion for \eqref{eq:umatrixisrational}:
\begin{equation}\label{eq:Uexpansionconcreteexample}
\begin{aligned}
	&U(t)=\begin{pmatrix}
		1+t & 4t \\
		0 & 0 
	\end{pmatrix}
    +p\cdot\begin{pmatrix}
		t & t \\
		2t & 1 
	\end{pmatrix} \\[3mm]
    &+p^2\frac{\begin{pmatrix}
t^2 + 2t^3 + t^4 + 2t^5 + 2t^6
&
2t^2 + t^5 + t^7
\\[2mm]
2t + t^2 + t^3 + 2t^4 + 2t^5
&
t + 2t^2 + t^5
\end{pmatrix}
    }{(1-16t)^{3\cdot2}}\\[3mm]
    &+p^3
\frac{\begin{pmatrix}
t^2 + t^4 + t^6 + 2t^7 + 2t^8 + 2t^9 + 2t^{10}
&
2t^2 + t^3 + t^4 + 2t^5 + 2t^6 + 2t^7 + 2t^9 + t^{10}
\\[2mm]
2t^2 + t^5 + 2t^6 + t^7 + t^9 + t^{10}
&
t + 2t^2 + t^3 + t^4 + t^5 + t^7 + t^8 + t^9 + t^{10}
\end{pmatrix}}{(1-16t)^{3\cdot3}}\\
    &+\dots
\end{aligned}
\end{equation}

\newpage

Looking at the $\delta_N$-s more carefully shows an interesting pattern which is summarised in the following table:
\[
\renewcommand{\arraystretch}{1.5}
\begin{array}{c|cccccccccccccccccccccc}
    \hline
    N
    & 0 & 1 & 2
    & \color{Green}3\normalcolor & 4 & 5
    & 6 & 7 & 8 & \color{Green}9\normalcolor & 10 & 11 & \dots 
    \\[1mm]
    \hline
    N-\delta_N
    & 0 & 1
    & \color{Green}0\normalcolor & \color{Green}0\normalcolor & 1
    & 0 & 0 & \color{Green}-1\normalcolor & \color{Green}-1\normalcolor & \color{Green}0\normalcolor & 1 & 0 & \dots 
    \\[1mm]
    \hline
    \hline
    N
    & 18 & 19 & 20 & 21 & 22 & 23 & 24 & 25 & 26 & \color{Green}27\normalcolor & 28 & 29 & \dots 
    \\[1mm]
    \hline
    N-\delta_N
    & 0 & 1 & 0 & 0 & 1 & 0 & \color{Green}-2\normalcolor & \color{Green}-2\normalcolor
    & \color{Green}-1\normalcolor & \color{Green}0\normalcolor & 1
    & 0 & \dots 
    \\[1mm]
    \hline
    \hline
    N
    & 72 & 73 & 74 & 75 & 76 & 77 & 78 & 79 & 80 & \color{Green}81\normalcolor & 82 & 83 & \dots 
    \\[1mm]
    \hline
    N-\delta_N
    & 0 & 1 & 0 & 0 & 1 & \color{Green}-3\normalcolor & \color{Green}-3\normalcolor & \color{Green}-2\normalcolor
    & \color{Green}-1\normalcolor & \color{Green}0\normalcolor & 1
    & 0 & \dots 
    \\[1mm]
    \hline
    \hline
    N
    & 234 & 235 & 236 & 237 & 238 & 239 & 240 & 241 & 242 & \color{Green}243\normalcolor & 244 & 245 & \dots 
    \\[1mm]
    \hline
    N-\delta_N
    & 0 & 1 & 0 & 0 & \color{Green}-4\normalcolor & \color{Green}-4\normalcolor & \color{Green}-3\normalcolor & \color{Green}-2\normalcolor
    & \color{Green}-1\normalcolor & \color{Green}0\normalcolor & 1
    & 0 & \dots 
    \\[1mm]
    \hline
\end{array}
\]
\normalsize

This pattern can be observed for each prime $p\neq2$. For $p=2$ we have a much simpler rule for the $\delta_N$-s: The matrix $U(t)\equiv U_N(t)\bmod p^{N+1}$ is already polynomial for each $N$, i.e. $\delta_N=0$ for each $N\geq0$.

The reason is that $p=2$ is special is that it is a (the only) ''bad prime'' of the family, i.e. dividing either a numerator or denominator of the local exponents or the polynomial discriminant of $\Delta(t)$. 

A similar behavior can be observed for all the hypergeometric and Beauville cases. 

\subsection{The platonic cases}

Somewhat to our surprise, the ansatz \eqref{eq:U0naiv} for $U(0)$ does not work for the critically scaled platonic operators (Table \ref{tab:platonicoperators}). Taking the diagonal shape does not lead to the $p$-adic expansion as in \eqref{eq:umatrixisrational}. As a consequence it is not possible to evaluate the matrix $U([t_0])\bmod p^N$ for $N\geq2$.

It turns out that one needs to return to the more general ansatz for $U(0)$:
\begin{equation}\label{eq:u0general}
	U(0)=
	\begin{pmatrix}
		1 & 0 \\
		x_p & p 
	\end{pmatrix}.
\end{equation}

We were able to determine the unknown terms $x_p$ in the $U$-matrices of our platonic operators: they are, like for the unscaled Legendre operator, all log-terms of the form
\begin{equation}
	x_p=(p-1)\frac{1}{n}\log_p(c),
\end{equation}
where $n$ is the width of the fibre $I_n$ at $0$.
The $c$ that we detected are summarised in the following table:

\[
\renewcommand{\arraystretch}{1.5}
\begin{array}{c|cc|ccc|ccc}
    \hline
    \text{Operator (Tab. \ref{tab:platonicoperators})}
    & \T(3) & \T(2)
    & \O(4) & \O(3) & \O(2)
    & \I(5) & \I(3) & \I(2)
    \\[1mm]
    \hline
    c
    & 2^3 & 3
    & 3^2 & 2^2 & 3^4
    & 2^3 & 2^5\!\cdot 5^2 & 5^3
    \\[1mm]
    \hline
\end{array}
\]

With those entries for $x_p$ in the $U(0)$, the corresponding matrices $U(t)$ do admit the $p$-adic expansion \eqref{eq:umatrixisrational} and the computation of $U([t_0])\bmod p^N$ is possible.

\subsection{Tate and Katz}
It is surprising that in spite of critical scaling the log-terms occur for the platonic operators. The origin lies in comparison 
between two different natural $q$-coordinates that can be attached to the situation.

To explain it, recall the $q$ expansion ot the normalised Eisenstein series
\[E_4 = 1+240 \sum_{n=1}^{\infty} \frac{n^3 q^n}{1-q^n}, \;\;\;\;E_6= 1-504 \sum_{n=1}^{\infty}\frac{n^5q^n}{1-q^n},\]
from which one can compute the monstrous $q$ expansion of the $j$-function:
\[j=1728\frac{E_4^3}{E_4^3-E_6^2}=\frac{1}{q}+744+196884q+21493760q^2+\ldots\]
In fact these series can be considered as 
$g_2, g_3$ and $j$ of a special elliptic curve
over $\Z[[q]]$, called the {\em Tate curve}
$Tate(q)$, that serves as a model of the universal
deformation of a (split) one-nodal elliptic curve, which arises in the limiting case where $q=0$.
In such a degeneration the j-invariant goes to
$\infty$, visible in the polar part $1/q$ of $j$.
The function $1/j$ goes to $0$ and can be expanded
in a power series:
\[\frac{1}{j(q)}= q - 744 q^2  + 356652 q^3  - 140361152 q^4  + 49336682190 q^5+\ldots.\]
Hence $1/j(q)$ is a local parameter for the Tate curve.
The $q$-parameter then can be expanded in the variable $1/j$:
\[q(1/j) :=\frac{1}{j} + 744 \frac{1}{j^2}  + 750420 \frac{1}{j^3}  + 872769632\frac{1}{j^4}  + 1102652742882\frac{1}{j^5}+\dots \]
For an  elliptic curve given by a Weierstrass equation, we can readily compute $j$ and via the above series one obtains the Tate canonical coordinate $q$. For an elliptic surface $\mathcal{E}\lra \P^1$ with an $I_n$-fibre at $0$,
given in Weierstrass form $y^2=4x^3-g_2(t)x-g_3(t)$,
one obtains a power-series
\[q_{Tate}(t):=q\left(\frac{\Delta(t)}{(12g_2(t))^3}\right) \in \Q[[t]].\]
For instance, for the icosahedral family  $2 3 \underline{5} II$ we have 
\[g_2(t)=-180t^3+\frac{130}{3}t^3+-\frac{10}{3}t+\frac{1}{12}\]
and
\[g_3(t)=432t^5-325t^4+\frac{1855}{27}t^3-\frac{115}{18}t^2+\frac{5}{18}t-\frac{1}{216},\]
so we get
\[\frac{1}{j(t)}= \frac{-t^5  (16 t - 1)^2  (27 t - 2)^3}{(12g_2(t))^3}= 8t^5+380t^6+11510t^7+\ldots \]
and
\[q_{Tate}(t)=q\left(\frac{1}{j(t)}\right)=8t^5+380 t^6+11510t^7+284045t^8+\ldots \in \Z[[t]].\]
It is obvious that the series starts with $t^5$, 
as in the example we have an $I_5$ at the origin.
However, the coefficient $8=2^3$ was not clear
a priori.

On the other hand, for a Picard-Fuchs operator with
a normalised holomorphic and logarithmic solution 
$\phi(t)$ and $\psi(t)=\phi(t)\log(t)+f(t)$ one can
form the $q$-coordinate of the Picard-Fuchs operator
by inversion of the power series
\[q_{PF}(t):=\exp(\psi(t)/\phi(t))=t \cdot \exp(f(t)/\phi(t))=t+\ldots \in \Q[[t]].\]
For our example we find
\[q_{PF}(t)=t(1 + 19/2 t + 429/4 t^2 + 10487/8 t^3+268571/16 t^4 + 35476533/160 t^5 + \ldots).\]
Note that the coefficients are not integral, but can be made integral by the transformation $t\mapsto10t$. The fifth power of $q_{PF}$, multiplied with $8$ is exactly the series $q_{Tate}(t)$. This is a general
fact that is not difficult to prove:
With an  $I_n$ fibre at the origin we have:
\[
\boxed{ q_{Tate}(t)=c\cdot q_{PF}(t)^n},
\]
where $n=\text{ord}_0(1/j(t))$, so that
\[\log(q_{Tate})= \log(c)+n \log(q_{PF}).\]
We call $c$ the {\em TPF-constant}. It can conveniently be computed es
\begin{equation}\label{eq:coeffofj}
	c=\left(\frac{1}{j(t) \exp\left(n \psi(t)/\phi(t)\right)}\right)\Bigg\rvert_{t=0}.
\end{equation}

It is this relation between the two canonical $q$-coordinates that causes the specific form of the
matrix $F(0)$. In Appendix 1 and 2 of the paper of Katz \cite{Katz2} one finds a precise
determination of the Gauss-Manin connection $\nabla$
and the Frobenius $F$ for the Tate elliptic curve, in the $q$-coordinate, in the basis of {\bf H} given by that canonical $1$-form $\omega_{can}$ and its dual form $\eta_{can}$ ($\langle \omega_{can},\eta_{can} \rangle=1$). To describe it let
\[P:=E_2=1-24\sum_{n=1}^{\infty}\frac{n t^n}{1-q^n},\]
\[Q:=E_4=P^2-12\theta(P),\;\;\;\theta=q\frac{d}{dq}.\]
Then one has (see A.1.4.6, p.103 of \cite{Katz2}):

\[\nabla \left(\begin{array}{c}\omega_{can}\\ \eta_{can} \end{array} \right) = \left(\begin{array}{cc}
-P/12&0\\ Q/144&P/12 \end{array} 
\right) \left( \begin{array}{c} \omega_{can}\\\eta_{can} 
\end{array}\right),\]
where $\nabla=\nabla_{\theta}$, and
\[\left(\begin{array}{c}F (\omega_{can})\\ F(\eta_{can}) \end{array} \right) = \left(\begin{array}{cc}
p&0\\ \frac{p\sigma(P)-P}{12}&1 \end{array} 
\right) \left( \begin{array}{c} \omega_{can}\\\eta_{can} 
\end{array}\right).\]
However, in the derivative basis, the matrix of $F$
is {\em diagonal}:
\[\left(\begin{array}{c}F (\omega_{can})\\ F(\nabla(\omega_{can})) \end{array} \right) = \left(\begin{array}{cc}
p&0\\0&1 \end{array} 
\right) \left( \begin{array}{c} \omega_{can}\\\nabla(\omega_{can}) 
\end{array}\right),\]
(see \cite{Katz2}, A2.2.7, p. $109$).
Converting this to the $U$-matrix and performing the
transformation $q_{Tate}=c q_{PF}^n$ leads to:

{\bf Proposition:} {\em Let $\mathcal{E}\rightarrow\P^1$ be an elliptic surface with multiplicative fibre type $I_n$ at $0$ and TPF-constant $c$. Then the constant
$x_p$ in \eqref{eq:u0general} is given by
\begin{equation}
	\boxed{x_p=(p-1)\frac{1}{n}\log_p(c)}.
\end{equation}
}

\section{Monodromy calculations}\label{ch:monodromycalculations}

We expect the monodromy matrices with respect to
the Frobenius basis to be affected by an analogous
logarithmic term. Numerical calculations of the monodromy matrices of the differential equations can be performed conveniently with the ORE-algebra package in SAGE. We used the routines for analytic continuation written by Mezzaroba (\cite{Mezzaroba}).
After selecting a polygonal path in the complex plane avoiding the singularities of the differential equation, the program can perform a
numerical analytic continuation with a {\em certified
precision}. The result is a matrix with floating point numbers. It tuns out that one can identify
all entries and come up with putative exact
matrices. We performed these calculations for all
members of the platonic family.

\begin{figure}[ht]
	\centering

	\begin{subfigure}[t]{0.49\textwidth}
		\centering
		\includegraphics[width=\textwidth]{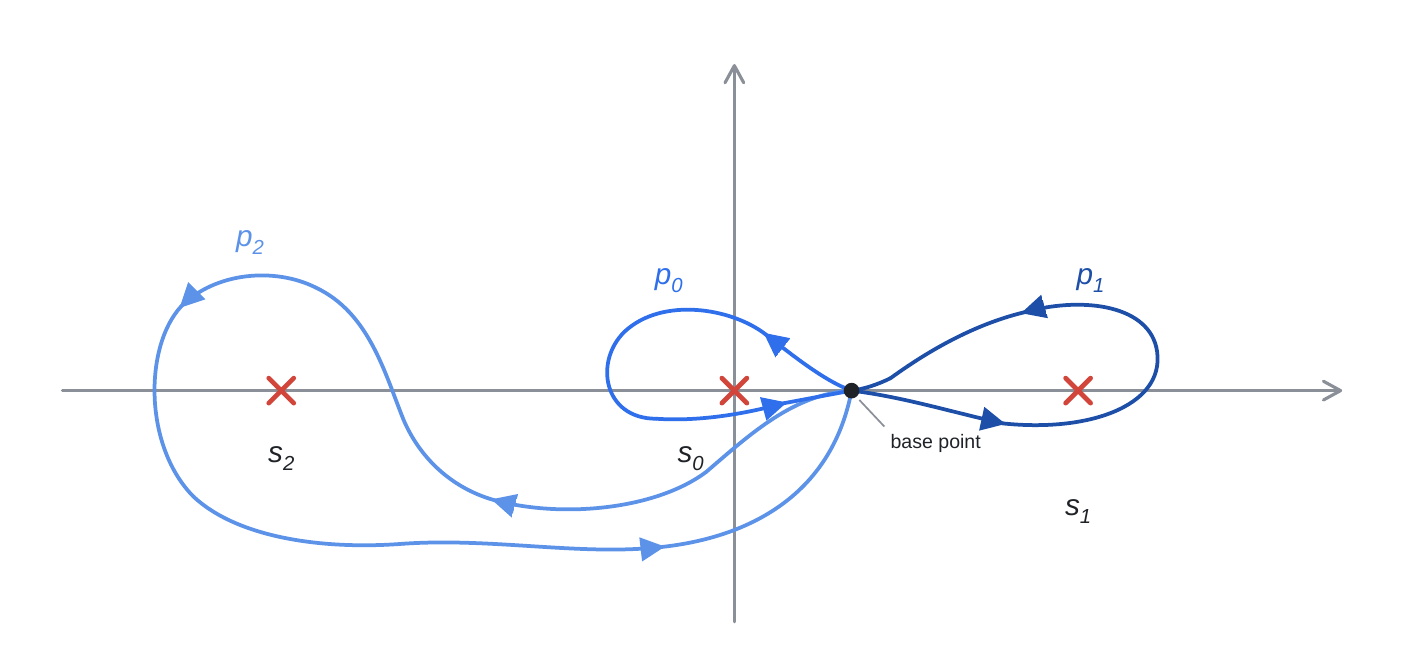}
		\label{fig:config1}
	\end{subfigure}
	\hfill
	\begin{subfigure}[t]{0.49\textwidth}
		\centering
		\includegraphics[width=\textwidth]{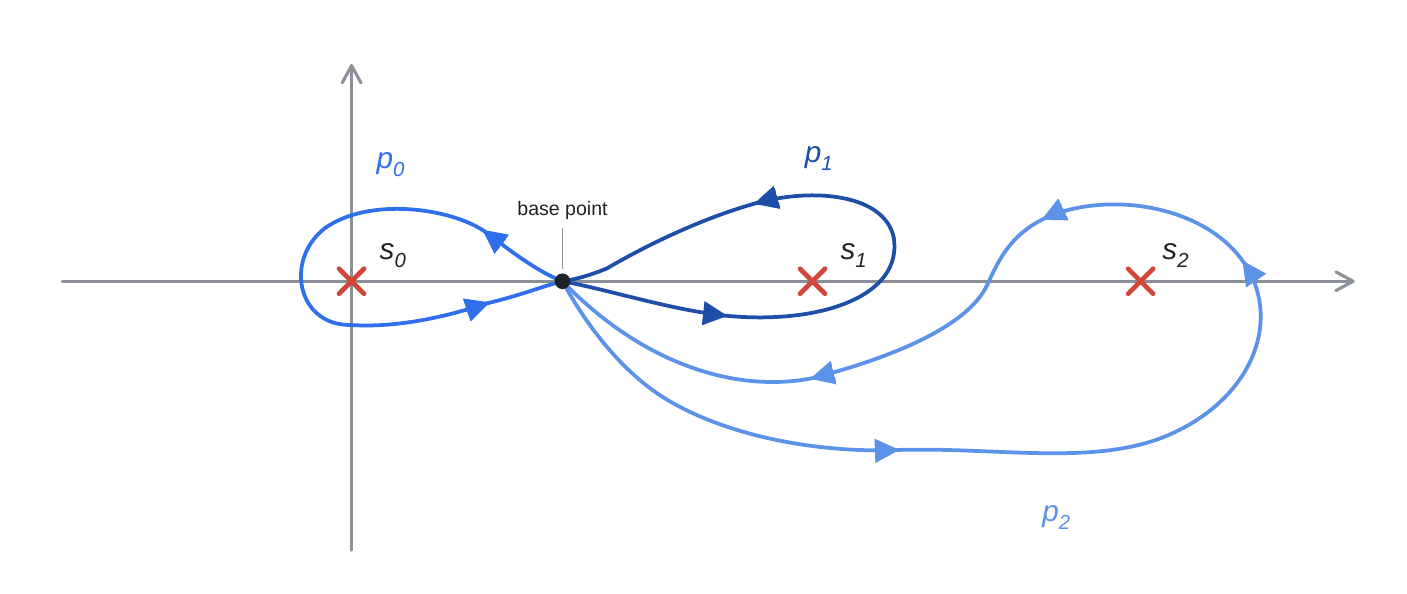}
		\label{fig:config2}
	\end{subfigure}
	\hfill
	\begin{subfigure}[t]{0.49\textwidth}
		\centering
		\includegraphics[width=\textwidth]{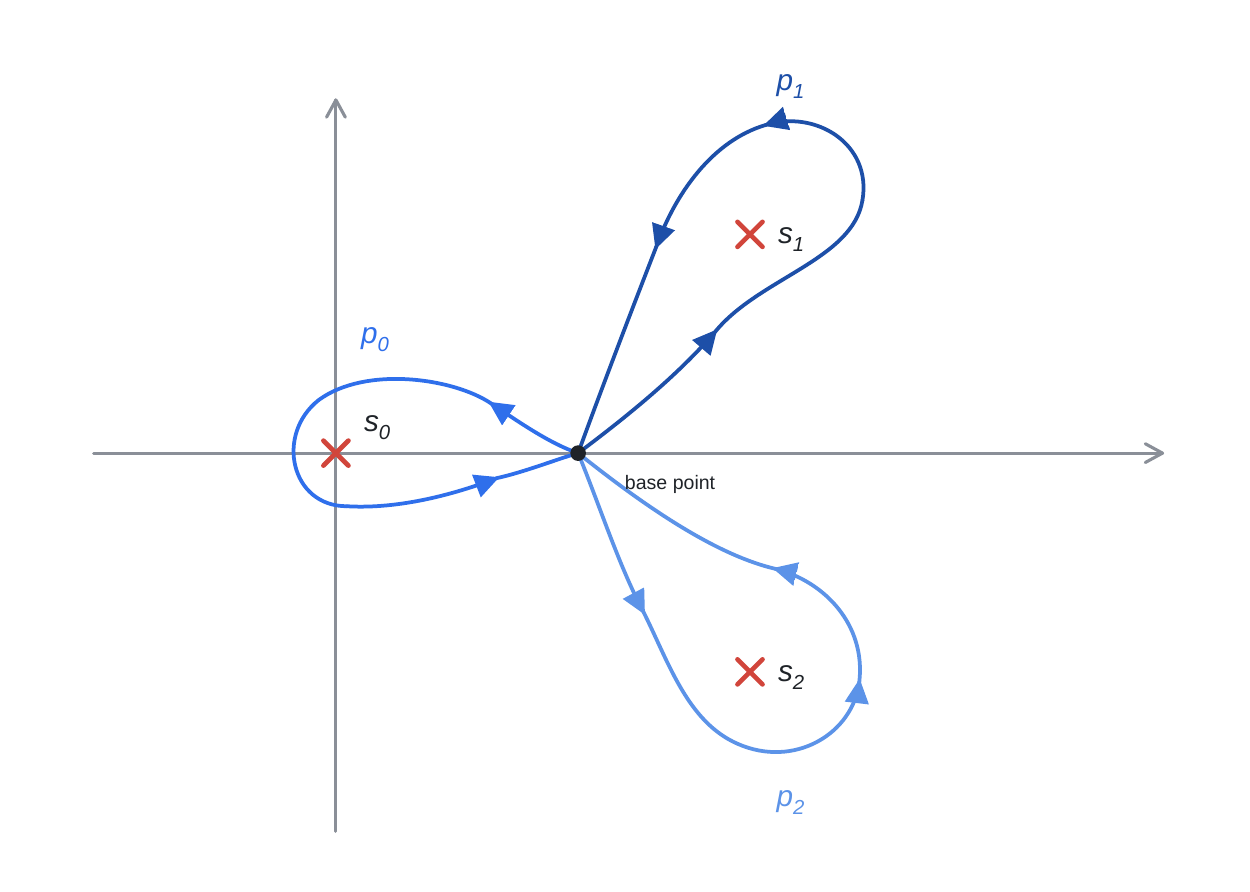}
		\label{fig:config3}
	\end{subfigure}
	
	\caption{The three different configurations of the singular points with
		counterclockwise generators of
		$\pi_1\bigl(\mathbb{P}^1\setminus\{s_0,s_1,s_2,\infty\}\bigr)$.}
	\label{fig:config}
\end{figure}

There are three different configurations of the singular fibres appearing in our examples (Figure \ref{fig:config}). 
With $s_0$ we denote the singular point $0$, with $s_1$ the singular point that is nearest to $0$ in the real case, and the complex point in the upper half plane in the complex conjugate case, and with $s_2$ the remaining singularity $\neq\infty$. The  base point choosen was $0+$, i.e. at a point near $0$ on the positive real axis. The paths we used in each of the three configuration are depicted in Figure \ref{fig:config}.

\begin{figure}[ht]
\includegraphics[width=13cm]{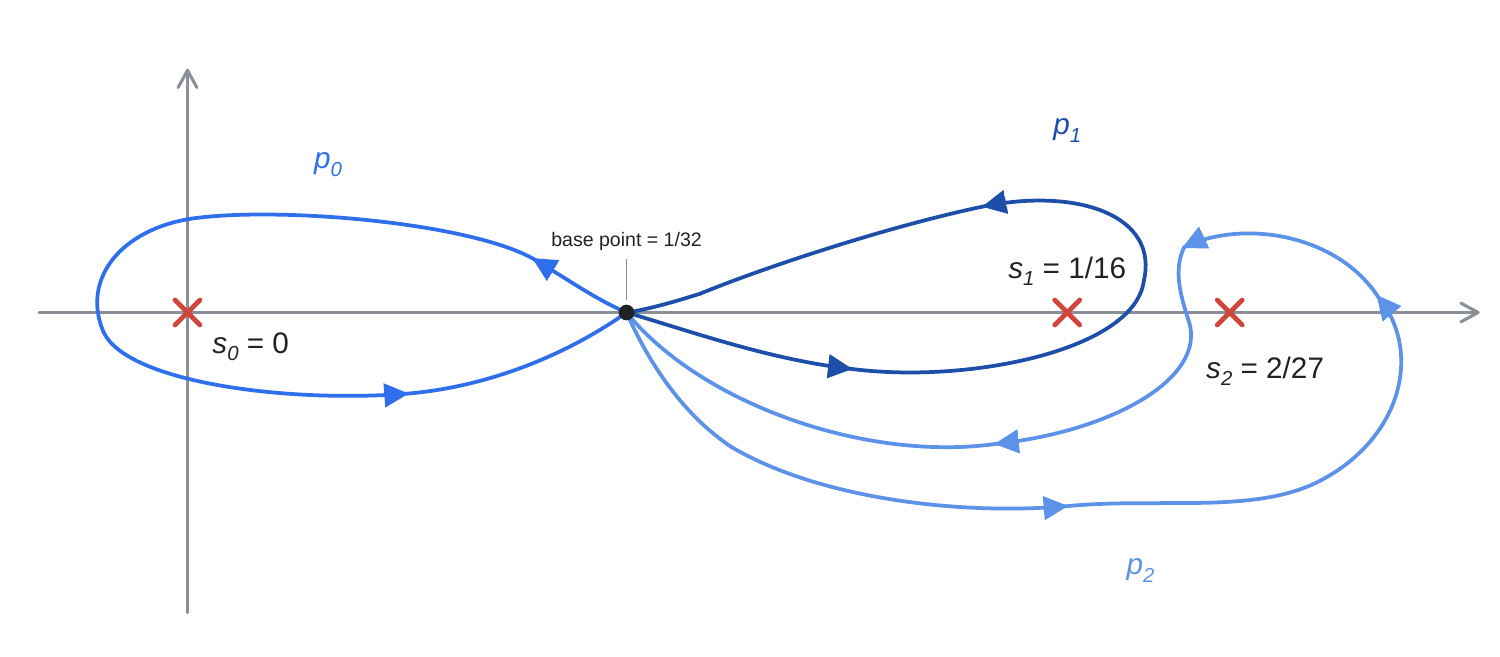}
\caption{Paths around the singular points of $\I(5)$.}
\label{fig:configicosa}
\end{figure}

As an example we give monodromy matrices for the icosahedral operator $\I(5)$ in the Frobenius basis. The singular fibres are $I_5$ at $0$, $I_2$ at $\frac{1}{16}$, $I_3$ at $\frac{2}{27}$ and $II$ at $\infty$ and the paths around the finite singular points we used
are polygonal versions of those depicted in Figure \ref{fig:configicosa}. The result after
identifying the entries of the matrices is:

\begin{equation}
\begin{aligned}
\widetilde{M}_0 =& \begin{pmatrix} 1 & 0 \\ 1 & 1 \end{pmatrix}, \;
\widetilde{M}_1=
\begin{pmatrix}
1+\dfrac{2\log(2^3)}{2\pi i}
&
-10
\\[2mm]
\dfrac{2}{5}
\left(\dfrac{\log(2^3)}{2\pi i}\right)^2
&
1-\dfrac{2\log(2^3)}{2\pi i}
\end{pmatrix}, \\[4mm]
&\widetilde{M}_2=
\begin{pmatrix}
4+\dfrac{3\log(2^3)}{2\pi i}
&
-15
\\[2mm]
\dfrac{3}{5}
\left(
1+\dfrac{\log(2^3)}{2\pi i}
\right)^2
&
-2-\dfrac{3\log(2^3)}{2\pi i}
\end{pmatrix}.
\end{aligned}
\end{equation}

There is a matrix that transforms the above matrices simultaneously into integral matrices:
\[ T=
\begin{pmatrix}
	1 & 0 \\
	-\frac{1}{2\pi i} \cdot \color{Green} \log(2^3)\normalcolor & 5
\end{pmatrix}.
\]

The integer transformed matrices $T\cdot \widetilde{M}\cdot T^{-1}$ are:
\[
M_0=\begin{pmatrix} 1 & 0 \\ 5 & 1 \end{pmatrix}, \quad
M_{1/16}=\begin{pmatrix} 1 & -2 \\ 0 & 1 \end{pmatrix}, \quad
M_{2/27}=\begin{pmatrix} 4 & -3 \\ 3 & -2 \end{pmatrix}.
\]

The computation was performed for all operators from the platonic family (Table \ref{tab:herfurtner-operators}) and the results are summarised in Appendix \ref{ap:herfurtnersurfaces}, Table \ref{tab:platonicoperators}. We observe that in all cases the non-rational entry appearing in the transformation matrix $T$ is equal to the TPF constant $c$ \eqref{eq:coeffofj} of the corresponding family.

\subsection{The platonic monodromy groups} 
Given some matrices in $\mathrm{SL}_2(\Z)$ like the monodromy matrices $M_0$, $M_{s_1}$, $M_{s_2}$, one may ask what sort of subgroup $\Gamma$ of 
$\mathrm{SL}_2(\Z)$ they generate. There are algorithms
like {\em Todd-Coxeter coset enumeration} \cite{ToddCoxeter} that can be used to determine the index of $\Gamma$ and in fact a permutation representation of $\mathrm{SL}_2(\Z)$ may be constructed that defines the
subgroup $\Gamma$. If $-I \in \Gamma$, then $\Gamma$ is the
pull-back of a subgroup $\overline{\Gamma}$ of the modular group $\mathrm{PSL}_2(\Z)$, with the same index.
The study of such subgroups and the enumeration
of those of low index are a classical topic of
research. The {\em signature} of such a subgroup is the tuple $(n,g;\nu_2,\nu_3,\nu_{\infty})$, where $n$ denotes the index, $g$ the genus (equal to the genus of $\H/\overline{\Gamma}$),  $\nu_2$, $\nu_3$ denoting the number of elliptic points of order $2$, $3$ and $\nu_{\infty}$ the
number of cusps. These invariants are connected by the {\em Riemann-Hurwitz relation}
\[g=1+\frac{n}{12}-\frac{\nu_2}{4}-\frac{\nu_3}{3}-\frac{\nu_{\infty}}{2}.\]
A further invariant is the list of {\em cusp widths} $(m_1,\dots,m_{\nu_\infty})$ of each of the $\nu_{\infty}$ cusps.

As we have explained before, the index $n$ divides the degree $\deg(J)$ of the $J$-map (see Section \ref{ch:ellipticsurfaces}). Recall that for a modular elliptic surface the index is equal to the degree of the $J$-map.

After performing the computation for all $11$ elliptic surfaces of the platonic family, we make the following interesting observations: 

The monodromy groups of the $4$ elliptic surfaces with fibre type $II$ are of index $10$ and genus $0$. They have cusp widths $(1,1,8)$, $(1,2,7)$, $(1,4,5)$ and $(2,3,5)$. These make
up all subgroups of $\mathrm{PSL}_2(\Z)$ of index $10$ of genus $0$ with three cusps.

Something similar holds for the $4$ elliptic surfaces with fibre type $III$: they have index $9$, genus $0$ and cusp widths $(1,1,7)$, $(1,2,6)$, $(1,3,5)$ and $(2,3,4)$. These make up all subgroups of $\mathrm{PSL}_2(\Z)$ of index $9$ of genus $0$ with three cusps.

Hence these $8$ surfaces with fibre type $II$ and $III$ from the platonic family  are  modular: in each case $\deg J=[\mathrm{PSL}:\bar\Gamma]$, so $\widetilde{\tau}$ from \eqref{eq:factorj} is an
isomorphism and the base $\mathbb{P}^1$ is the modular curve $X(\bar\Gamma)$ of the monodromy group. In fact these surfaces are modular.

This is a similar phenomenon as for the $6$ Beauville surfaces where the monodromy groups are in one-to-one correspondence to the conjugacy classes of subgroups of $\mathrm{PSL}_2(\Z)$ of index $12$, genus $0$ and four cusps. The difference is that in this case
the groups are congruence subgroups, whereas
in the $8$ platonic cases we are dealing with {\em non-congruence subgroups}. 

Every larger family has its strange members that have habits that run against traditions. In this case the outliers are the surfaces with fibre type $IV$. For the tetrahedral surface
$223 IV$ and its isogenous partner $116IV$ one
has monodromy groups   $\Gamma_0(3)\subset\mathrm{SL}_2(\Z)$, which is a congruence subgroup of index $8$, but with $-I\notin\Gamma_0(3)$; hence its image in $\mathrm{PSL}_2(\Z)$ has index $4<8=\deg(J)$. 
The black sheep in the family is the surface $125IV$. Its monodromy group is the full group $\mathrm{PSL}_2(\Z)$: the matrices for $12\underline{5}IV$ are given in terms of the standard generators $S,T$ of $\mathrm{SL}_2(\Z)$:
\begin{align*}
 M_{s_1} &= \smat{ 1 & -1 \\ 0 & 1 } = T^{-1}, \\
 M_{s_2} &= \smat{ 5 & -2 \\ 8 & -3 } = S\,T^{-2}S^{-1}\,T^{-2}\,S\,T^{2}S^{-1}, \\
 M_{0} &= \smat{ 1 & 0 \\ 5 & 1 } = S\,T^{-5}S^{-1}.
\end{align*}
Furthermore we have $T=M_{s_1}^{-1}$ and
$S=M_{s_1}^{2}M_0^{-1}M_{s_2}^{-1}M_{s_1}^{-1}M_0^{-1}M_{s_2}$.
Hence the three matrices generate $\mathrm{SL}_2(\Z)$ and their images generate $\mathrm{PSL}_2(\Z)$.

This is consistent with the fact that there are no subgroups of $\mathrm{PSL}_2(\Z)$ of index $8$ of genus $0$ and three cusps: such a group would satisfy $3\nu_2+4\nu_3=2$, which is impossible.

\section{Outlook}\label{ch:outlook}
Starting from the reflection groups of the platonic solids in three-space we were led to an interesting class of elliptic surfaces. There are several ways one can look for variations and generalisations. We briefly mention a few further directions of inquiry.

\subsection{Other sections}
In Section \ref{ch:platonicellipticsurfaces} we restricted the elliptic threefold to the lines $P=r$, where $r$ is a non-zero constant. A different choice of the lines in the $P-Q$ plane lead to different elliptic surfaces. Table \ref{tab:platonic-specializations} shows all possible configurations of lines. Here $r$ denotes a constant and $t$ the parameter of the family.

\setlength{\extrarowheight}{1pt}

\setcellgapes{3pt}
\makegapedcells

\begin{table}[ht]
    \centering
    \small
    \renewcommand{\arraystretch}{1.85}
    \setlength{\extrarowheight}{2pt}
    \setlength{\tabcolsep}{4pt}

    \begin{tabularx}{\textwidth}{
        >{\raggedright\arraybackslash}p{0.22\textwidth}
        |>{\centering\arraybackslash}X
        >{\centering\arraybackslash}X
        >{\centering\arraybackslash}X
    }
        \hhline{|----|}
        \textbf{Specialisation}
        &
        \textbf{Tetrahedron}
        &
        \textbf{Octahedron}
        &
        \textbf{Icosahedron}
        \\
        \hhline{|====|}

        $P=r,\ Q=t$, $r\neq 0$
        &
        $\mathrm{I}_2\enspace
         \mathrm{I}_3\enspace
         \mathrm{I}_3\enspace
         \mathrm{IV}$
        &
        $\mathrm{I}_4\enspace
         \mathrm{I}_2\enspace
         \mathrm{I}_3\enspace
         \mathrm{III}$
        &
        $\mathrm{I}_2\enspace
         \mathrm{I}_5\enspace
         \mathrm{I}_3\enspace
         \mathrm{II}$
        \\
        \hhline{|----|}

        $P=0,\ Q=t$
        &
        $\mathrm{IV}^{*}\enspace \mathrm{IV}$
        &
        $\mathrm{III}^{*}\enspace \mathrm{III}$
        &
        $\mathrm{II}^{*}\enspace \mathrm{II}$
        \\
        \hhline{|----|}

        $P=t,\ Q=r$, $r\neq 0$
        &
        $\mathrm{I}_3\enspace
         \mathrm{I}_3\enspace
         \mathrm{I}_3\enspace
         \mathrm{I}_3$
        &
        \makecell[tc]{
            $\mathrm{I}_2\enspace
             \mathrm{I}_2\enspace
             \mathrm{I}_3\enspace
             \mathrm{I}_3\enspace\mathrm{I}_8^{*}$
        }
        &
        \makecell[tc]{
            $\mathrm{I}_5\enspace
             \mathrm{I}_5\enspace
             \mathrm{I}_5\enspace
             \mathrm{I}_3$\\
            $\mathrm{I}_3\enspace
             \mathrm{I}_3\enspace
             \mathrm{I}_6^{*}$
        }
        \\
        \hhline{|----|}

        $P=t,\ Q=0$
        &
        \multicolumn{3}{c}{Singular generic fibre}
        \\
        \hhline{|----|}

        \makecell[tl]{
            $P=mt+r,\ Q=t$,\\
            $m,r\neq0$, generic
        }
        &
        \makecell[tc]{
            $m^2r\neq4$\\[1mm]
            $\mathrm{I}_2\enspace
             \mathrm{I}_3\enspace
             \mathrm{I}_3\enspace
             \mathrm{I}_3\enspace
             \mathrm{I}_1$
        }
        &
        \makecell[tc]{
            $mr\notin\left\{1,\frac34\right\}$\\[1mm]
            $\mathrm{I}_4\enspace
             \mathrm{I}_2\enspace
             \mathrm{I}_2\enspace
             \mathrm{I}_3\enspace\mathrm{I}_3\enspace
             \mathrm{I}_4^{*}$
        }
        &
        \makecell[tc]{
            $mr^2\notin
             \left\{\frac4{27},-\frac45\right\}$\\[1mm]
            $\mathrm{I}_2\enspace
             \mathrm{I}_5\enspace
             \mathrm{I}_5\enspace
             \mathrm{I}_5$\\
            $\mathrm{I}_3\enspace
             \mathrm{I}_3\enspace
             \mathrm{I}_3\enspace
             \mathrm{I}_4^{*}$
        }
        \\
        \hhline{|----|}

        \makecell[tl]{
            $P=mt+r,\ Q=t$,\\
            tangent
        }
        &
        \makecell[tc]{
            $m^2r=4$\\[1mm]
            $\mathrm{I}_2\enspace
             \mathrm{I}_6\enspace
             \mathrm{I}_3\enspace
             \mathrm{I}_1$
        }
        &
        \makecell[tc]{
            $mr=1$\\
            $\mathrm{I}_4\enspace
             \mathrm{I}_4\enspace
             \mathrm{I}_3\enspace
             \mathrm{I}_3\enspace
             \mathrm{I}_4^{*}$\\[2mm]
            $mr=\frac34$\\
            $\mathrm{I}_4\enspace
             \mathrm{I}_2\enspace
             \mathrm{I}_2\enspace
             \mathrm{I}_6\enspace
             \mathrm{I}_4^{*}$
        }
        &
        \makecell[tc]{
            $mr^2=\frac4{27}$\\
            $\mathrm{I}_2\enspace
             \mathrm{I}_{10}\enspace
             \mathrm{I}_5\enspace
             \mathrm{I}_3$\\
            $\mathrm{I}_3\enspace
             \mathrm{I}_3\enspace
             \mathrm{I}_4^{*}$\\[2mm]
            $mr^2=-\frac45$\\
            $\mathrm{I}_2\enspace
             \mathrm{I}_5\enspace
             \mathrm{I}_5\enspace
             \mathrm{I}_5$\\
            $\mathrm{I}_6\enspace
             \mathrm{I}_3\enspace
             \mathrm{I}_4^{*}$
        }
        \\
        \hhline{|----|}

        $P=mt,\ Q=t$, $m\neq0$
        &
        $\mathrm{IV}^{*}\enspace
         \mathrm{I}_3\enspace
         \mathrm{I}_1$
        &
        \makecell[tc]{
            $\mathrm{III}^{*}\enspace
             \mathrm{I}_2\enspace
             \mathrm{I}_3\enspace\mathrm{I}_4^{*}$
        }
        &
        \makecell[tc]{
            $\mathrm{II}^{*}\enspace
             \mathrm{I}_5\enspace
             \mathrm{I}_5\enspace
             \mathrm{I}_3\enspace\mathrm{I}_3\enspace
             \mathrm{I}_4^{*}$
        }
        \\
        \hhline{|----|}
    \end{tabularx}

    \caption{Kodaira fibre configurations of the tetrahedral,
    octahedral, and icosahedral families for general lines in the $P-Q$ plane.}
    \label{tab:platonic-specializations}
\end{table}

\subsection{Other reflection groups}
One may wonder if there are other rings of invariants that lead to elliptic
fibrations. It turns out there are quite a number of further cases that deserve
closer study. One can look at invariant rings of rotation subgroups of  complex reflection groups. We mention the group $G_{168}$ (nr. 24 in \cite{ShephardTodd}),  with (semi)-invariants of degree $4$, $6$, $14$ and $21$ 
which leads to an elliptic surface of K3-type with fibres $I_4\, I_4\,I_3\, I_3\, II^*$. For the Valentiner group (nr. 27 in \cite{ShephardTodd}), there are
(semi)-invariants of degrees $6$, $12$, $30$ and $45$, producing a K3 elliptic surface
with singular fibres $I_4\, I_5\, I_3\, I_3\, III^*$.

\subsection{Other invariant rings} But also invariant rings of more general algebraic groups may be considered. For example, for the action of $\mathrm{SL}_2(\C)$ acting on binary forms of degree
five, it is a classical fact that its ring of invariants of binary quintics is generated by invariants $\mathrm{I}_4, \mathrm{I}_8, \mathrm{I}_{12}, \mathrm{I}_{18}$, with a single relation, \cite{duPlessisWall}.
In the paper \cite{Abdesselam} simpler invariants $J,K,L,H$ are defined, among which there
exists the relation (eq. (32) in \cite{Abdesselam}):
\[ 16H^2 = -432L^3 - 72L^2 K J + 8L K^3 - 2LK^2 J^2 + L^2 J^3 + K^4J.\] 
For fixed value of $J \neq 0$ we obtain over the $K$-line the cousin $135III$ of the octahedral case $\O$  as family of elliptic curves written in the variables $L$ and $H$.

\subsection{Other Herfurtner operators}
In the list of Herfurtner, additionally to the elliptic surfaces that we summarised as the platonic family, there are some further cases where the Picard-Fuchs operator has degree $2$. For example, there are
$13$ cases with two semi-stable and two reduced non$-I_n$-fibres. 
By bringing these non semi-simple fibres to $0$ and $\infty$ one obtains an operator of degree $2$; of course
these operators have finite local monodromy at $0$ and $\infty$
but lead to two-term recursion for the coefficients.

\subsection{Other dimensions}
In this paper we were mainly concerned with
elliptic surfaces and associated second-order operators, which belong to a well studied field. Families of
K3 surfaces with Picard number $19$ lead to interesting third-order operators and families
of Calabi-Yau threefolds with $h^{12}=1$ lead 
to fourth-order operators with very nice arithmetical
properties. Nice examples of Calabi-Yau operators
can be constructed as Hadamard products from the hypergeometric cases (\ref{tab:euclidean-elliptic-surfaces}) and Beauville cases (\ref{tab:beauville-elliptic-surfaces}).
Bradley Klee asked if the platonic operators
may be used to build interesting higher order
operators. So far we have not been able to capitalise
on his idea and use it to construct new Calabi-Yau operators.
On the other hand it may be rewarding to take a look
at higher rank reflection groups which may give
rise to interesting pencils of K3 surfaces or Calabi-Yau threefolds.

\subsection{Ap\'ery-like limits}
In \cite{Apery} Ap\'ery used the solution $(a_n)$ to
his recursion (with $a_{-1}=0$, $a_0=1$), considered a second solution $(b_n)_n$ of the same recursion (with $b_0=0$, $b_1=1$) and showed
that the limit
\[
	L:=\lim_{n\to\infty}\frac{b_n}{a_n}
\]
is equal to $\zeta(2)/5$ and proved thereby its irrationality.
But also for much more general Picard-Fuchs operators
and corresponding recursions one can form the
sequences $(a_n)$ and $(b_n)$ and try to identify the limit $L$. This has been
done for the Beauville operators by various authors.
To our knowledge it was not done yet for the platonic family.
Table \ref{tab:herfurtner-operators} in Appendix \ref{ap:herfurtnersurfaces} collects for each platonic operator the limit. In some cases, the number was identified and turns out to be a logarithmic term. In other cases, for example for $\O(2)$ and $\I(2)$, the limit was computed to high precision, but we, nor {\tt ChatGPT} or {\tt Claude}, were not able to identify the constant up to now. These cases are marked by ${\bf ?}$ in Table \ref{tab:herfurtner-operators}.

\appendix
\section{Goursat invariants}\label{ap:goursatinvariants}

The invariant $P$ always has the form
\[
P=x^2+y^2+z^2.
\]

\emph{Tetrahedron:}
\begin{align*}
	Q&= x^2y^2+y^2z^2+x^2z^2,\\
	R&= xyz,\\
	S&= (x^2-y^2)(y^2-z^2)(z^2-x^2).
\end{align*}

\emph{Octahedron:}
\begin{align*}
	Q&= x^2y^2+y^2z^2+x^2z^2, \\
	R&= x^2y^2z^2, \\
	S&= xyz(x^2-y^2)(y^2-z^2)(z^2-x^2).
\end{align*}

\emph{Icosahedron:}
\begin{align*}
	Q&= z (-5 z^3 (x^2 + y^2) + 5 z (x^2 + y^2)^2 - 2 x (x^4 - 10 x^2 y^2 + 5 y^4) + z^5), \\
	R&= (4 x^2+6 x z+z^2) (x^4+8 x^3 z-10 x^2 y^2+14 x^2 z^2-8 x
	z^3+5 y^4-10 y^2 z^2+z^4) \\ & \quad \quad \quad (x^4-2 x^3 z-10 x^2 y^2-x^2 z^2+30 x
	y^2 z+2 x z^3+5 y^4-25 y^2 z^2+z^4), \\
	S&= 8 y (5 x^4-10 x^2 y^2+y^4) (x^2-x z-z^2) \\
	& \quad \quad \quad(x^4-10 x^2 y^2+5 y^4-12 x^3 z+20
	x y^2 z+44 x^2 z^2-20 y^2 z^2-48 x z^3+16 z^4) \\
	& \quad \quad \quad(x^4-10 x^2 y^2+5 y^4+8
	x^3 z-40 x y^2 z+24 x^2 z^2-40 y^2 z^2+32 x z^3+16 z^4).
\end{align*}

\section{The cousins: surfaces with three
\texorpdfstring{$I_n$}{I-n}-fibres and one
\texorpdfstring{$II$}{II},
\texorpdfstring{$III$}{III} or
\texorpdfstring{$IV$}{IV}-fibre}\label{ap:herfurtnersurfaces}

\subsection{Picard-Fuchs operators}
By multiplying the solutions with an elementary factor, the solutions of any differential equation with four regular singular singularities can be reduced to the {\em Heun equation}
\[ \frac{d^2y}{dx^2} +\left(\frac{\gamma}{x}+\frac{\delta}{x-1}+\frac{\epsilon}{x-t}\right)\frac{dy}{dx}+\frac{\alpha\beta x-q}{x(x-1)(x-t)}y=0\]
with Riemann symbol
\[\left\{
\begin{array}{cccc}
0&1&t&\infty\\
\hline
0&0&0&\alpha\\
1-\gamma&1-\delta&1-\epsilon&\beta\\ \end{array}
\right\}\]
and accessory parameter $q$.
The Picard-Fuchs operators which arise from elliptic surfaces from Herfurtner's list are all 
of Heun type and can be found in \cite{MovesatiReiter}.
Here we collect the basic data for the platonic family (Table \ref{tab:platonic-family}), i.e. those elliptic surfaces with four {\em reduced} fibres, three of which are semi-stable.
The Picard-Fuchs operators that annihilate the normalised period integrals are of degree $2$,
hence lead to two-term recursions. The operators can be written as
\begin{equation}\label{genformeq}
	a\theta^2-t(b\theta^2+b\theta+\lambda)+ct^2(\theta+\alpha)(\theta+\beta),
\end{equation}
where $a,b,c$ are integers and the local exponents are
\begin{align*}
	\left(\alpha,\beta\right)
	=\begin{cases} \left(\frac{5}{6},\frac{7}{6}\right), \; \text{for} \; II, \\
		\left(\frac{3}{4},\frac{5}{4}\right), \; \text{for} \; III, \\
		\left(\frac{2}{3},\frac{4}{3}\right), \; \text{for} \; IV.
	\end{cases}
\end{align*}

The Riemann symbol for such an operator has the form
\begin{equation}
	\begin{Bmatrix}
		0 & * & * & \infty \\
		\hline
		0 & 0 & 0 & \alpha\\
		0 & 0 & 0 & \beta
	\end{Bmatrix},
      \end{equation}

so we are dealing with a {\em Heun equation} with a special exponent structure: $${\gamma=\delta=\epsilon=1}.$$ Around $0$ and around the two singularities denoted by $*$, there is a holomorphic solution $\phi$ and a logarithmic solution $\psi$, so that the local monodromy matrices are conjugate to the Jordan block
\begin{equation}
	\begin{pmatrix}
		1 & 1 \\
		0 & 1
              \end{pmatrix}.
            \end{equation}
At $\infty$, the monodromy is  of finite order $=:e$ and conjugate to the matrix
\begin{equation}         
\begin{pmatrix}
	e^{2\pi i \alpha} & 0 \\
	0 & e^{2\pi i \beta} 
\end{pmatrix}.
\end{equation}
Indeed, pulling back the elliptic surface near $\infty$ by an $e$-fold cover, 
the fibre at $\infty$ is a smooth elliptic curve with an automorphism of order $e$.

Table \ref{tab:herfurtner-operators} collects some information on these operators: we
record the Weierstrass covariants $g_2$ and $g_3$ for the elliptic surfaces,\footnote{In \cite{Herfurtner} there is a typo in the table on p. $337$. For the case $I_1I_3I_5III$ one should replace $100X^3Y$ by $100XY^3$ in $G_3$. } the critically scaled
Picard-Fuchs operator and the location of the singular points,  the first few coefficients of the integer sequence, a Laurent polynomial generating that sequence, the Ap\'ery constant $L$ (if available), the integral monodromy matrices, the monodromy groups $\Gamma_{i;(m_1,m_2,m_3)}\subset \mathrm{PSL}_2(\Z)$ generated by those matrices, where $i$ denotes the index and $m_k$ the cusp widths, the transformation matrix $T$ connecting the Frobenius and integral basis and the TPF-constant $c$.

\newgeometry{top=1.5cm,bottom=1cm,left=1cm,right=1cm}

\begingroup
\setlength{\tabcolsep}{1.5pt}
\setlength{\LTleft}{\fill}
\setlength{\LTright}{\fill}

\setlength{\LTpre}{\fill}
\setlength{\LTpost}{\fill}
\renewcommand{\arraystretch}{2}
\fontsize{9.2pt}{11pt}\selectfont

\newsavebox{\npformula}
\newsavebox{\nppicture}
\newlength{\npfree}
\newcommand{\Fnp}[2]{%
	\sbox{\npformula}{#1}%
	\sbox{\nppicture}{\raisebox{-5pt}{#2}}%
	\setlength{\npfree}{\linewidth}%
	\addtolength{\npfree}{-\wd\nppicture}%
	\addtolength{\npfree}{-0.6em}%
	\hbox to \linewidth{%
		\ifdim\wd\npformula>\npfree
		\resizebox{\npfree}{!}{\usebox{\npformula}}%
		\else
		\usebox{\npformula}%
		\fi
		\hfill\usebox{\nppicture}}%
}

\begin{longtable}{|c|p{7.75cm}|c|c|c|}
	\caption{The platonic family.}
	\label{tab:herfurtner-operators}\\
	\hline
	{\bf Notation} & \begin{tabular}{@{}p{7.45cm}@{}}
		{\bf $g_2(t),g_3(t)$ for the Weierstrass form} \\ {\bf Picard-Fuchs operator}\\ {\bf Critically scaled $(a_n)_n$}\\ {\bf Laurent polynomial}\\
		{\bf Decimal approximation of $L=\lim_{n\to\infty} b_n/a_n$ and exact value (if identified)}\\
		{\bf Monodromy matrices $M_0$, $M_{s_1}$, $M_{s_2}$}
	\end{tabular}  & \begin{tabular}{c}
		{\bf Configuration of} \\ {\bf Singular fibres} \\
		$0$, $s_1$, $s_2$, $\infty$
	\end{tabular} & \begin{tabular}{c} {\bf TPF-constant $c$} \\ \\
		\hline \\ {\bf Transformation matrix} \end{tabular} & \multirow{2}{*}{\makecell[c]{
        {\bf Monodromy}\\
        {\bf group} \\
        $\Gamma_{i;(m_1,m_2,m_3)}$
        }} \\
	\triplehline
	\endfirsthead

	\multicolumn{5}{c}{\tablename\ \thetable{} -- continued}\\
	\hline
	\endhead

	\endlastfoot
	$11\underline{8}II$ & \begin{tabular}{@{}p{7.75cm}@{}}
		$g_2=-64t^3 + 20t^2 - 2t + 1/12$ \\
		$g_3=64t^5 - 80t^4 + 64/3t^3 - 8/3t^2 + 1/6t - 1/216$ \\
		\hline $3\theta^2-t\cdot(56\theta^2 + 56\theta +
		18)+t^2\cdot432\left(\theta+\frac{7}{6}\right)\left(\theta+\frac{5}{6}\right)$ \\
		\hline
		$a_n=1, 6, 30, 12, -2250, -35244, -329364,\dots$ (See \cite{NolanDalbyCondon}) \\
		\hline
		\Fnp{$F=\frac{3x + 4x^2 - xy^2 + 4y + 6xy}{xy}$}{\includegraphics{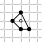}} \\
		\hline
		$L:$ does not exist\\
		\hline
		$M_{s_1}=\smat{ -1 & -1 \\ 4 & 3}$, $M_{s_2}=\smat{ 2 & -1 \\ 1 & 0 }$, $M_0=\smat{ 1 & 0 \\ 8 & 1 }$ \end{tabular} & \begin{tabular}{cccc} $I_1$ & $I_1$ & $I_8$
		& $II$ \\ \multicolumn{2}{c}{{\tiny $\frac{7\pm4i\sqrt{2}}{108}$}} & $0$ & $\infty$ \end{tabular} & \begin{tabular}{c} \\ $-2^8\cdot3$ \\  \\ \hline \\ $T=
		\smat{ 1 & 0 \\ -\frac{1}{2\pi i}\,{\color{Green}\log\!\left(-2^8\cdot 3\right)} & 8}$\\ \end{tabular} & \makecell[c]{$\Gamma_{10;(1,1,8)}$}
	\\
	\hline
	\hline
	$12\underline{7}II$ & \begin{tabular}{@{}p{7.75cm}@{}}
		$g_2=-28t^3 + 14t^2 - 2t + 1/12$ \\
		$g_3=16t^5 - 29t^4 + 37/3t^3 - 13/6t^2 + 1/6t - 1/216$ \\
		\hline
		$6\theta^2-t\cdot(113\theta^2 + 113\theta +
		36)+t^2\cdot432\left(\theta+\frac{7}{6}\right)\left(\theta+\frac{5}{6}\right)$ \\
		\hline
		$a_n=1, 6, 48, 444, 4500, 48456, 543816, \dots$ \\
		\hline
		\Fnp{$F=\frac{3x + 2x^2 + 3x^2y + xy^2 + y + 6xy}{xy}$}{\includegraphics{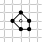}} \\
		\hline
		$L=0.27279178641\dots=\frac{1}{7}\log\!\left(\frac{3^3}{2^2}\right)$\\
		\hline
		$M_{s_1}=\smat{ 1 & -1 \\ 0 & 1 }$, $M_{s_2}=\smat{ 5 & -2 \\ 8 & -3 }$, $M_0=\smat{ 1 & 0 \\ 7 & 1 }$ \end{tabular} & \begin{tabular}{cccc} $I_1$ & $I_2$ &
		$I_7$ & $II$ \\ $\frac{2}{27}$ & $\frac{3}{16}$ & $0$ & $\infty$
	\end{tabular} & \begin{tabular}{c} \\ $2\cdot3^2$ \\  \\ \hline \\ $T=\smat{ 1 & 0 \\ -\frac{1}{2\pi i}\,{\color{Green}\log\!\left(2\cdot 3^{2}\right)} & 7 }$\\ \end{tabular} & \multirow{1}{*}{\makecell[c]{$\Gamma_{10;(1,2,7)}$}}
	\\
	\cline{1-5}
	$1\underline{2}7II$ & \begin{tabular}{@{}p{7.75cm}@{}}
		Substitute $t\mapsto -49t+3/16$ \\
		\hline $3\theta^2-t\cdot(2080\theta^2 + 2080\theta +
		684)$\\$+t^2\cdot338688\left(\theta+\frac{7}{6}\right)\left(\theta+\frac{5}{6}\right)$ \\
		\hline
		$a_n=1, 228, 64596, 20133456, 6666855300,\dots$  \\
		\hline
		\Fnp{$F=\frac{N}{xy^2}$, where}{\includegraphics{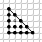}} \\
		$\begin{aligned}
				N={}&-729+36y+186y^2+4y^3-9y^4-324x\\
				&-100xy+228xy^2-12xy^3+450x^2-36x^2y\\
				&+50x^2y^2+108x^3+36x^3y-81x^4
			\end{aligned}$ \\
		\hline
		$L=0.00903720141\dots=\displaystyle
		\frac{1}{16\sqrt{7}}
		\log\!\left(\frac{29+4\sqrt{7}}{3^3}\right)$\\
		\hline
		$M_{s_1}=\smat{ -1 & -4 \\ 1 & 3 }$, $M_0=\smat{ 1 & 0 \\ 2 & 1 }$, $M_{s_2}=\smat{ 1 & -7 \\ 0 & 1 }$ \end{tabular} & \begin{tabular}{cccc} $I_1$ & $I_2$ & $I_7$ & $II$ \\ $\frac{1}{432}$ &
		$0$ & $\frac{3}{784}$ & $\infty$ \end{tabular} & \begin{tabular}{c} $-3^7$ \\ \\ \hline \\ $T=\smat{ 1 & 0 \\ -\frac{1}{2\pi i}\,{\color{Green}\log\!\left(-3^{7}\right)} & 2 }$ \\ \end{tabular} & \multirow{2}{*}{\makecell[c]{$\Gamma_{10;(1,2,7)}$}}
	\\
	\cline{1-4}
	$\underline{1}27II$ & \begin{tabular}{@{}p{7.75cm}@{}}
		Substitute $t\mapsto -49t+2/27$ \\
		\hline $2\theta^2-t\cdot(459\theta^2 + 459\theta +
		132)-t^2\cdot571536\left(\theta+\frac{7}{6}\right)\left(\theta+\frac{5}{6}\right)$ \\
		\hline
		$a_n=1,66, 78120, 20849556, 16193337300,\dots$  \\
		\hline
		\Fnp{$F=\frac{N}{x^2y^3}$, where}{\includegraphics{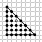}} \\
		$\begin{aligned}
				N={}&1024+768y+528y^2+400y^3+132y^4+48y^5\\
				&+16y^6+768x+480xy+264xy^2+150xy^3\\
				&+33xy^4+6xy^5+3264x^2+1632x^2y\\
				&+333x^2y^2+66x^2y^3-57x^2y^4+1552x^3\\
				&+582x^3y+75x^3y^2-76x^3y^3+3264x^4\\
				&+816x^4y-228x^4y^2+768x^5+96x^5y\\
				&+1024x^6
			\end{aligned}$ \\
		\hline
		$L=0.00416061797\dots={\bf ?}$ \\
		\hline
		$M_0=\smat{ 1 & 0 \\ 1 & 1 }$, $M_{s_2}=\smat{ 5 & -8 \\ 2 & -3 }$, $M_{s_1}=\smat{ 1 & -7 \\ 0 & 1 }$ \end{tabular} & \begin{tabular}{cccc} $I_1$ & $I_2$ & $I_7$ & $II$ \\ $0$ & $-\frac{1}{432}$ &  $\frac{2}{1323}$ & $\infty$ \end{tabular} & \begin{tabular}{c} $2^7$ \\ \\ \hline \\ $T=\smat{ 1 & 0 \\ -\frac{1}{2\pi i}\,{\color{Green}\log\!\left(2^7\right)} & 1 }$\\ \end{tabular} &
	\\
	\hline
	$14\underline{5}II$ & \begin{tabular}{@{}p{7.75cm}@{}}
		$g_2=640t^3 + 140t^2 - 10t + 1/12$ \\
		$g_3=-1024t^5 - 3200t^4 + 440/3t^3 - 110/3t^2$ \\ $+ 5/6t - 1/216$ \\
		\hline
		$3\theta^2-t\cdot(328\theta^2 + 328\theta +
		90)+t^2\cdot432\left(\theta+\frac{7}{6}\right)\left(\theta+\frac{5}{6}\right)$ \\
		\hline
		$a_n=1, 30, 1830, 137580, 11391750, 997411140,\dots$ \\
		\hline
		\Fnp{$F=\frac{9x^2y^2 + 12x^2y + 4x^2 + 16xy^2 + 30xy + 12x + 16y + 9}{xy}$}{\includegraphics{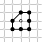}}\\
		\hline
		$L=0.03397980736\dots=\frac{1}{5}\log\!\left(\frac{2^5}{3^3}\right)$\\
		\hline
		 $M_{s_1}=\smat{ 1 & -1 \\ 0 & 1 }$, $M_{s_2}=\smat{ 9 & -4 \\ 16 & -7 }$, $M_0=\smat{ 1 & 0 \\ 5 & 1 }$ \end{tabular} & \begin{tabular}{cccc} $I_1$ & $I_4$ &
		$I_5$ & $II$ \\ $\frac{1}{108}$ & $\frac{3}{4}$ & $0$ & $\infty$
	\end{tabular} & \begin{tabular}{c} \\ $2^{10}\cdot3^4$ \\  \\ \hline \\ $T=\smat{ 1 & 0 \\ -\frac{1}{2\pi i}\,{\color{Green}\log\!\left(2^{10}\cdot 3^{4}\right)} & 5 }$\\ \end{tabular} & \multirow{2}{*}{\makecell[c]{$\Gamma_{10;(1,4,5)}$}}
	\\
	\cline{1-4}
	$1\underline{4}5II$ & \begin{tabular}{@{}p{7.75cm}@{}}
		$g_2=-40960t^3 + 1580t^2 - 20t + 1/12$ \\
		$g_3=4194304t^5 - 450560t^4 + 45640/3t^3 - 695/3t^2$ \\ $+ 5/3t - 1/216$ \\
		\hline $15\theta^2-t\cdot(2576\theta^2 + 2576\theta
		+ 900)$\\$+t^2\cdot110592\left(\theta+\frac{7}{6}\right)\left(\theta+\frac{5}{6}\right)$ \\
		\hline
		$a_n=1,60, 4260, 320880, 24920100, 1972024560,\dots$ \\
		\hline \Fnp{$F=\frac{3x^2y-x^2+22xy^2+60xy+6x-45y^3-51y^2+49y-9}{xy}$}{\includegraphics{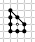}} \\
		\hline
		$L=0.06587495912\dots=\frac{5}{2^7}\log\!\left(\frac{3^3}{5}\right)$\\
		\hline
		$M_{s_1}=\smat{ -1 & -4 \\ 1 & 3 }$, $M_0=\smat{ 1 & 0 \\ 4 & 1 }$, $M_{s_2}=\smat{ 1 & -5 \\ 0 & 1 }$ \end{tabular} & \begin{tabular}{cccc} $I_1$ & $I_4$ & $I_5$ & $II$ \\ $\frac{5}{432}$ &
		$0$ & $\frac{3}{256}$ & $\infty$ \end{tabular} & \begin{tabular}{c} $-3^5\cdot5$ \\ \\ \hline \\ $T=\smat{ 1 & 0 \\ -\frac{1}{2\pi i}\,{\color{Green}\log\!\left(-3^{5}\cdot 5\right)} & 4 }$\\ \end{tabular} & 
	\\
	\cline{1-5}
	$\underline{1}45II$ & \begin{tabular}{@{}p{7.75cm}@{}}
		Substitute $t\mapsto -64t+1/108$ \\
		\hline $5\theta^2-t\cdot(34128\theta^2 + 34128\theta+ 9300)$ \\
		$-t^2\cdot2985984\left(\theta+\frac{7}{6}\right)\left(\theta+\frac{5}{6}\right)$\\
		\hline
		\resizebox{\linewidth}{!}{$a_n=1,1860, 7357860, 35492084880, 188281083032100,\dots$}  \\
		\hline
		\Fnp{$F=\frac{N}{x^2y^3}$, where}{\includegraphics{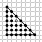}} \\
		$\begin{aligned}
				N={}&-1953125+1218750y-433875y^2+107300y^3\\
				&-17355y^4+1950y^5-125y^6+1875000x\\
				&-705000xy+139440xy^2-17808xy^3+24xy^4\\
				&+120xy^5-834375x^2+174300x^2y-17622x^2y^2\\
				&+1860x^2y^3-3x^2y^4+214000x^3-17520x^3y\\
				&+48x^3y^2+16x^3y^3-33375x^4+30x^4y\\
				&-3x^4y^2+3000x^5+120x^5y-125x^6
			\end{aligned}$ \\
		\hline
		$L=0.00052717285\dots={\bf ?}$ \\
		\hline
		$M_0=\smat{ 1 & 0 \\ 1 & 1 }$, $M_{s_2}=\smat{ 9 & -16 \\ 4 & -7 }$, $M_{s_1}=\smat{ 1 & -5 \\ 0 & 1 }$ \end{tabular} & \begin{tabular}{cccc} $I_1$ & $I_4$ & $I_5$ & $II$ \\ $0$ &
		$-\frac{5}{432}$ & $\frac{1}{6912}$ & $\infty$ \end{tabular} & \begin{tabular}{c} $5^4$ \\ \\ \hline \\ $T=\smat{ 1 & 0 \\ -\frac{1}{2\pi i}\,{\color{Green}\log\!\left(5^4\right)} & 1 }$\\ \end{tabular} & \multirow{1}{*}{\makecell[c]{$\Gamma_{10;(1,4,5)}$}}
	\\
	\hline
	\hline
    \begin{tabular}{c}
	$23\underline{5}II$\\
    ($\I(5)$)
    \end{tabular} & \begin{tabular}{@{}p{7.75cm}@{}}
		$g_2=-180t^3 + 130/3t^2 - 10/3t + 1/12$ \\
		$g_3=432t^5 - 325t^4 + 1855/27t^3$ \\ $- 115/18t^2 + 5/18t - 1/216$ \\
		\hline
		$2\theta^2-t\cdot(59\theta^2 + 59\theta +
		20)+t^2\cdot432\left(\theta+\frac{7}{6}\right)\left(\theta+\frac{5}{6}\right)$ \\
		\hline
		$a_n=1, 10, 120, 1540, 20500, 279480, 3876600,\dots$ \\
		\hline
		\Fnp{$F=\frac{2 - 3x - 2x^2 + 3x^2y + 2x^2y^2 - xy^2 + y + 10xy}{xy}$}{\includegraphics{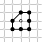}} \\
		\hline
		$L=0.27725887222\dots=\frac{2}{5}\log 2$\\
		\hline
		$M_{s_1}=\smat{ 1 & -2 \\ 0 & 1 }$, $M_{s_2}=\smat{ 4 & -3 \\ 3 & -2 }$, $M_0=\smat{ 1 & 0 \\ 5 & 1 }$ \end{tabular} & \begin{tabular}{cccc} $I_2$ & $I_3$ & $I_5$
		& $II$ \\ $\frac{1}{16}$ & $\frac{2}{27}$ & $0$ & $\infty$ \end{tabular} & \begin{tabular}{c} \\ $2^3$ \\  \\ \hline \\ $T=\smat{ 1 & 0 \\ -\frac{1}{2\pi i}\,{\color{Green}\log\!\left(2^{3}\right)} & 5 }$\\ \end{tabular} & \multirow{1}{*}{\makecell[c]{$\Gamma_{10;(2,3,5)}$}}
	\\
	\cline{1-5}
	\begin{tabular}{c}
	$2\underline{3}5II$\\
    ($\I(3)$)
    \end{tabular} & \begin{tabular}{@{}p{7.75cm}@{}}
		$g_2=14580t^3 + 270t^2 - 10t + 1/12$ \\
		$g_3=-314928t^5 - 120285t^4 + 2835t^3 - 95/2t^2$ \\ $+ 5/6t - 1/216$ \\
		\hline $10\theta^2-t\cdot(999\theta^2 + 999\theta +
		300)$\\$+t^2\cdot11664\left(\theta+\frac{7}{6}\right)\left(\theta+\frac{5}{6}\right)$ \\
		\hline
		$a_n=1,30, 1440, 85260, 5606100, 391231080,\dots$ \\
		\hline
		\Fnp{$F=\frac{2x^3-3x^2y-12xy^2+20y^3-3x^2+30xy-21y^2-12x-21y+20}{xy}$}{\includegraphics{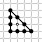}} \\
		\hline
		$L=0.04132287987\dots=\frac{5}{3^3}\log\!\left(\frac{5}{2^2}\right)$\\
		\hline
		$M_{s_1}=\smat{ 1 & -2 \\ 0 & 1 }$, $M_0=\smat{ 1 & 0 \\ 3 & 1 }$, $M_{s_2}=\smat{ 6 & -5 \\ 5 & -4 }$ \end{tabular} & \begin{tabular}{cccc} $I_2$ & $I_3$ & $I_5$ & $II$ \\ $\frac{5}{432}$ &
		$0$ & $\frac{2}{27}$ & $\infty$ \end{tabular} & \begin{tabular}{c} $2^5\cdot5^2$ \\ \\ \hline \\ $T=\smat{ 1 & 0 \\ -\frac{1}{2\pi i}\,{\color{Green}\log\!\left(2^{5}\cdot 5^{2}\right)} & 3 }$\\ \end{tabular} & \multirow{2}{*}{\makecell[c]{$\Gamma_{10;(2,3,5)}$}}
	\\
	\cline{1-4}
	\begin{tabular}{c}
	$\underline{2}35II$\\
    ($\I(2)$)
    \end{tabular} & \begin{tabular}{@{}p{7.75cm}@{}}
		Substitute $t\mapsto t+1/16$ \\
		\hline $5\theta^2-t\cdot(352\theta^2 + 352\theta +
		100)$\\$-t^2\cdot6912\left(\theta+\frac{7}{6}\right)\left(\theta+\frac{5}{6}\right)$ \\
		\hline
		$a_n=1,20, 1140, 68240, 4572100, 321427920,\dots$ \\
		\hline
		\Fnp{$F=\frac{N}{xy^2}$, where}{\includegraphics{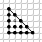}} \\
		$\begin{aligned}
				N={}&5-4y-2y^2-4y^3+5y^4-20x-4xy+20xy^2\\
				&-12xy^3+70x^2+4x^2y+6x^2y^2-100x^3\\
				&-60x^3y+125x^4
			\end{aligned}$ \\
		\hline
		$L=0.03746947912\dots={\bf ?}$ \\
		\hline
		$M_0=\smat{ 1 & 0 \\ 2 & 1 }$, $M_{s_1}=\smat{ 1 & -3 \\ 0 & 1 }$, $M_{s_2}=\smat{ 6 & -5 \\ 5 & -4 }$ \end{tabular} & \begin{tabular}{cccc} $I_2$ & $I_3$ & $I_5$ & $II$ \\ $0$ &
		$\frac{5}{432}$ & $-\frac{1}{16}$ & $\infty$ \end{tabular} & \begin{tabular}{c} $5^3$ \\ \\ \hline \\ $T=\smat{ 1 & 0 \\ -\frac{1}{2\pi i}\,{\color{Green}\log\!\left(5^{3}\right)} & 2 }$\\ \end{tabular} &
	\\
	\hline
	$11\underline{7}III$ & \begin{tabular}{@{}p{7.75cm}@{}}
		$g_2=-4t^3 + 10/3t^2 - 2/3t + 1/12$ \\
		$g_3=-7/3t^4 + 35/27t^3 - 7/18t^2 + 1/18t - 1/216$ \\
		\hline
		$2\theta^2-t\cdot(13\theta^2 + 13\theta + 4)+t^2\cdot64\left(\theta+\frac{5}{4}\right)\left(\theta+\frac{3}{4}\right)$ \\
		\hline
		$a_n=1, 2, 0, -28, -140, -168, 2184, 16200,\dots$ (See \cite{NolanDalbyCondon}) \\
		\hline
		\Fnp{$F=\frac{(x + y)(2x + 1 - xy)}{xy}$}{\includegraphics{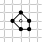}} \\
		\hline
		$L:$ does not exist\\
		\hline
		$M_{s_1}=\smat{ -1 & -1 \\ 4 & 3}$, $M_{s_2}=\smat{ 2 & -1 \\ 1 & 0 }$, $M_0=\smat{ 1 & 0 \\ 7 & 1 }$ \end{tabular} & \begin{tabular}{cccc} $I_1$ & $I_1$ & $I_7$
		& $III$ \\ \multicolumn{2}{l}{{\tiny$\frac{13\pm i\sqrt{7}}{128}$}} & $0$ & $\infty$ \end{tabular} & \begin{tabular}{c} \\ $-2$ \\  \\ \hline \\ $T=\smat{ 1 & 0 \\ -\frac{1}{2\pi i}\,{\color{Green}\log\!\left(-2\right)} & 7 }$\\ \end{tabular} & \makecell[c]{$\Gamma_{9;(1,1,7)}$}
	\\
	\hline
	\hline
	\begin{tabular}{l}$12\underline{6}III$ \\ $\left(\stackrel{\text{isog.}}{\simeq} 2\underline{3}4III\right)$ \end{tabular} & \begin{tabular}{@{}p{7.75cm}@{}}
		$g_2=-16t^3 + 12t^2 - 2t + 1/12$ \\
		$g_3=-16t^4 + 28/3t^3 - 2t^2 + 1/6t - 1/216$ \\
		\hline
		$\theta^2-t\cdot(20\theta^2 + 20\theta + 6)+t^2\cdot64\left(\theta+\frac{5}{4}\right)\left(\theta+\frac{3}{4}\right)$ \\
		\hline
		$a_n=1, 6, 54, 588, 7110, 91476, 1224636,\dots$ \\
		\hline
		\Fnp{$F=\frac{(2x+y+1)(xy+x+2y)}{xy}$}{\includegraphics{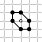}} \\
		\hline
		$L:$ see $2\underline{3}4III$\\
		\hline
		$M_{s_1}=\smat{ 1 & -1 \\ 0 & 1 }$, $M_{s_2}=\smat{ 5 & -2 \\ 8 & -3 }$, $M_0=\smat{ 1 & 0 \\ 6 & 1 }$ \end{tabular} & \begin{tabular}{cccc} $I_1$ & $I_2$ & $I_6$ & $III$ \\
		$\frac{1}{16}$ & $\frac{1}{4}$ & $0$ & $\infty$ \end{tabular} & \begin{tabular}{c} \\ $2^4$ \\  \\ \hline \\ $T=\smat{ 1 & 0 \\ -\frac{1}{2\pi i}\,{\color{Green}\log\!\left(2^{4}\right)} & 6 }$\\ \end{tabular} & $\Gamma_{9;(1,2,6)}$
	\\
	\cline{1-5}
	\begin{tabular}{l}$1\underline{2}6III$ \\ $\left(\stackrel{\text{isog.}}{\simeq}23\underline{4}III\right)$ \end{tabular} & \begin{tabular}{@{}p{7.75cm}@{}}
		Substitute $t\mapsto -4t+1/4$ \\
		\hline
		$3\theta^2 -t\cdot(112\theta^2+112\theta+36)$\\$+t^2\cdot1024\left(\theta+\frac{5}{4}\right)\left(\theta+\frac{3}{4}\right)$ \\
		\hline
		$a_n=1, 12, 180, 2928, 49860, 875952,\dots$ \\
		\hline
		\Fnp{$F=\frac{P\,Q}{x^2y^2}$, where}{\includegraphics{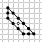}} \\
		$\begin{aligned}
				P={}&x^2-2xy-x-3y^2+y,\\
				Q={}&x^3+x^2y-x^2-5xy^2-2xy+3y^3-y^2
			\end{aligned}$ \\
		\hline
		$L:$ see $23\underline{4}III$\\
		\hline
		$M_{s_1}=\smat{ -1 & -4 \\ 1 & 3 }$, $M_0=\smat{ 1 & 0 \\ 2 & 1 }$, $M_{s_2}=\smat{ 1 & -6 \\ 0 & 1 }$ \end{tabular}& \begin{tabular}{cccc} $I_1$ & $I_2$ & $I_6$ & $III$ \\
		$\frac{3}{64}$ & $0$ & $\frac{1}{16}$ & $\infty$  \end{tabular} & \begin{tabular}{c} $-3$ \\ \\ \hline \\ $T=\smat{ 1 & 0 \\ -\frac{1}{2\pi i}\,{\color{Green}\log\!\left(-3\right)} & 2 }$\\ \end{tabular} & \multirow{2}{*}{$\Gamma_{9;(1,2,6)}$}
	\\*
	\cline{1-4}
	\begin{tabular}{l}$\underline{1}26III$ \\ $\left(\stackrel{\text{isog.}}{\simeq}\underline{2}34III\right)$ \end{tabular} & \begin{tabular}{@{}p{7.75cm}@{}}
		Substitute $t\mapsto -4t+1/16$ \\
		\hline
		$3\theta^2 -t\cdot(128\theta^2+128\theta+36)$\\$-t^2\cdot4096\left(\theta+\frac{5}{4}\right)\left(\theta+\frac{3}{4}\right)$ \\
		\hline
		$a_n=1, 12, 612, 25392, 1298340, 67041072,\dots$ \\
		\hline
		\Fnp{$F=\frac{N}{x^2y^3}$, where}{\includegraphics{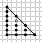}} \\
		$\begin{aligned}
				N={}&-9+6y-7y^2+20y^3-7y^4+6y^5-9y^6\\
				&-27x^2+12x^2y-2x^2y^2+12x^2y^3\\
				&+5x^2y^4-27x^4+6x^4y+5x^4y^2-9x^6
			\end{aligned}$ \\
		\hline
		$L:$ see $\underline{2}34III$ \\
		\hline
		$M_0=\smat{ 1 & 0 \\ 1 & 1 }$, $M_{s_2}=\smat{ 5 & -8 \\ 2 & -3 }$, $M_{s_1}=\smat{ 1 & -6 \\ 0 & 1 }$ \end{tabular}& \begin{tabular}{cccc} $I_1$ & $I_2$ & $I_6$ & $III$ \\
		$0$ & $-\frac{3}{64}$ & $\frac{1}{64}$ & $\infty$  \end{tabular} & \begin{tabular}{c} $3^2$ \\ \\ \hline \\ $T=\smat{ 1 & 0 \\ -\frac{1}{2\pi i}\,{\color{Green}\log\!\left(3^2\right)} & 1 }$\\ \end{tabular} &
	\\
	\hline
	$13\underline{5}III$ & \begin{tabular}{@{}p{7.75cm}@{}}
		$g_2=4t^3 + 6t^2 - 2t + 1/12$ \\
		$g_3=-5t^4 + 5/3t^3 - 3/2t^2 + 1/6t - 1/216$ \\
		\hline
		$6\theta^2-t\cdot(131\theta^2 + 131\theta +
		36)+t^2\cdot64\left(\theta+\frac{5}{4}\right)\left(\theta+\frac{3}{4}\right)$ \\
		\hline
		$a_n=1, 6, 72, 1068, 17460, 301896, 5414136,\dots$ \\
		\hline
		\Fnp{$F=\frac{1 + 2x + x^2 + 4x^2y + 4x^2y^2 + 3xy^2 + 2y + 6xy}{xy}$}{\includegraphics{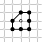}} \\
		\hline
		$L=0.17260924347\dots=\frac{3}{5}\log\!\left(\frac{2^2}{3}\right)$\\
		\hline
		 $M_{s_1}=\smat{ 1 & -1 \\ 0 & 1 }$, $M_{s_2}=\smat{ 7 & -3 \\ 12 & -5 }$, $M_0=\smat{ 1 & 0 \\ 5 & 1 }$ \end{tabular} & \begin{tabular}{cccc} $I_1$ & $I_3$ &
		$I_5$ & $III$ \\ $\frac{3}{64}$ & $2$ & $0$ & $\infty$ \end{tabular} & \begin{tabular}{c} \\ $2^3\cdot3$ \\  \\ \hline \\ $T=\smat{ 1 & 0 \\ -\frac{1}{2\pi i}\,{\color{Green}\log\!\left(2^{3}\cdot 3\right)} & 5 }$\\ \end{tabular} & \multirow{2}{*}{$\Gamma_{9;(1,3,5)}$}
	\\
	\cline{1-4}
	$1\underline{3}5III$ & \begin{tabular}{@{}p{7.75cm}@{}}
		$g_2=-12500t^3 + 750t^2 - 14t + 1/12$ \\
		$g_3=-78125t^4 + 14375/3t^3 - 223/2t^2$ \\ $+ 7/6t - 1/216$ \\
		\hline $2\theta^2-t\cdot(253\theta^2 + 253\theta +
		84)$\\$+t^2\cdot8000\left(\theta+\frac{5}{4}\right)\left(\theta+\frac{3}{4}\right)$ \\
		\hline
		$a_n=1, 42, 2160, 118740, 6750900, 391745592,\dots$ \\
		\hline
		\Fnp{$F=\frac{(x+y-2)(4x^2-17xy+4x+4y^2+4y+1)}{xy}$}{\includegraphics{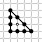}} \\
		\hline
		$L=0.09241962407\dots=\frac{2}{15}\log 2$\\
		\hline
		$M_0=\smat{ 1 & 0 \\ 3 & 1 }$, $M_{s_1}=\smat{ -1 & -4 \\ 1 & 3 }$, $M_{s_2}=\smat{ 1 & -5 \\ 0 & 1 }$ \end{tabular} & \begin{tabular}{cccc} $I_1$ & $I_3$ & $I_5$ & $III$ \\ $\frac{1}{64}$ &
		$0$ & $\frac{2}{125}$ & $\infty$  \end{tabular} & \begin{tabular}{c} $-2^5$ \\ \\ \hline \\ $T=\smat{ 1 & 0 \\ -\frac{1}{2\pi i}\,{\color{Green}\log\!\left(-2^{5}\right)} & 3 }$\\ \end{tabular} & 
	\\
	\cline{1-5}
	$\underline{1}35III$ & \begin{tabular}{@{}p{7.75cm}@{}}
		Substitute $t\mapsto -125t+3/64$ \\
		\hline $3\theta^2-t\cdot(7808\theta^2 + 7808\theta + 2124)$\\$-t^2\cdot512000\left(\theta+\frac{5}{4}\right)\left(\theta+\frac{3}{4}\right)$ \\
		\hline
		$a_n=1, 708, 1086660, 2023823760, 4143653100900,\dots$ \\
		\hline
		\Fnp{$F=\frac{N}{x^2y^3}$, where}{\includegraphics{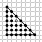}} \\
		$\begin{aligned}
				N={}&14348907-2243862y-161595y^2+58860y^3\\
				&-1995y^4-342y^5+27y^6+3188646x\\
				&-135594xy-14724xy^2+1676xy^3-34xy^4\\
				&+30xy^5+767637x^2-38556x^2y-1794x^2y^2\\
				&+708x^2y^3+5x^2y^4+84564x^3+1836x^3y\\
				&-4x^3y^2+4x^3y^3+9477x^4+18x^4y\\
				&+5x^4y^2+486x^5+30x^5y+27x^6
			\end{aligned}$ \\
		\hline
		$L=0.00136176272\dots={\bf ?}$ \\
		\hline
		$M_0=\smat{ 1 & 0 \\ 1 & 1 }$, $M_{s_2}=\smat{ 7 & -12 \\ 3 & -5 }$, $M_{s_1}=\smat{ 1 & -5 \\ 0 & 1 }$ \end{tabular} & \begin{tabular}{cccc} $I_1$ & $I_3$ & $I_5$ & $III$ \\ $0$ &
		$-\frac{1}{64}$ & $\frac{3}{8000}$ & $\infty$ \end{tabular} & \begin{tabular}{c} $3^5$ \\ \\ \hline \\ $T=\smat{ 1 & 0 \\ -\frac{1}{2\pi i}\,{\color{Green}\log\!\left(3^5\right)} & 1 }$\\ \end{tabular} & \multirow{1}{*}{$\Gamma_{9;(1,3,5)}$}
	\\
	\hline
	\hline
	\begin{tabular}{l} $23\underline{4}III$ ($\O(4)$) \\ $\left(\stackrel{\text{isog.}}{\simeq}1\underline{2}6III\right)$ \end{tabular} & \begin{tabular}{@{}p{7.75cm}@{}}
		$g_2=1024t^3 + 144t^2 + 11/4t + 1/768$ \\
		$g_3=-4096t^4 - 1024/3t^3 - 10t^2 - 1/24t + 1/110592$ \\
		\hline
		$3\theta^2-t\cdot(112\theta^2 + 112\theta +
		36)$\\$+t^2\cdot1024\left(\theta+\frac{5}{4}\right)\left(\theta+\frac{3}{4}\right)$ \\
		\hline
		$a_n=1, 12, 180, 2928, 49860, 875952,\dots$ \\
		\hline
		\Fnp{$F=\frac{(xy-x+3y+1)(xy+3x-y+1)}{xy}$}{\includegraphics{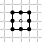}} \\
		\hline
		$L=0.20598980413\dots=\frac{3}{2^4}\log 3$\\
		\hline
		$M_{s_1}=\smat{ 1 & -2 \\ 0 & 1 }$, $M_{s_2}=\smat{ 4 & -3 \\ 3 & -2 }$, $M_0=\smat{ 1 & 0 \\ 4 & 1 }$ \end{tabular} & \begin{tabular}{cccc} $I_2$ & $I_3$ & $I_4$
		& $III$ \\ $\frac{3}{64}$ & $\frac{1}{16}$ & $0$ & $\infty$  \end{tabular} & \begin{tabular}{c} \\ $3^2$ \\  \\ \hline \\ $T=\smat{ 1 & 0 \\ -\frac{1}{2\pi i}\,{\color{Green}\log\!\left(3^{2}\right)} & 4 }$\\ \end{tabular} & \multirow{1}{*}{$\Gamma_{9;(2,3,4)}$}
	\\
	\cline{1-5}
	\begin{tabular}{l}$2\underline{3}4III$ ($\O(3)$) \\ $\left(\stackrel{\text{isog.}}{\simeq}12\underline{6}III\right)$ \end{tabular} & \begin{tabular}{@{}p{7.75cm}@{}} $g_2=64t^3 + 12t^2 - 2t + 1/12$ \\
		$g_3=-128t^4 + 56/3t^3 - 2t^2 + 1/6t - 1/216$ \\
		\hline
		$\theta^2-t\cdot(20\theta^2 + 20\theta + 6)+t^2\cdot64\left(\theta+\frac{5}{4}\right)\left(\theta+\frac{3}{4}\right)$ \\
		\hline
		$a_n=1, 6, 54, 588, 7110, 91476, 1224636,\dots$ \\
		\hline
		\Fnp{$F=\frac{(x + y - 1)(2x^2 -4xy + x + 2y^2 + y - 1)}{xy}$}{\includegraphics{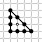}} \\
		\hline
		$L=0.23104906019\dots=\frac{1}{3}\log 2$\\
		\hline
		$M_{s_1}=\smat{ 1 & -2 \\ 0 & 1 }$, $M_0=\smat{ 1 & 0 \\ 3 & 1 }$, $M_{s_2}=\smat{ 5 & -4 \\ 4 & -3 }$ \end{tabular} &
	\begin{tabular}{cccc} $I_2$ & $I_3$ & $I_4$ & $III$ \\ $\frac{1}{16}$ &
		$0$ & $\frac{1}{4}$ & $\infty$ \end{tabular} & \begin{tabular}{c} \\ $2^2$ \\  \\ \hline \\ $T=\smat{ 1 & 0 \\ -\frac{1}{2\pi i}\,{\color{Green}\log\!\left(2^2\right)} & 3 }$\\ \end{tabular} & \multirow{2}{*}{$\Gamma_{9;(2,3,4)}$}
	\\*
	\cline{1-4}
	\begin{tabular}{l}$\underline{2}34III$ ($\O(2)$) \\ $\left(\stackrel{\text{isog.}}{\simeq}\underline{1}26III\right)$ \end{tabular} & \begin{tabular}{@{}p{7.75cm}@{}}
		Substitute $t\mapsto t+3/64$ \\
		\hline $3\theta^2-t\cdot(128\theta^2 + 128\theta +
		36)$\\$-t^2\cdot4096\left(\theta+\frac{5}{4}\right)\left(\theta+\frac{3}{4}\right)$ \\
		\hline
		$a_n=1, 12, 612, 25392, 1298340, 67041072,\dots$ \\
		\hline
		\Fnp{$F=-\frac{P\,Q}{x^2y^2}$, where}{\includegraphics{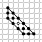}} \\
		$\begin{aligned}
				P={}&3x^2-2xy-x+3y^2-y,\\
				Q={}&9x^3-5x^2y-3x^2-5xy^2+2xy+9y^3-3y^2
			\end{aligned}$ \\
		\hline
		$L=0.04872722313\dots={\bf ?}$\\
		\hline
		$M_0=\smat{ 1 & 0 \\ 2 & 1 }$, $M_{s_1}=\smat{ 1 & -3 \\ 0 & 1 }$, $M_{s_2}=\smat{ 5 & -4 \\ 4 & -3 }$ \end{tabular} & \begin{tabular}{cccc} $I_2$ & $I_3$ & $I_4$ & $III$ \\ $0$ &
		$\frac{1}{64}$ & $-\frac{3}{64}$ & $\infty$ \end{tabular} & \begin{tabular}{c} $3^4$ \\ \\ \hline \\ $T=\smat{ 1 & 0 \\ -\frac{1}{2\pi i}\,{\color{Green}\log\!\left(3^4\right)} & 2 }$\\ \end{tabular} &
	\\
	\hline
	\begin{tabular}{l}$11\underline{6}IV$ \\ $\left(\stackrel{\text{isog.}}{\simeq}\underline{2}33IV\right)$ \end{tabular} & \begin{tabular}{@{}p{7.75cm}@{}}
		$g_2=-6t^2 + 1/12$ \\
		$g_3=-9t^4 + 1/2t^2 - 1/216$ \\
		\hline
		$\theta^2-t^2\cdot81\left(\theta+\frac{4}{3}\right)\left(\theta+\frac{2}{3}\right)$ \\
		\hline
		$a_n=1, 0, 18, 0, 810, 0, 45360, 0, 2806650,\dots$ \\
		\hline
		\Fnp{$F=\frac{(1+x)(1 + x + 3xy^2)}{xy}$}{\includegraphics{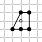}} \\
		\hline
		$L:$ does not exist (HG pullback) \\
		\hline
		$M_{s_1}=\smat{ 1 & -1 \\ 0 & 1 }$, $M_{s_2}=\smat{ 4 & -1 \\ 9 & -2 }$, $M_0=\smat{ 1 & 0 \\ 6 & 1 }$ \end{tabular} & \begin{tabular}{cccc} $I_1$ & $I_1$ & $I_6$ & $IV$ \\ $\frac{1}{9}$ &
		$-\frac{1}{9}$ & $0$ & $\infty$ \end{tabular} & \begin{tabular}{c} \\ $3^3$ \\  \\ \hline \\ $T=\smat{ 1 & 0 \\ -\frac{1}{2\pi i}\,{\color{Green}\log\!\left(3^{3}\right)} & 6 }$\\ \end{tabular} & \multirow{2}{*}{$\Gamma_0(3)$}
	\\
	\cline{1-4}
	\begin{tabular}{l}$1\underline{1}6IV$ \\ $\left(\stackrel{\text{isog.}}{\simeq}23\underline{3}IV\right)$ \end{tabular} & \begin{tabular}{@{}p{7.75cm}@{}}
		Substitute $t\mapsto -3t+1/9$ \\
		\hline
		$2\theta^2 -t\cdot(81\theta^2+81\theta+24)+ t^2\cdot729\left(\theta+\frac{4}{3}\right)\left(\theta+\frac{2}{3}\right)$ \\
		\hline
		$a_n=1, 12, 198, 3720, 75690, 1626912,\dots$ \\
		\hline
		\Fnp{$F_2=\frac{N}{x^2y}$, where}{\includegraphics{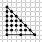}} \\
		$\begin{aligned}
				N={}&-1-6y+3y^2+52y^3+12y^4-96y^5-64y^6\\
				&+3xy+21xy^2+18xy^3-96xy^4-96xy^5\\
				&+12x^2y+9x^2y^2-72x^2y^3-96x^2y^4\\
				&-2x^3+3x^3y-24x^3y^2-56x^3y^3-6x^4y\\
				&-24x^4y^2-6x^5y-x^6
			\end{aligned}$ \\
		\hline
		$L:$ see $23\underline{3}IV$ \\
		\hline
		$M_{s_1}=\smat{ -17 & -81 \\ 4 & 19 }$, $M_0=\smat{ 1 & 0 \\ 1 & 1 }$, $M_{s_2}=\smat{ 1 & -6 \\ 0 & 1 }$ \end{tabular}&
	\begin{tabular}{cccc} $I_1$ & $I_1$ & $I_6$ & $IV$ \\ $\frac{2}{27}$ &
		$0$ & $\frac{1}{27}$ & $\infty$ \end{tabular} & \begin{tabular}{c} $2$ \\ \\ \hline \\ $T=\smat{ 1 & 0 \\ -\frac{1}{2\pi i}\,{\color{Green}\log\!\left(2\right)} & 1 }$\\ \end{tabular} & \\*
	\cline{1-4}
	\begin{tabular}{l}$\underline{1}16IV$ \\ $\left(\stackrel{\text{isog.}}{\simeq}2\underline{3}3IV\right)$ \end{tabular} & \begin{tabular}{@{}p{7.75cm}@{}}
		Same as for $1\underline{1}6IV$  \end{tabular}&
	\begin{tabular}{cccc} $I_1$ & $I_1$ & $I_6$ & $IV$ \\ $0$ &
		$\frac{2}{27}$ & $\frac{1}{27}$ & $\infty$ \end{tabular} & \begin{tabular}{c} Same as for $1\underline{1}6IV$  \end{tabular} & \\
	\hline

	$12\underline{5}IV$ & \begin{tabular}{@{}p{7.75cm}@{}}
		$g_2=10/3t^2 - 4/3t + 1/12$ \\
		$g_3=-t^4 + 28/27t^3 - 13/18t^2 + 1/9t - 1/216$ \\
		\hline
		$2\theta^2-t\cdot(29\theta^2 + 29\theta + 8)+t^2\cdot27\left(\theta+\frac{4}{3}\right)\left(\theta+\frac{2}{3}\right)$ \\
		\hline
		$a_n=1, 4, 30, 280, 2890, 31584, 358176, \dots$ \\
		\hline
		\Fnp{$F=\frac{(1+x)(1 + x + 2 y + 2 x y + x y^2)}{xy}$}{\includegraphics{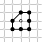}} \\
		\hline
		$L=0.27725887222\dots=\frac{2}{5}\log 2$\\
		\hline
		$M_{s_1}=\smat{ 1 & -1 \\ 0 & 1 }$, $M_{s_2}=\smat{ 5 & -2 \\ 8 & -3 }$, $M_0=\smat{ 1 & 0 \\ 5 & 1 }$ \end{tabular} & \begin{tabular}{cccc} $I_1$ & $I_2$ & $I_5$
		& $IV$ \\ $\frac{2}{27}$ & $1$ & $0$ & $\infty$ \end{tabular} & \begin{tabular}{c} \\ $2$ \\  \\ \hline \\ $T=\smat{ 1 & 0 \\ -\frac{1}{2\pi i}\,{\color{Green}\log\!\left(2\right)} & 5 }$\\ \end{tabular} & \multirow{2}{*}{$\mathrm{PSL}_2(\Z)$}
	\\
	\cline{1-4}
	$1\underline{2}5IV$ & \begin{tabular}{@{}p{7.75cm}@{}}
		Substitute $t\mapsto -25t+1$ \\
		\hline $\theta^2-t\cdot(52\theta^2 + 52\theta +
		16)+t^2\cdot675\left(\theta+\frac{4}{3}\right)\left(\theta+\frac{2}{3}\right)$ \\
		\hline
		$a_n=1, 16, 330, 7360, 170650, 4050816,\dots$ \\
		\hline
		\Fnp{$F=\frac{N}{xy^2}$, where}{\includegraphics{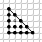}} \\
		$\begin{aligned}
				N={}&1+y-4y^2+y^3+y^4+2x+4xy+16xy^2\\
				&+3xy^3-x^2+2x^2y-6x^2y^2-2x^3-3x^3y\\
				&+x^4
			\end{aligned}$ \\
		\hline
		$L=0.21520447048\dots=\frac{1}{\sqrt{5}}\log\!\left(\frac{1+\sqrt{5}}{2}\right)=\frac{1}{2}L(1,\chi_5)$\\
		\hline
		$M_{s_1}=\smat{ -1 & -4 \\ 1 & 3 }$, $M_0=\smat{ 1 & 0 \\ 2 & 1 }$, $M_{s_2}=\smat{ 1 & -5 \\ 0 & 1 }$ \end{tabular} & \begin{tabular}{cccc} $I_1$ & $I_2$ & $I_5$ & $IV$ \\ $\frac{1}{27}$ & $0$
		& $\frac{1}{25}$ & $\infty$ \end{tabular} & \begin{tabular}{c} $-1$ \\ \\ \hline \\ $T=\smat{ 1 & 0 \\ -\frac{1}{2\pi i}\,{\color{Green}\log\!\left(-1\right)} & 2 }$\\ \end{tabular} & 
	\\
	\cline{1-5}
	$\underline{1}25IV$ & \begin{tabular}{@{}p{7.75cm}@{}}
		Substitute $t\mapsto -25t+2/27$ \\
		\hline $2\theta^2-t\cdot(621\theta^2 + 621\theta +
		168)$\\$-t^2\cdot18225\left(\theta+\frac{4}{3}\right)\left(\theta+\frac{2}{3}\right)$ \\
		\hline
		\resizebox{\linewidth}{!}{$a_n=1, 84, 16830, 3971640, 1030948650, 282553385184,\dots$}  \\
		\hline
		\Fnp{$F=\frac{N}{x^2y^3}$, where}{\includegraphics{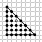}} \\
		$\begin{aligned}
				N={}&-1024-768y+240y^2+400y^3+60y^4-48y^5\\
				&-16y^6+768x+96xy-72xy^2-6xy^3-3xy^4\\
				&-18xy^5-960x^2-288x^2y+27x^2y^2+84x^2y^3\\
				&-9x^2y^4+400x^3-66x^3y+6x^3y^2-14x^3y^3\\
				&-240x^4-12x^4y-9x^4y^2+48x^5-18x^5y\\
				&-16x^6
			\end{aligned}$ \\
		\hline
		$L=0.01058254145\dots={\bf ?}$ \\
		\hline
		$M_0=\smat{ 1 & 0 \\ 1 & 1 }$, $M_{s_2}=\smat{ 5 & -8 \\ 2 & -3 }$, $M_{s_1}=\smat{ 1 & -5 \\ 0 & 1 }$ \end{tabular} & \begin{tabular}{cccc} $I_1$ & $I_2$ & $I_5$ & $IV$ \\ $0$ &
		$-\frac{1}{27}$ & $\frac{2}{675}$ & $\infty$ \end{tabular} & \begin{tabular}{c} \\ $2^5$ \\ \\ \hline \\ $T=\smat{ 1 & 0 \\ -\frac{1}{2\pi i}\,{\color{Green}\log\!\left(2^5\right)} & 1 }$\\ \end{tabular} & \multirow{1}{*}{$\mathrm{PSL}_2(\Z)$}
	\\
	\hline
	\hline
	\begin{tabular}{l}$23\underline{3}IV$ ($\T(3)$) \\ $\left(\stackrel{\text{isog.}}{\simeq}1\underline{1}6IV\right)$ \end{tabular} & \begin{tabular}{@{}p{7.75cm}@{}}
		$g_2=54t^2 - 4t + 1/12$ \\
		$g_3=-729t^4 + 108t^3 - 17/2t^2 + 1/3t - 1/216$ \\
		\hline $2\theta^2-t\cdot(81\theta^2 +
		81\theta + 24)+t^2\cdot729\left(\theta+\frac{4}{3}\right)\left(\theta+\frac{2}{3}\right)$ \\
		\hline
		$a_n=1, 12, 198, 3720, 75690, 1626912,\dots$ \\
		\hline
		\Fnp{$F_1=\frac{(x -2y - 1)(x + y + 2)(y - 2x - 1)}{xy}$}{\includegraphics{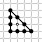}} \\
		
		\hline
		$L=0.15403270679\dots=\frac{2}{3^2}\log 2$\\
		\hline
		$M_{s_1}=\smat{ 1 & -2 \\ 0 & 1 }$, $M_{s_2}=\smat{ 4 & -3 \\ 3 & -2 }$, $M_0=\smat{ 1 & 0 \\ 3 & 1 }$ \end{tabular} & \begin{tabular}{cccc} $I_2$ & $I_3$ & $I_3$ & $IV$ \\ $\frac{1}{27}$ &
		$\frac{2}{27}$ & $0$ & $\infty$ \end{tabular} & \begin{tabular}{c} \\ $2^3$ \\ \\ \hline \\ $T=\smat{ 1 & 0 \\ -\frac{1}{2\pi i}\,{\color{Green}\log\!\left(2^3\right)} & 3 }$\\ \end{tabular} & \multirow{2}{*}{$\Gamma_0(3)$}
	\\
	\cline{1-4}
	\begin{tabular}{l}$2\underline{3}3IV$ \\ $\left(\stackrel{\text{isog.}}{\simeq}\underline{1}16IV\right)$ \end{tabular} & \begin{tabular}{@{}p{7.75cm}@{}}
		Same as for $23\underline{3}IV$  \end{tabular} &
	\begin{tabular}{cccc} $I_2$ & $I_3$ & $I_3$ & $IV$ \\ $\frac{1}{27}$ &
		$0$ & $\frac{2}{27}$ & $\infty$ \end{tabular} & \begin{tabular}{c}  Same as for $23\underline{3}IV$  \\ \end{tabular} &
	\\*
	\cline{1-5}
	\begin{tabular}{l}$\underline{2}33IV$ ($\T(2)$) \\ $\left(\stackrel{\text{isog.}}{\simeq}11\underline{6}IV\right)$ \end{tabular} & \begin{tabular}{@{}p{7.75cm}@{}}
		Substitute $t\mapsto -(1/3)t+1/27$ \\
		\hline $\theta^2-t^2\cdot81\left(\theta+\frac{4}{3}\right)\left(\theta+\frac{2}{3}\right)$ \\
		\hline
		$a_n=1, 0, 18, 0, 810, 0, 45360, 0, 2806650,\dots$ \\
		\hline
		\Fnp{$F_2=\frac{9x^4-3x^2y^2+3x^2y+6x^2+y^4+y^3+y+1}{xy^2}$}{\includegraphics{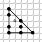}} \\
		\hline
		$L:$ does not exist (HG pullback)\\
		\hline
		$M_0=\smat{ 1 & 0 \\ 2 & 1 }$, $M_{s_1}=\smat{ 1 & -3 \\ 0 & 1 }$, $M_{s_2}=\smat{ 4 & -3 \\ 3 & -2 }$ \end{tabular} &
	\begin{tabular}{cccc} $I_2$ & $I_3$ & $I_3$ & $IV$ \\ $0$ &
		$-\frac{1}{9}$ & $\frac{1}{9}$ & $\infty$ \end{tabular} & \begin{tabular}{c} $3$ \\ \\  \hline \\ $T=\smat{ 1 & 0 \\ -\frac{1}{2\pi i}\,{\color{Green}\log\!\left(3\right)} & 2 }$\\ \end{tabular} & \multirow{1}{*}{$\Gamma_0(3)$}
	\\
	\hline
\end{longtable}
\endgroup

\restoregeometry

\end{document}